\documentclass[a4paper]{article}

\pdfoutput=1

\usepackage{epic}
\usepackage{eepic}
\usepackage{amsmath}
\usepackage{amssymb}
\usepackage{theorem}
\usepackage{enumerate}
\usepackage{tikz-cd}
\usepackage{mathtools}

\let\OLDthebibliography\thebibliography
\renewcommand\thebibliography[1]{
	\OLDthebibliography{#1}
	\setlength{\itemsep}{0pt}
}

\newcommand{\al}{\alpha}
\newcommand{\bt}{\beta}

\newcommand{\ov}{\overline}
\newcommand{\la}{\langle}
\newcommand{\ra}{\rangle}
\newcommand{\C}{\mathbb C}

\newcommand{\Q}{\mathbb Q}
\newcommand{\R}{\mathbb R}
\newcommand{\Z}{\mathbb Z}
\newcommand{\co}{\cal O}
\newcommand{\g}{\mathfrak g}
\newcommand{\h}{\mathfrak h}

\newcommand{\Aut}{\operatorname{Aut}}

\newcommand{\mult}{\operatorname{mult}}
\newcommand{\sign}{\operatorname{sign}}

\newcommand{\tr}{\operatorname{tr}}
\newcommand{\rk}{\operatorname{rk}}

\newcommand{\q}{\operatorname{q}}
\newcommand{\im}{\operatorname{im}}
\newcommand{\af}{\operatorname{af}}
\newcommand{\ch}{\operatorname{ch}}
\newcommand{\vac}{\mathbf 1}

\newcommand{\eop}{\hspace*{\fill} $\Box$}

\theoremstyle{break}
\newtheorem{thm}{Theorem}[section]
\newtheorem{prp}[thm]{Proposition}

\begin{document}

\begin{center}
{\Large\bf Reflective modular varieties and their cusps}\\[12mm]
Thomas Driscoll-Spittler, Nils R.\ Scheithauer, Janik Wilhelm\\
Fachbereich Mathematik, Technische Universit\"at Darmstadt\\
Darmstadt, Deutschland
\end{center}

\vspace*{18mm}
\noindent
We classify reflective automorphic products of singular weight
under certain regularity assumptions.
Using obstruction theory we show that there are exactly 11 such functions. They are naturally related to certain conjugacy classes in Conway's group $\text{Co}_0$. The corresponding modular varieties have a very rich geometry. We establish a bijection between their $1$-dimensional type-$0$ cusps and the root systems in Schellekens' list.
We also describe a $1$-dimensional cusp along which the restriction of the automorphic product is given by the eta product of the corresponding class in $\text{Co}_0$. Finally we apply our results to give a complex-geometric proof of Schellekens' list.

\vspace*{14mm}
\begin{tabular}{ll}
1  & Introduction \\
2  & Modular forms for the Weil representation\\
3  & Bounds for reflective modular forms\\
4  & Automorphic forms on orthogonal groups\\
5  & Reflective forms of singular weight\\
6  & Cusps of orthogonal modular varieties\\
7  & Holomorphic vertex operator algebras of central charge $24$\\
\multicolumn{2}{l}{Appendix}
\end{tabular}
\vspace*{14mm}

\section{Introduction}

Automorphic forms for the orthogonal group $\text{O}_{n,2}(\R)$ are functions on the hermitian symmetric space
$\text{SO}_{n,2}(\R)^+/( \text{SO}_n(\R) \times \text{SO}_2(\R))$
which are up to a cocycle invariant under an integral subgroup of $\text{O}_{n,2}(\R)$.  They naturally generalise the classical modular forms on $\text{SL}_2(\R)$. Such a form is called reflective if its divisor is generated
by rational quadratic divisors corresponding to hyperplanes of reflections in the integral subgroup.
Reflective automorphic forms have many applications, for example in the representation theory of infinite-dimensional Lie algebras, in the study of moduli spaces and in differential geometry (see e.g.\ \cite{B92, B98, GN1, GN2, M21, S06, S17}, \cite{B96, GHS07, GHS10} and \cite{Y04, Y13}, a nice overview is given in \cite{G18}). As we will see, they can also be used to classify holomorphic vertex operator algebras of central charge $24$.

Borcherds' singular theta correspondence \cite{B98} is a map from modular forms for the Weil representation of $\text{SL}_2(\Z)$ to automorphic forms on orthogonal groups. These automorphic forms have nice product expansions at the $0$-dimensional cusps and therefore are called automorphic products. The simplest example is Dedekind's eta function
which is the lift of the theta function of the $A_1$-lattice.
The divisor of an automorphic product is a linear combination of rational quadratic divisors. Bruinier \cite{Br02, Br14} showed that under some weak assumptions an automorphic form
with such a divisor is an automorphic product.

The reflective automorphic forms which occur in Lie theory usually have zeros of order $1$. We include this in our definition of reflectivity. With regard to our application to vertex operator algebras we also assume that the roots of length $2$ contribute to the divisor. (For the precise definition see Section 5.) Ma \cite{M18} showed that up to scalings
there are only finitely many lattices carrying such forms.
Our first main result gives explicit values for the possible levels (see Theorem \ref{thebeatles} and the following table).

{\em Let $L$ be a regular even lattice of signature $(n,2)$, $n>2$ and even, carrying a reflective automorphic product. Then the level of $L$ is one of the following:}
\[
\renewcommand{\arraystretch}{1.2}
\begin{array}{r|l}
n  & \text{level}\\[0.5mm] \hline 
   &   \\[-4mm]
4  & 3, 4, 6, 7, 8, 9, 11, 12, 14, 15, 16, 18, 19, 20, 21, 22, 23, 24, 27, 28\\
   & 30, 33, 35, 36, 40, 42, 44, 45, 49, 60, 63, 72, 75, 98, 100, 121 \\
6  & 2, 3, 4, 5, 6, 7, 8, 9, 10, 11, 12, 14, 15, 18, 20, 25, 36       \\
8  & 3, 4, 6, 7, 8, 9, 12                \\
10 & 1, 2, 3, 4, 5, 6, 9 \\
12 & 3, 4  \\
14 & 2, 3, 4 \\     
18 & 1, 2 \\
26 & 1
\end{array}            
\]

The theorem generalises results of Scheithauer \cite{S17} and Dittmann \cite{D} to arbitrary levels. The idea of the proof is as follows. Let $F$ be a vector-valued modular form which lifts to a reflective automorphic product. The regularity condition implies that
$L'/L$ contains sufficiently many non-trivial isotropic elements.
This can be used to construct a linear combination of the components of $F$ which has small pole orders at the cusps. The Riemann-Roch theorem then implies the restrictions on the levels.

The smallest possible weight of a non-constant holomorphic automorphic form on $\text{O}_{n,2}(\R)$ is $(n-2)/2$, the so-called singular weight. Forms of this weight are of particular interest. Their Fourier expansions at $0$-dimensional cusps are supported only on isotropic vectors.
This implies that there are no non-trivial cusp forms of singular weight. Our second main result shows that reflective automorphic products of singular weight are very rare (see Theorems \ref{cabedelo} and \ref{blazefoley}).

{\em There are exactly 11 regular even lattices $L$ of signature $(n,2)$, $n>2$ and even, splitting $I\!I_{1,1} \oplus I\!I_{1,1}$ which carry a reflective automorphic product
  of singular weight. They are given in the following table:}
\[
\renewcommand{\arraystretch}{1.2}
\begin{array}{r|l}
n  & L \\[0.5mm] \hline 
   &   \\[-4mm]
  26 & I\!I_{26,2} \\[0.7mm]
  18 & I\!I_{18,2}(2_{I\!I}^{+10}) \\[0.7mm]
  14 & I\!I_{14,2}(2_{I\!I}^{-10} 4_{I\!I}^{-2}), \, I\!I_{14,2}(3^{-8}) \\[0.7mm]
  12 & I\!I_{12,2}(2_2^{+2} 4_{I\!I}^{+6}) \\[0.7mm]
  10 & I\!I_{10,2}(2_{I\!I}^{+6} 3^{-6}), \, I\!I_{10,2}(5^{+6}) \\[0.7mm]
  8  & I\!I_{8,2}(2_{I\!I}^{+4} 4_{I\!I}^{-2} 3^{+5}), \, I\!I_{8,2}(2_7^{+1} 4_7^{+1} 8_{I\!I}^{+4}), \, I\!I_{8,2}(7^{-5}) \\[0.7mm]
 6  & I\!I_{6,2}(2_{I\!I}^{-2} 4_{I\!I}^{-2} 5^{+4}) 
\end{array}
\]
{\em In each case the corresponding automorphic product is unique up to $\text{O}(L)^+$ and corresponds naturally to a unique conjugacy class in $\text{Co}_0$.}

The proof is based on obstruction theory and explicit construction. We describe this in more detail. Using the first main result we get a list of $474$ lattices which meet the assumptions and potentially carry a reflective automorphic product. For each of them we check whether it can satisfy the Eisenstein condition for singular weight. Then we are left with $132$ lattices. We construct obstructions coming from cusp forms to eliminate another
$121$ lattices. The remaining $11$ lattices are in natural correspondence with certain conjugacy classes in $\text{Co}_0$, the orthogonal group of the Leech lattice $\Lambda$. This allows us to construct explicitly for each case a vector-valued modular form $F$ which lifts to a reflective automorphic product $\psi_F$ of singular weight. The uniqueness follows again from obstruction theory.  

Note that after splitting $I\!I_{1,1} \oplus I\!I_{1,1}$ we obtain exactly the $11$ genera found by H\"ohn \cite{H17} in his investigation of the genus of the Moonshine module.

Let $L$ be an even lattice of signature $(n,2)$, $n>2$ and $\Gamma \subset \text{O}(L)$ a subgroup containing the discriminant kernel of $L$. Then the orthogonal modular variety $\Gamma^+\backslash {\cal H}$ can be compactified by adding $0$- and $1$-dimensional cusps. We associate to such a cusp two invariants, the type and the associated lattice, and define the notion of a splitting cusp. We show (see Theorem \ref{classcuspsgen}):

{\em The splitting cusps of a given type and associated lattice are parametrised by the double quotient $\ov{\Gamma} \backslash \text{O}(D) / \ov{\text{O}(L)_S^+}$ where $S \subset L$ is any isotropic submodule with
these invariants.}

Special cases of this result were described by Attwell-Duval \cite{At14, At16} and Kiefer \cite{Ki}.

Now we consider the
modular varieties $\text{O}(L,F)^+\backslash {\cal H}$ corresponding to the $11$ re\-flec\-tive automorphic products $\psi_F$ constructed above. To a $1$-dimensional cusp ${\cal C}$ of $\text{O}(L,F)^+\backslash {\cal H}$ we associate a set $R_{\cal C}$ which together with the type of ${\cal C}$ determines the first term in the expansion of $\psi_F$ at ${\cal C}$. For a cusp of type $0$ the set $R_{\cal C}$ is either empty or a scaled root system. In the latter case $\psi_F$ vanishes along ${\cal C}$. Applying the above parametrisation of the splitting cusps we show (see Theorem \ref{heaven17}):

{\em The scaled root systems corresponding to the $1$-dimensional cusps of type $0$ in the reflective modular varieties $\text{O}(L,F)^+\backslash {\cal H}$ are precisely
  those described by Schellekens in his classification of meromorphic conformal field theories of central charge $24$.}

The root systems determine the lowest order approximation of $\psi_F$ at the corresponding cusp (cf.\ Theorem \ref{onedimexp}):

{\em The lowest order term in the expansion of $\psi_F$ at a $1$-dimensional cusp ${\cal C}$ of type $0$ with root system $R_{\cal C}$ is the denominator function of the affine Kac-Moody algebra corresponding to $R_{\cal C}$.}

Since $\psi_F$ has singular weight, the modular variety $\text{O}(L,F)^+\backslash {\cal H}$ has a $1$-dimensional cusp on which $\psi_F$
is not identically zero.
For each of the $11$ cases we locate such a cusp (see Theorem \ref{falco}).

{\em Let $g \in \text{Co}_0$ be in the class corresponding to $\psi_F$, $m$ the order of $g$ and $\eta_g$ the eta product associated with $g$. Then there exists a cusp ${\cal C}$ of $\text{O}(L,F)^+\backslash {\cal H}$ such that the restriction of $\psi_F$ to ${\cal C}$ is $\eta_g$. The cusp ${\cal C}$ has type a cyclic subgroup of $D^{N/m}$ of order $m$ and associated lattice $\Lambda^g$. If $N/m =1$, i.e.\ in $8$ out of the $11$ cases, ${\cal C}$ is the unique splitting cusp with these invariants. In the remaining cases ${\cal C}$ does not split.} 

A pictorial description of the situation for $I\!I_{12,2}(2_2^{+2} 4_{I\!I}^{+6})$ is given at the end of Section \ref{chemicalbrothers}.

  In \cite{G12} Gritsenko considered the automorphic product $\psi_F$ on the unimodular lattice $I\!I_{26,2}$ and determined its expansions at the $24$ $1$-dimensional cusps (which are necessarily of type $0$).
  This paper was one of the main motivations for our investigations.

  Finally we describe the relation
  of the above results to vertex algebras.
Let $V$ be a holomorphic vertex operator algebra of central charge $24$. Then the weight-$1$ subspace $V_1$ is a Lie algebra. This Lie algebra is either trivial, abelian of rank $24$ or non-trivial and semisimple. In the second case $V$ is the vertex algebra associated to the Leech lattice $\Lambda$. We consider the third case. Then $V_1$ can be written as
\[  V_1 = \g_{1,k_1} \oplus \ldots \oplus \g_{m,k_m}  \]
with simple factors $\g_i$ and scalings $k_i \in \Z_{>0}$. This decomposition is called the affine structure of $V$. The subalgebra $\la V_1 \ra$ of $V$ generated by $V_1$ is isomorphic to the tensor product
\[ \la V_1 \ra \simeq L_{\g_1,k_1} \otimes \ldots \otimes L_{\g_m,k_m} \]
where $L_{\g_i,k_i}$ is the simple affine vertex algebra associated with $\g_i$ of level $k_i$. The vertex algebra $V$ decomposes into finitely many modules over $\la V_1 \ra$. Its character
\[    \chi_V = \tr_V e^{2\pi i v_0} q^{L_0-1}  \]
is a Jacobi form of weight $0$ and lattice index $M = \bigoplus_{i=1}^m Q_i^{\vee}(k_i)$ where $Q_i^{\vee}(k_i)$ is the coroot lattice $Q_i^{\vee}$ of $\g_i$ with the bilinear form rescaled by $k_i$. Adding the cominimal simple currents in $V$ to $M$ we obtain the lattice $K$ associated with $V$. Then (see Theorem \ref{radiohead}):

{\em We can decompose the character $\chi_V$ as 
  \[  \chi_V = \sum_{\gamma \in K'/K}  F_{\gamma} \theta_{\gamma}  \]
  where $F = \sum_{\gamma \in K'/K}  F_{\gamma} e^{\gamma}$ is a modular form for the Weil representation of the lattice $K$ and
  $\theta = \sum_{\gamma \in K'/K}  \theta_{\gamma} e^{\gamma}$ the Jacobi theta function of $K$.}

The proof is based on the representation theory of affine Kac-Moody algebras. (Strictly speaking, here and in the following statements, we would have to assume for technical reasons that $K$ has even rank.)

We associate a Lie algebra $\g(V)$ to $V$ describing the physical states of a bosonic string moving on a torus orbifold.
If $V$ is unitary, the Lie algebra $\g(V)$ is a generalised Kac-Moody algebra. From this we derive (see Theorem \ref{everythingcounts}):

{\em The modular form $F$ defines a reflective automorphic product $\psi_F$ of singular weight on $L = K \oplus I\!I_{1,1} \oplus I\!I_{1,1}$.}

We can recover $\g(V)$ and the affine structure of $V$ from $\psi_F$. The modular variety $\text{O}(L,F)^+\backslash {\cal H}$ has a unique $0$-dimensional cusp ${\cal C}$ of type $0$. The expansion of $\psi_F$ at ${\cal C}$ is the denominator function of $\g(V)$. For the affine structure of $V$ we find (see Theorem \ref{aw}):

{\em The decomposition $L = K \oplus I\!I_{1,1} \oplus I\!I_{1,1}$ defines
  a $1$-dimensional cusp ${\cal C}$ of $\text{O}(L,F)^+\backslash {\cal H}$ of type $0$ with associated lattice $K$. The scaled root system $R_{\cal C}$ associated with ${\cal C}$ is the root system of the affine structure of $V$ together with its scaling.}

Combining this result with the above classification of type-$0$ cusps we obtain:

{\em Let $V$ be a holomorphic vertex operator algebra of central charge $24$ with non-trivial, semisimple weight-$1$ space. Suppose $V$ is unitary and the lattice associated with $V$ is regular and of even rank. Then the affine structure of $V$ is one of the $69$ non-trivial structures given in Theorem \ref{heaven17}}.

\medskip

Our results give a natural explanation of the $11$ classes found by H\"ohn in his investigation of the genus of the Moonshine module \cite{H17} and a complex-geometric proof of Schellekens' list \cite{ANS} under the stated conditions. It is surprising to us that the only reflective automorphic products of singular weight are those coming from holomorphic vertex operator algebras of central charge $24$ and that Schellekens' list accounts for all type-$0$ cusps of the corresponding reflective modular varieties.

\medskip

The paper is organised as follows. In Section 2 we
recall some results on
modular forms for the Weil representation and define reflective modular forms. Then we show that the Riemann-Roch theorem imposes strong restrictions on the level and the weight of a weakly reflective modular form associated with a regular discriminant form. In Section 4 we review automorphic forms on orthogonal groups. In particular we describe Borcherds' singular theta correspondence and Kudla's product expansion at $1$-dimensional cusps. Then we define reflective automorphic products and show that there are exactly $11$ reflective automorphic products of singular weight under certain regularity conditions. In Section 6 we classify the $1$-dimensional cusps of type $0$ of the corresponding modular varieties. We show that they are in natural bijection with the root systems in Schellekens' list. We also describe a $1$-dimensional cusp on which the reflective automorphic product restricts to the eta product of the
associated
class in Conway's group. Finally we relate our classification results to the theory of vertex operator algebras. In particular we obtain a complex-geometric proof of Schellekens' classification of affine structures  of holomorphic vertex operator algebras of central charge $24$. In the Appendix we give upper bounds for the weights of reflective automorphic products and list the cusp forms we used to construct obstructions in Section 5.

 \medskip

 The authors thank Tomoyuki Arakawa, Richard Borcherds, Jan Bruinier, Eberhard Freitag, Jens Funke,
 Va\-le\-ry Gritsenko, Gerald H\"ohn, Victor Kac, Martin M\"oller, Sven M\"oller, Manuel M\"uller, Riccardo Salvati-Manni and Don Zagier
for valuable discussions and helpful comments.

The authors acknowledge support by the LOEWE research unit \emph{Uniformized Structures in Arithmetic and Geometry} and by the DFG through the CRC \emph{Geometry and Arithmetic of Uniformized Structures}, project number 444845124.

\section{Modular forms for the Weil representation} \label{modforms}

We recall some results on modular forms for the Weil representation and define reflective modular forms.

\medskip

A discriminant form is a finite abelian group $D$ with a quadratic form $\q : D \to \Q/\Z$ such that $(\bt,\gamma) = \q(\bt + \gamma)- \q(\bt) - \q(\gamma) \! \mod 1$ is a non-degenerate symmetric bilinear form. The level of $D$ is the smallest positive integer $N$ such that $N\q(\gamma) = 0 \! \mod 1$ for all $\gamma \in D$. If $L$ is an even lattice, then $L'/L$ is a discriminant form with the quadratic form given by $\q(\gamma) = \gamma^2/2 \! \mod 1$. Conversely every discriminant form can be obtained in this way. The corresponding lattice can be chosen to be positive-definite. The signature $\sign(D) \in \Z/8\Z$ of a discriminant form $D$ is defined as the signature modulo $8$ of any even lattice with that discriminant form.

Let $c$ be an integer. Then $c$ acts by multiplication on $D$ and we have an exact sequence
$0 \to D_c \to D \to D^c \to 0$
where $D_c$ is the kernel and $D^c$ the image of this map. Note that $D^c$ is the orthogonal complement of $D_c$.
The set $D^{c*} = \{ \gamma \in D \, | \, c\q(\alpha) + (\alpha,\gamma) = 0 \, \text{ for all } \alpha \in D_c \}$ is a coset of $D^c$. After a choice of Jordan decomposition there is a canonical coset representative $x_c \in D^{c*}$ with $2x_c=0$. We can write $\gamma \in D^{c*}$ as $\gamma = x_c + c \mu$. Then the reduced norm $\q_c(\gamma) = c \q(\mu) + x_c \mu \! \mod 1$ is well-defined. We will also use the notations $D_{c,x} = \{ \gamma \in D_c \, | \, \q(\gamma) = x \! \mod 1 \}$, $D_{c,x,n} = \{ \gamma \in D_{c,x} \, | \, \text{order}(\gamma) = n \}$ and $D^{c*}_x = \{ \gamma \in D^{c*} \, | \, \q_c(\gamma) = x \! \mod 1 \}$.

A discriminant form $D$ of level $N$ is called {\em regular} if $D^{N/p}$ contains a non-trivial isotropic element for each prime $p|N$.

For more details on discriminant forms we refer the reader to \cite{AGM}, \cite{B00}, \cite{CS}, \cite{N} and \cite{S09}. 

Let $L$ be an even lattice. We define the level of $L$ as the level of its discriminant form $L'/L$ and we say that $L$ is regular if $L'/L$ is regular.

\medskip

Let $D$ be a discriminant form of even signature. We define a scalar product on the group ring $\C [D]$ which is linear in the first and antilinear in the second variable by $(e^{\gamma},e^{\bt})= \delta_{\gamma \bt}$. Then there is a unitary action of the group $\text{SL}_2(\Z)$ on $\C [D]$ defined by 
\begin{align*} 
\rho_D(T) e^{\gamma}  & = e(-\q(\gamma))\, e^{\gamma} \\
\rho_D(S) e^{\gamma}  & = \frac{e(\sign(D)/8)}{\sqrt{|D|}}
                  \sum_{\bt\in D} e((\gamma,\bt))\, e^{\bt}   
\end{align*}
where 
$S = \left( \begin{smallmatrix} 0 & -1 \\ 1 & 0 \end{smallmatrix} \right)$ and 
$T = \left( \begin{smallmatrix} 1 &  1 \\ 0 & 1 \end{smallmatrix} \right)$
are the standard generators of $\text{SL}_2(\Z)$. This rep\-re\-sen\-tation is called the Weil representation of $\text{SL}_2(\Z)$ on $\C[D]$. It commutes with the orthogonal group $\text{O}(D)$ of $D$. For a general matrix $M = \left( \begin{smallmatrix} a & b \\ c & d \end{smallmatrix} \right) \in \text{SL}_2(\Z)$ we have
\[ \rho_D(M) e^{\gamma} = \xi \frac{\sqrt{|D_c|}}{\sqrt{|D|}}  
                    \sum_{\bt \in D^{c*}} 
                    e(-a\q_c(\bt)) \, e(-b(\bt,\gamma)) \, e(-bd \q(\gamma)) \,  
                    e^{d\gamma + \bt}
\]
with $\xi = e(\sign(D)/4) \prod \xi_p$. The local factors $\xi_p$ can be expressed in terms of the Jordan components of $D$
(see Theorem 4.7 in \cite{S09}).

\medskip

Let 
\[ F(\tau) = \sum_{\gamma \in D} F_{\gamma}(\tau) e^{\gamma}  \]
be a holomorphic function on the complex upper halfplane
with values in $\C [D]$ and $k$ an integer. Then $F$ is a modular form for $\rho_D$ of weight $k$ if 
\[ F(M\tau) = (c\tau + d)^k \rho_D(M) F(\tau)  \]
for all $M = \left( \begin{smallmatrix} a & b \\ c & d \end{smallmatrix} \right) \in \text{SL}_2(\Z)$ and $F$ is meromorphic at $\infty$. We say that $F$ is symmetric if it is invariant under $\text{O}(D)$. 

We can easily construct modular forms for the Weil representation by symmetrising scalar-valued modular forms on congruence subgroups.

Let $D$ be a discriminant form of even signature and level dividing $N$. Let $f$ be a scalar-valued modular form on $\Gamma_0(N)$ of weight $k$ and character $\chi_D$ and $H$ an isotropic subset of $D$ which is invariant under $(\Z/N\Z)^*$ as a set. Then
\[ F_{\, \Gamma_0(N),f,H} = \sum_{M \in \Gamma_0(N)\setminus \text{SL}_2(\Z) } \, 
             \sum_{\gamma \in H} 
             f|_{k,M}  \, \rho_D(M^{-1}) e^\gamma  \]
is a modular form for $\rho_D$ of weight $k$. Here $|_k$ denotes the Petersson slash operator of weight $k$.
If $\gamma \in D$ and $f$ is a scalar-valued modular form on $\Gamma_1(N)$ of weight $k$ and character $\chi_{\gamma}$, then
\[ F_{\, \Gamma_1(N),f,\gamma} = \sum_{M \in \Gamma_1(N)\setminus \text{SL}_2(\Z) } \, 
             f|_{k,M}  \, \rho_D(M^{-1}) e^{\gamma}  \]
is a modular form for $\rho_D$ of weight $k$. 
An analogous result holds for $\Gamma(N)$. Every modular form for $\rho_D$  can be written as a linear combination of liftings from $\Gamma_1(N)$ or $\Gamma(N)$.
Explicit formulas for the components of the symmetrisations are given in \cite{S09}, \cite{S15}.

We observe that for
\[  F_{\, \Gamma_0(N),f,0} = \sum_{\gamma \in D} F_{\gamma} e^{\gamma}  \]
the components $F_{\gamma}$ only depend on the invariants $\{ (c,\q_c(\gamma)) \, | \, c|N, \gamma \in D^{c*} \}$. 

Among other things we will use this result to construct the Eisenstein series for the dual Weil representation. Let $D$ be a discriminant form of even signature and level dividing $N$. Then 
\[ E_k = \frac{1}{2} \sum_{M \in \Gamma_{\infty}^+  \setminus \Gamma_1(N)} 1|_{k,M} \]
with $\Gamma_{\infty}^+ = \{ \, T^n \, | \, n \in \Z \, \}$ is an Eisenstein series for $\Gamma_1(N)$ of weight $k$. Let $\gamma \in D$ be isotropic. Then 
\[ E_{\gamma} = 
  \sum_{M \in \Gamma_1(N)\setminus \text{SL}_2(\Z) } \, E_k|_{k,M}  \, \ov{\rho}_D(M^{-1}) e^{\gamma} \]
is an Eisenstein series for the dual Weil representation $\ov{\rho}_D$. For an equivalent definition see \cite{Br02}. In the case $\gamma = 0$ we can write 
\[ E_0 = \sum_{M \in \Gamma_0(N)\setminus \text{SL}_2(\Z) } \, E_{k,\chi}|_{k,M}  \, \ov{\rho}_D(M^{-1}) e^0 \]
where
\[ E_{k,\chi} 
= \sum_{M \in \Gamma_1(N) \setminus \Gamma_0(N)} \chi(M) \, E_k|_{k,M}    \]
is an Eisenstein series for $\Gamma_0(N)$ of weight $k$ and character $\ov{\chi} = \chi = \chi_D$. We denote $E = E_0$.
There are explicit formulas for the Fourier coefficients of $E$ (cf.\ \cite{BrK}, \cite{S06} and \cite{KY}) and Opitz \cite{O} wrote a program which computes these numbers.

The dimension of the space of holomorphic modular forms for the Weil representation can be worked out using the Riemann-Roch theorem \cite{F} or the Selberg trace formula \cite{ES}, \cite{B00}. 

Let $D$ be a discriminant form of even signature. If $F$ is a modular form for $\rho_D$ of weight $2-k$ and $G$ a modular form for $\ov{\rho}_D$ of weight $k$ then the product $(F,\ov{G}) = \sum_{\gamma \in D} F_\gamma G_{\gamma}$ is a modular form for $\text{SL}_2(\Z)$ of weight $2$ with a possible pole at $\infty$. By the residue theorem the constant coefficient of $(F,\ov{G})$ has to vanish. 

\medskip

Finally we define reflective modular forms. Let $D$ be a discriminant form of level $N$
and even signature. A modular form $F = \sum_{\gamma \in D} F_{\gamma} e^{\gamma}$ for $\rho_D$ of weight $k$ is called {\em weakly reflective} if the following conditions are satisfied:

\begin{enumerate}[i)]
\item There exists an even lattice $L$ of signature $(n,2)$ with $n>2$ and $k=1-n/2$ such that $D=L'/L$. 
\item If $F_{\gamma}$ has a pole at $\infty$, then there is a divisor $d|N$ such that
$\gamma \in D_{d,1/d}$ and 
  \[ F_{\gamma}(\tau) = c_{\gamma} q^{-1/d} + \ldots  \]
  for some $c_{\gamma}\in \Q$.
\end{enumerate}

We call $F$ {\em reflective} if in addition $F_0(\tau) = q^{-1} + \ldots $ and $F_{\gamma}(\tau) = q^{-1/d} + \ldots $ if $F_{\gamma}$ is singular. For the mo\-ti\-vation of these definitions we refer to Section 5.

\section{Bounds for reflective modular forms}

Let $D$ be a discriminant form of even signature and level $N$ and $F : H \to \C[D]$ a weakly reflective modular form of weight $k$ for the Weil representation. If $D$ is regular, we can construct a modular form $g$ for $\Gamma_0(N)$ which has small pole orders at the cusps by adding the components of $F$ appropriately. Then the Riemann-Roch theorem implies strong restrictions on $N$ and $k$.

\medskip

Let $D$ be a discriminant form of even signature and level $N=p^{\nu}$, $p$ prime and $v = \sum_{\bt \in D} v_{\bt} e^{\bt}$ some element in $\C[D]$ which is supported only on $I^{N/p} = I \cap D^{N/p}$ where $I$ is the set of isotropic elements in $D$.

\begin{prp}  \label{ironmaiden}
Let $M = \left( \begin{smallmatrix} a & b \\ c & d \end{smallmatrix} \right) \in \text{SL}_2(\Z)$ such that $\nu_p(c)=0$ and $\gamma \in D_{N/p}$. Then 
  \[ (\rho_D(M)e^{\gamma},v) = d_{\gamma} \sum_{\bt \in D} \ov{v_{\bt}}  \]
for some $d_{\gamma} \in \C$.
\end{prp}
{\em Proof:} By the explicit formula for the Weil representation we have
\[
\rho_D(M) e^{\gamma} = \xi \frac{\sqrt{|D_c|}}{\sqrt{|D|}}  
                    \sum_{\mu \in D^c} 
                    e(-a\q_c(\mu)) \, e(-b(\mu,\gamma)) \, e(-bd \q(\gamma)) \,  
                    e^{d\gamma + \mu}  \, .
                  \]
Let $c^{-1}$ be the inverse of $c$ modulo $p^{\nu}$ and $\bt \in I^{N/p}$. Then for $\mu = \bt - d\gamma \in D^c$
\begin{align*}
-a\q_c(\mu) -b(\mu,\gamma) - bd \q(\gamma)
&= -a c^{-1} \q(\bt) + c^{-1} (\bt,\gamma) -c^{-1}d \q(\gamma)  \mod 1 \\
&= c^{-1} (\bt,\gamma) -c^{-1}d \q(\gamma)  \mod 1 \, . 
\end{align*}
The orthogonal complement of $D_{N/p}$ is $D^{N/p}$. Hence
\[ -a\q_c(\mu) -b(\mu,\gamma) - bd \q(\gamma)  = - c^{-1}d \q(\gamma)  \mod 1 \, \]
and
\[   (\rho_D(M)e^{\gamma},e^{\bt})
  = \xi \frac{\sqrt{|D_c|}}{\sqrt{|D|}} e(- c^{-1}d \q(\gamma))  \, . \]
Since the right hand side is independent of $\beta$ this implies the statement. \eop
\medskip

\noindent
Now we consider the inner products $(\rho_D(M)e^{\gamma},v)$ in the case where the valuation $\nu_p(c)$ is positive.

\begin{prp}  \label{slater}
  Suppose $p$ is odd. Let $M = \left( \begin{smallmatrix} a & b \\ c & d \end{smallmatrix} \right) \in \text{SL}_2(\Z)$ such that $\nu_p(c) = j$ with $1 \leq j \leq \nu -2$ and $\gamma \in D_{p^k}$ with
  $\q(\gamma) = x/p^k \bmod 1$ where $\nu_p(x)=0$. Then
\[ (\rho_D(M)e^{\gamma},v) = d_{\gamma} \sum_{\bt \in D} \ov{v_{\bt}}  \]
for some $d_{\gamma} \in \C$ with $d_{\gamma} = 0$ if $\gamma \neq 0$.
\end{prp}  
{\em Proof:} As above
\[
\rho_D(M) e^{\gamma} = \xi \frac{\sqrt{|D_c|}}{\sqrt{|D|}}  
                    \sum_{\mu \in D^c} 
                    e(-a\q_c(\mu)) \, e(-b(\mu,\gamma)) \, e(-bd \q(\gamma)) \,  
                    e^{d\gamma + \mu}  \, . 
\]
Let $\bt \in I^{N/p}$. If $d\gamma + \mu = \bt$ or equivalently $d \gamma = \bt - \mu$ with $\mu \in D^c$, then $\gamma \in D^c$ because $d$ is invertible modulo $p$. Then the condition on the norm of $\gamma$ can only hold if $\gamma = 0$. Under this assumption
\[   (\rho_D(M)e^{\gamma},e^{\bt})
    = \xi \frac{\sqrt{|D_c|}}{\sqrt{|D|}} e(- a \q_c(\bt))  \, . \]
  Write $\bt = p^{\nu-1} \al$. Then $\q_{p^j}(\bt) = p^{2\nu-2-j} \q(\al) = 0 \! \mod 1$. This implies $\q_c(\bt) = 0 \! \mod 1$. \eop

\medskip
\noindent
For $p=2$ the situation is slightly more complicated because of the odd 2-adic Jordan components.

\begin{prp}  \label{slater2}
  Suppose $p=2$. Let $M = \left( \begin{smallmatrix} a & b \\ c & d \end{smallmatrix} \right) \in \text{SL}_2(\Z)$ such that $\nu_p(c) =j$ with $1 \leq j \leq \nu -2$ and $\gamma \in D_{p^k}$ with
  $\q(\gamma) = x/p^k  \! \mod 1$ where $\nu_p(x) = 0$. Then except in the case when $D$ is of type
\[  2_{l_2}^{\epsilon_2 n_2} 4_{l_4}^{\epsilon_4 n_4} 8_{I\!I}^{\epsilon_8 n_8}  \]
with $n_2 = n_4 = 1 \! \mod 2$ and $k = 2$ we have
\[ (\rho_D(M)e^{\gamma},v) = d_{\gamma} \sum_{\bt \in D} \ov{v_{\bt}}  \]
for some $d_{\gamma} \in \C$ with $d_{\gamma} = 0$ if $\gamma \notin D^{c*}$.
\end{prp}
{\em Proof:} As above
\[
\rho_D(M) e^{\gamma} = \xi \frac{\sqrt{|D_c|}}{\sqrt{|D|}}  
                    \sum_{\mu \in D^{c*}} 
                    e(-a\q_c(\mu)) \, e(-b(\mu,\gamma)) \, e(-bd \q(\gamma)) \,  
                    e^{d\gamma + \mu} \, .
                  \]
                  We consider the coefficient at $e^{\bt}$ for $\bt \in I^{N/2}$. It vanishes if $\gamma \notin D^{c*}$. Suppose $\gamma \in D^{c*}$. Then $\mu = \bt -d\gamma \in D^{c*}$ and
\begin{align*}
-a\q_c(\mu) -b(\mu,\gamma) -bd \q(\gamma)
  &= -a \q_c( \mu ) - b(\bt,\gamma) + bd \q(\gamma) \mod 1 \\
  &= -a \q_c( \mu ) + bd \q(\gamma) \mod 1  \, . 
\end{align*}
Now write $\mu = x_c + c \al$. Then $d\gamma = \bt - x_c -c \al$ so that
\[ d^2 \q(\gamma) = \q(x_c) + c^2 \q(\al)  \mod 1 \]
and
\[ d^2 \q_c(\gamma) = \q_c(\mu) - (\bt,\al) \mod 1 \, . \]
We show $\al \in D_{2^{\nu -1}}$. The first equation implies
\[  2^{\nu -j-1} d^2\q(\gamma) = 2^{\nu -j-1} \q(x_c)  \mod 1  \, . \]
Hence $2^{\nu -j-1}\q(\gamma) = 0 \! \mod 1$ except in the excluded case. It follows $k \leq {\nu -j-1}$ and
$2^{\nu -j-1}c\al  = 2^{\nu -j-1} (\bt - x_c - d \gamma) = 0$. This proves the claim.
We get $\q_c(\mu) = d^2 \q_c(\gamma)$ and
\[  (\rho_D(M) e^{\gamma} , e^{\bt}) = \xi \frac{\sqrt{|D_c|}}{\sqrt{|D|}} e(-ad^2\q_c(\gamma) + bd \q(\gamma))  \, . \]
This proves the proposition. \eop
  
\medskip
\noindent
Finally we consider the case $\nu_p(c)\geq \nu$.

\begin{prp}  \label{iggypop}
Let $M = \left( \begin{smallmatrix} a & b \\ c & d \end{smallmatrix} \right) \in \text{SL}_2(\Z)$ such that $\nu_p(c)\geq\nu$ and $\gamma \in D$. Then
\[ (\rho_D(M)e^{\gamma},v) =
\begin{cases}
      \, \chi_D(a) \ov{v_{d\gamma}} & \text{if $\gamma \in I^{N/p}$,} \\
      \, 0                         & \text{otherwise.}
\end{cases}
\]
\end{prp}

\medskip

The group $\Gamma_0(N)$ has index $N \prod_{p|N} (1 + 1/p)$ in $\text{SL}_2(\Z)$ and $\sum_{c|N}\phi((c,N/c))$ classes of cusps. Let $s = a/c \in \Q$ with $(a,c)=1$. Then the equivalence class of $s$ as a cusp of $\Gamma_0(N)$ is determined by the invariants $(c,N)$ (a divisor of $N$) and $ac/(c,N)$ (a unit in $\Z/(c,N/(c,N))\Z\,$). The width of $s$ is $t_s = N/(c^2,N)$. A complete set of representatives of the cusps of $\Gamma_0(N)$ is given by the numbers $a/c \in \Q$, $(a,c)=1$ where $c$ ranges over the divisors of $N$ and $a$ over the units in $\Z/(c,N/c)\Z$.

Let $f : H \to \C$ be a modular form of weight $k$ for $\Gamma_0(N)$ with quadratic Dirichlet character $\chi$ modulo $N$. We assume that $f$ is holomorphic on $H$ and possibly has poles at cusps. We say that $f$ has {\em small pole orders at cusps} if for each cusp $s=a/c \in \Q$ 
\[
  f_s(\tau) = f|_{k, M_s}(\tau) = c_s q^{-1/m_s t_s}  + \ldots  \]
with
\[ m_s =
\begin{cases}
  \, 1  & \text{if $\chi(T_s)=1$,}\\
  \, 2  & \text{if $\chi(T_s)=-1$}
\end{cases}
\]
where $M_s=\left( \begin{smallmatrix} a & b \\ c & d \end{smallmatrix} \right) \in \text{SL}_2(\Z)$ is any matrix sending $\infty$ to $s$, $t_s$ is the width of $s$ and $T_s = M_s T^{t_s} M_s^{-1}$ is a generator of the stabiliser of $s$ in $\Gamma_0(N)$ (cf.\ Section 5 in \cite{S09}).

\begin{thm} \label{ledzeppelin}
Let $D$ be a regular discriminant form of even signature and level $N = \prod_{p|N}p^{\nu_p}$ and $F = \sum_{\gamma \in D} F_{\gamma} e^{\gamma}$ a weakly reflective modular form of weight $k$ for $\rho_D$ with
singular $F_0$. Then there is a linear combination $g$ of the components of $F$ with the following properties:
\begin{enumerate}[i)]
\item The function $g$ is a non-zero modular form of weight $k$ for $\Gamma_0(N)$ with character $\chi_D$ and small pole orders at cusps.
\item If $p|N$ and $p$ odd and $g$ has a pole at the cusp $s=a/c$, then $\nu_p(c) \in \{ 0, \nu_p-1, \nu_p \}$. 
\item If $2|N$ and the $2$-adic part of $D$ is not of type $2_{l_2}^{\epsilon_2 n_2} 4_{l_4}^{\epsilon_4 n_4} 8_{I\!I}^{\epsilon_8 n_8}$ with $n_2 = n_4 = 1 \! \mod 2$ and $g$ has a pole at the cusp $s=a/c$, then $\nu_2(c) \in \{ 0, \nu_2-1, \nu_2 \}$.
\end{enumerate}  
\end{thm}
{\em Proof:} We decompose
\[  D = \bigoplus_{p|N} D_{p^{\nu_p}}  \]
and $\rho_D = \otimes_{p|N} \, \rho_{D_{p^{\nu_p}}}$.
Let $v = \otimes_{p|N} v_p$ with $v_p \in \C[D_{p^{\nu_p}}]$. 
We assume that 
\[ v_p = \sum_{\mu \in D_{p^{\nu_p}}} v_{p,\mu} e^{\mu}  \]
is supported only on $I^{p^{\nu_p-1}}$, invariant under $(\Z/p^{\nu_p}\Z)^*$ and satisfies
\[  \sum_{\mu \in D_{p^{\nu_p}}} v_{p, \mu} = 0 \]
with $v_{p, 0} \neq 0$.
The existence of such an element is ensured by the condition that $D$ is regular.
We define
\[ g = (F,v)
= \sum_{\gamma \in D} F_{\gamma} \ov{v_{\gamma}}  
= \sum_{\gamma \in D} F_{\gamma} \, \prod_{p|N} \, (e^{\gamma_p},v_p)
\]
where $\gamma_p$ denotes the projection of $\gamma$ on $D_{p^{\nu_p}}$. For $M =\left( \begin{smallmatrix} a & b \\ c & d \end{smallmatrix} \right) \in \text{SL}_2(\Z)$ we have
\[  g|_{k,M} = \sum_{\gamma \in D} F_{\gamma} \, \prod_{p|N} \, (\rho_{D_{p^{\nu_p}}}(M)e^{\gamma_p},v_p) \, . 
\]
By Proposition \ref{iggypop} the function $g$ is a modular form for $\Gamma_0(N)$ with character $\chi_D$.
The assumption on $F_0$ ensures that $g$ is non-zero.

First we want to prove the stated properties of the valuations $\nu_p(c)$ if $g$ has a pole at the cusp $s=M\infty=a/c$.
Let $p|N$ be an odd prime such that $1 \leq \nu_p(c) \leq \nu_p-2$. Since $F$ is weakly reflective we find that for every $F_{\gamma}$ with a pole the projection $\gamma_p$ satisfies the condition from Proposition \ref{slater}. This implies that no such $F_\gamma$ can contribute a pole to the expansion of $g|_{k,M}$ because the corresponding factor vanishes by the assumptions on $v$. In the same way we argue for $p=2$ using Proposition \ref{slater2}.

Finally we show that a pole at the cusp $s=M\infty=a/c$ has pole order bounded by $1/t_s$. Suppose $F_{\gamma}$ contributes a pole to $g|_{k,M}$. Then there is an integer $j|N$ such that $j\gamma =0$, $q(\gamma)=1/j \! \mod 1$ and $F_\gamma$ has a pole of order $1/j$.
We have to show that $1/j \leq 1/t_s$ or equivalently $\nu_p(t_s) \leq \nu_p(j)$ for all primes $p|N$. We start with the case $\nu_2\ne 3$. Since $\nu_p(c) \in \{ 0, \nu_p-1, \nu_p \}$, we have $\nu_p(t_s)=\nu_p$ if $\nu_p(c)=0$ and $\nu_p(t_s)=0$ in the other cases. If $\nu_p(c)=0$, Proposition \ref{ironmaiden} implies that the local factor corresponding to $F_{\gamma}$ vanishes unless $\nu_p(j)=\nu_p$. We are left with the case $\nu_2= 3$. For the odd primes $p|N$ and for $p=2$ whenever $\nu_2(c)\neq 1$ we still have $\nu_p(t_s) = 0$ or $\nu_p$. So we can again argue as above. For $\nu_2(c)=1$ we have $\nu_2(t_s)=1$ and by Proposition \ref{slater2} we might get contributions of poles from $F_{\gamma}$
if $\nu_2(j)=2$.
But then $\nu_2(t_s) \leq \nu_2(j)$. \eop

\medskip
\noindent
For a positive integer $N = \prod_{p|N}p^{\nu_p}$ we introduce the local factors
\[
\mu_p(N)=
\begin{cases}
  \, p^{\nu_p-1}  & \text{if $\nu_p\geq 2$},\\
  \, p+1         & \text{if $\nu_p =1$},\\
  \, 1           & \text{if $\nu_p =0$}.
\end{cases}
\]

\begin{thm} \label{thebeatles}
Let $D$ be a regular discriminant form of even signature carrying a weakly reflective modular form $F$ of weight $k$ with singular $0$-component $F_0$. Then the level $N$ of $D$ and the weight $k$ satisfy the inequality
\begin{equation*}
	\frac{-k}{12}\prod_{p\vert N} \mu_p(N)\leq 2^{\omega(N)}
\end{equation*}
where $\omega(N)$ denotes the number of primes that are Hall divisors of $N$.
\end{thm}
{\em Proof:} Write $N = \prod_{p|N}p^{\nu_p}$. Using Theorem \ref{ledzeppelin} we can construct a non-zero modular form $g$ with the stated properties. Define $M_p = \{ 0, \nu_p-1 , \nu_p\}$ except when the $2$-adic part of $D$ is of type $2_{l_2}^{\epsilon_2 n_2} 4_{l_4}^{\epsilon_4 n_4} 8_{I\!I}^{\epsilon_8 n_8}$ with $n_2 = n_4 = 1 \! \mod 2$ in which case we put $M_2 = \{ 0,1,2,3 \}$. Applying the valence formula (see Theorem 4.1 in \cite{HBJ94}) to $g$ yields
\[  \sum_{s \in \Gamma_0(N) \backslash P} t_s \text{ord}_s(g) \leq \frac{k}{12} [ \text{SL}_2(\Z) : \Gamma_0(N) ] \]
where $P = \Q \cup \{ \infty \}$. Proposition 5.1 in \cite{S09} shows that the numbers $m_s$ for a cusp $s=a/c$
depend only on $c$ and are multiplicative. It follows
\[  \frac{-k}{12} \, N \, \prod_{p|N} \left(1+\frac{1}{p}\right)
  \leq \sum_{\substack{ c|N \\ \nu_p(c) \in M_p}} \phi((c,N/c))/m_c  \, . 
\]
We can write the sum on the right-hand side as
\[ \sum_{\substack{ c|N \\ \nu_p(c) \in M_p}} \phi((c,N/c))/m_c
= \, \prod_{p|N} \, \sum_{j \in M_p} \phi\big(p^{\text{min}(\nu_p-j,j)} \big)/m_{p^j} \, . \]
For $p|N$ we have
\[
\sum_{ j \in M_p} \phi\big( p^{\text{min}(\nu_p-j,j)} \big)/m_{p^j} \leq 
\begin{cases}
  \, p+1  & \text{if $\nu_p\geq 2$},\\
  \, 2    & \text{if $\nu_p =1$}.
\end{cases}
\]
To see this we bound $m_s$ from below by $1$ whenever the $2$-adic component of $D$ is not of type 
$2_{l_2}^{\epsilon_2 n_2} 4_{l_4}^{\epsilon_4 n_4} 8_{I\!I}^{\epsilon_8 n_8}$ with $n_1=n_2=1 \!\mod 2$.
If $D$ is of that type, then we have $m_s=2$ for all cusps $s=a/c$ with $\nu_2(c)=1$ or $2$ and so the inequality still holds. The statement of the theorem is now a direct consequence. \eop

\medskip

\noindent
The theorem generalises the results in \cite{S17} and \cite{D} to arbitrary levels. The bounds slightly improve the ones given in \cite{DS}. We also remark that the condition that $F_0$ is singular can be removed. Then the bounds get much larger (cf.\ \cite{DS}).

\medskip

We can easily determine the numbers $N$ which solve the inequality given in Theorem \ref{thebeatles}. When we remove the cases which do not correspond to existing lattices, we obtain the following table:

\[
\renewcommand{\arraystretch}{1.2}
\begin{array}{r|l}
n  & \text{level}\\[0.5mm] \hline 
   &   \\[-4mm]
4  & 3, 4, 6, 7, 8, 9, 11, 12, 14, 15, 16, 18, 19, 20, 21, 22, 23, 24, 27, 28\\
   & 30, 33, 35, 36, 40, 42, 44, 45, 49, 60, 63, 72, 75, 98, 100, 121 \\
6  & 2, 3, 4, 5, 6, 7, 8, 9, 10, 11, 12, 14, 15, 18, 20, 25, 36       \\
8  & 3, 4, 6, 7, 8, 9, 12                \\
10 & 1, 2, 3, 4, 5, 6, 9 \\
12 & 3, 4  \\
14 & 2, 3, 4 \\     
18 & 1, 2 \\
26 & 1
\end{array}            
\]
\vspace*{0.5mm}

\section{Automorphic forms on orthogonal groups} \label{brel}

We recall some results about automorphic forms on orthogonal groups. In particular we describe Borcherds' singular theta correspondence. It maps modular forms for the Weil representation to automorphic forms on orthogonal groups. Borcherds determined their expansions at $0$-dimensional cusps. The expansions at $1$-dimensional cusps were computed by Kudla.

\medskip

Let $L$ be an even lattice of signature $(n,2)$ with $n>2$ and $V = L \otimes_{\Z} \Q$. Then
\[ {\cal K} = \{ [z] \in P(V(\C)) \, | \, (z,z) = 0, \, (z,\ov{z}) < 0  \}
\]
is a complex manifold with two connected components which are exchanged under the map $z \mapsto \ov{z}$. Let ${\cal H}$ be one of these components. The group $\text{O}(V(\R))$ acts on ${\cal K}$
and the index-$2$ subgroup $\text{O}(V(\R))^+$ preserving ${\cal H}$ consists of the elements with positive spinor norm. 
Let $\Gamma \subset \text{O}(L)$ be a subgroup of finite index and $\Gamma^+ = \Gamma \cap \text{O}(V(\R))^+$. Then the quotient $\Gamma^+ \backslash {\cal H}$ is a complex analytic space that can be compactified by means of the Baily-Borel compactification \cite{BB66} by adding finitely many $0$- and $1$-dimensional rational boundary components. They correspond to the $\Gamma^+$-orbits of the $1$- and $2$-dimensional isotropic subspaces $U$ of $V$. More precisely the $0$- and $1$-dimensional cusps are given by the $\Gamma^+$-orbits of ${\cal C}(U) = \{ [z] \in P(V(\C)) \, | \, \text{$z$ generates $U(\C)$} \}$ and ${\cal C}(U) = \{ [z] \in P(V(\C)) \, | \, \text{$z,\ov{z}$ generate $U(\C)$} \}$, respectively.

Each cusp gives a realisation of ${\cal H}$ as a tube domain. We describe this first for the $0$-dimensional rational boundary components.
Let $U$ be a $1$-dimensional isotropic subspace of $V$ and $U'$ an isotropic subspace of $V$ dual to $U$ under $( \, , \, )$. Choose basis vectors $e_1$ and $e_1'$, respectively, such that $(e_1,e_1')=1$. Let $W$ be the orthogonal complement of $U + U'$. Then $V = W + U + U'$ and $W$ has signature $(n-1,1)$. We write $w + ae_1 + be_1' \in V$ as $(w,a,b)$. The map
\[   \{ z = x + iy \in W(\C) \, | \, (y,y) < 0 \}  \to {\cal K}, \;  z \mapsto [z_L]   \]
with $z_L = (z,-(z,z)/2,1)$
is biholomorphic and we define $H_{U,U'}$ as the component of $\{ z = x + iy \in W(\C) \, | \, (y,y) < 0 \}$ which is mapped to ${\cal H}$.
The $0$-dimensional boundary component associated to $U$ corresponds to the limit of $[(tz)_L]$ for $t\to \infty$.
The group $\text{O}(V(\R))^+$ acts naturally on ${\cal H}$ and we can define an action on $H_{U,U'}$ by making the diagram
\[
\begin{tikzcd}
    {\cal H}      \arrow[r] & {\cal H}  \\
    H_{U,U'} \arrow[u] \arrow[r] & H_{U,U'} \arrow[u] 
\end{tikzcd}
\]
commutative. We also define a cocycle $j$ on $\text{O}(V(\R))^+ \times H_{U,U'}$ by setting $j(\sigma,z) = (\sigma(z_L),e_1)$. Let $k \in \frac{1}{2} \Z$ and $\chi$ a multiplier system
of weight $k$ for $\Gamma^+$. A meromorphic function $\psi: H_{U,U'} \to \C$ is called an {\em automorphic form} of weight $k$ for $\Gamma^+$ with multiplier system $\chi$ if
\[  \psi(\sigma z) = \chi(\sigma) j(\sigma,z)^k \psi(z)  \]
for all $\sigma \in \Gamma^+$, $z \in H_{U,U'}$.

Let $\psi$ be a holomorphic automorphic form on $H_{U,U'}$. Then by the Koecher boundedness principle $\psi$ is also holomorphic on the boundary. Moreover $\psi$ is either constant or of weight $k \geq (n-2)/2$. If $\psi$ vanishes at the boundary, then $\psi$ is called a {\em cusp form}. In this case $\psi = 0$ or $k > (n-2)/2$. If $k = (n-2)/2$, we say that $\psi$ has {\em singular weight}. Then the Fourier expansion of $\psi$ at a $0$-dimensional cusp is supported only on isotropic vectors.

\medskip

Next we describe Borcherds' \cite{B98} construction of automorphic forms on $\text{O}_{n,2}(\R)$. Let $L$ be an even lattice of signature $(n,2)$, $n>2$, $n$ even and let $F$ be a modular form for the Weil representation of $\text{SL}_2(\Z)$ on $\C[D]$, $D=L'/L$ of weight $1-n/2$ with integral principal part. We denote the Fourier coefficients of $F_{\gamma}$ by $c_{\gamma}(m)$. Let $U$ be a $1$-dimensional isotropic subspace of $V$ (there is an assumption on $U'$ that we describe below). Borcherds computes the regularised theta lift of $F$ in the coordinates of $H_{U,U'}$ described above. In a neighbourhood of the cusp defined by $U$ the regularised theta integral is the logarithm of the absolute value of a holomorphic function $\psi$.
This function extends to a meromorphic function on $H_{U,U'}$ and is an automorphic form for $\text{O}(L,F)^+$ of weight
$c_0(0)/2 \in \frac{1}{2}\Z$
with respect to some multiplier system of finite order. The divisor of $\psi$ is determined by the principal part of $F$. More precisely the zeros or poles of $\psi$ lie on rational quadratic divisors $\lambda^{\perp}$ where $\lambda$ is a primitive vector of positive norm in $L$. The divisor $\lambda^{\perp}$ has order 
\[ \sum_{\substack{x \in \Q_{>0} \\x \lambda \in L'}} c_{x \lambda}(-x^2\lambda^2/2) \, . \]
Let $e_1$ be a generator of $U \cap L$ and $e_1'$ a generator of $U'$ such that $e_1' \in L' + U$. 

\begin{thm} \label{queen}
In a neighbourhood of the cusp ${\cal C}(U)$ the function $\psi$ has a product expansion which up to a constant is given by
\[ e((z_L,\rho)) 
\prod_{\substack{\al \in L' \\[0.03mm] (\alpha,e_1)=0 \\[0.03mm] (\al,C) > 0 \\[0.03mm] \hspace*{-2mm} \mod L\cap\Q e_1} } 
            \big( 1-e(-(\alpha,z_L)) \big)^{c_{\al}(-\al^2/2)} \, .  \]
Here $C$ is a Weyl chamber in $W(\R)$ and $\rho$ is the corresponding Weyl vector.          
\end{thm}

\noindent
Because of these expansions $\psi$ is also called the {\em automorphic product} corresponding to $F$. In the following we will often write $\psi_F$ for $\psi$.

We also remark that we slightly modified Borcherds' notation in order to
match Kudla's in \cite{K}. 

Bruinier's converse theorem states that if $L$ splits two hyperbolic planes over $\Z$, 
then an automorphic form for the discriminant kernel of $\text{O}(L)^+$ whose divisor 
is supported on a union of rational quadratic divisors is an automorphic product (see \cite{Br02}, \cite{Br14}).

\medskip

Kudla \cite{K} showed that $\psi_F$ also has a product expansion at the $1$-dimensional cusps which we will describe now in more detail. Let $U$ be a $2$-dimensional isotropic subspace of $V$ with basis $(e_1,e_2)$. Choose an isotropic subspace $U'$ of $V$ which is dual to $U$. Let $(e_1',e_2')$ be a basis of $U'$ such that $(e_i,e_j') = \delta_{ij}$. Define $W$ as the orthogonal complement of $U + U'$. Then $V = W + U + U'$ and $W$ is positive-definite.
We write $w + ae_1 + be_2 + ce_1' + d e_2' \in V$ as $(w,a,b,c,d)$ and define a biholomorphic map
\begin{gather*}
   \{ z = w - \tau_1 e_2 + \tau_2 e_2' \in W(\C) + \C e_2 + \C e_2' \,
     | \, \im(\tau_1) \im(\tau_2) > (\im(w),\im(w))/2  \}  \\
\hspace*{50mm} \to {\cal K} 
\end{gather*}
by 
\[   z \mapsto [z_L] = [(w, \tau_1 \tau_2 - (w,w)/2, -\tau_1 , 1, \tau_2) ]     \]
and choose $H_{U,U'}$ as the component of the domain which is mapped to ${\cal H}$.
We can assume that this choice corresponds to $\im(\tau_1)>0$ so that the $1$-dimensional boundary component associated to $U$ is given by the limit of $[z_L]$ for  $\im(\tau_1)\to \infty$. We define the Jacobi theta function 
\[\vartheta(z,\tau) = \sum_{n\in\Z} e^{i \pi \left( n+1/2 \right)^2 \tau + 2\pi i\left(n+1/2\right)\left(z-1/2\right)}  \, . \]
It is related to the Dedekind eta function $\eta$ by
\[  \vartheta(z,\tau) / \eta(\tau) = -i q^{1/12} (\zeta^{1/2} - \zeta^{-1/2}) \prod_{n = 1}^{\infty} (1-\zeta q^n)(1-\zeta^{-1}q^n)\]
where $ q = e(\tau) $ and $ \zeta = e(z) $. Now we can describe Kudla's product expansion of $\psi_F$. For this suppose $U \cap L = \Z e_1+\Z e_2$.
\begin{thm} \label{kudlaprod}
In a neighbourhood of the $1$-dimensional cusp ${\cal C}(U)$ the automorphic form $\psi_F$ is a product of the following four factors
\begin{enumerate}[i)]
\item
\[ \prod_{\substack{\alpha\in L' \\[0.03mm] (\alpha,e_1)=0 \\[0.03mm] (\al,e_2)>0 \\[0.03mm] \hspace*{-2mm} \mod L\cap \Q e_1 }}
  \big( 1-e(-(\alpha,z_L)) \big)^{c_\alpha(-\alpha^2/2)} \, , \]
\item
\[ \prod_{\substack{ \al\in L'\cap U^{\perp} \\[0.03mm] \hspace*{-4mm} \mod L\cap U \\[0.03mm] \hspace*{-2mm} (\alpha,C)>0}}
  \left(\frac{\vartheta(-(\alpha,z_L),\tau_2)}{\eta(\tau_2)}e((\alpha,z_L)
    -(\alpha_U,z_L)/2)^{(\alpha,e'_2)}\right)^{c_\alpha(-\alpha^2/2)} \, ,  \]
where $\al_U = (\al, e'_1)e_1 + (\al,e'_2)e_2$ and $C$ a Weyl chamber in $W(\R)$,
\item
\[ \prod_{\substack{\al\in L'\cap U \\[0.03mm] \hspace*{-3mm} \mod L\cap U \\[0.03mm] \al\ne 0}}
  \left(\frac{\vartheta(-(\alpha,z_L),\tau_2)}{\eta(\tau_2)}e((\alpha,z_L)/2)^{(\alpha,e'_2)}\right)^{c_\alpha(0)/2}  \]
\item
and
\[  \kappa \eta(\tau_2)^{c_0(0)}q_1^{I_0}\]
where $\kappa$ is a constant of absolute value $1$ and
\[ I_0  = - \sum_{m \in \Q} \,
\sum_{\substack{\alpha\in L'\cap U^\perp \\[0.3mm] \hspace*{-3mm}\mod L\cap U}} c_\al(-m) \sigma_1(m-\alpha^2/2) \, . \]
\end{enumerate}
\end{thm}
\noindent 
Here $\sigma_1$ is the usual sum-of-divisors function with $\sigma_1(r)=0$ if $r \notin \Z_{\geq 0}$ and $\sigma_1(0)=-1/24$. If we write the Fourier-Jacobi expansion of $\psi_F$ at ${\cal C}(U)$ as
\[   \psi_F(z) = q_1^{I_0} \sum_{m=0}^{\infty} \psi_m(w,\tau_2) q_1^m \]
then 
\begin{multline*}
\psi_0(w,\tau_2) = \kappa \eta(\tau_2)^{c_0(0)} \prod_{\substack{\al\in L'\cap U \\[0.03mm] \hspace*{-3mm} \mod L\cap U \\[0.03mm] \al\ne 0}}
\left(\frac{\vartheta(-(\alpha,z_L),\tau_2)}{\eta(\tau_2)}e((\alpha,z_L)/2)^{(\alpha,e'_2)}\right)^{c_\alpha(0)/2}  \\
\prod_{\substack{ \al\in L'\cap U^{\perp} \\[0.03mm] \hspace*{-4mm} \mod L\cap U \\[0.03mm] \hspace*{-2mm} (\alpha,C)>0}}
  \left(\frac{\vartheta(-(\alpha,z_L),\tau_2)}{\eta(\tau_2)}e((\alpha,z_L)
    -(\alpha_U,z_L)/2)^{(\alpha,e'_2)}\right)^{c_\alpha(-\alpha^2/2)}  \, . 
\end{multline*}
Special cases of Kudla's expansion were already described by Gritsenko and Nikulin (cf.\ \cite{GN1}, \cite{GN2} and \cite{G18}). Their result was generalised by Wang and Williams in \cite{WW22}.

\section{Reflective forms of singular weight} \label{DM}

Reflective automorphic products are automorphic products whose divisor is a union of simple zeros along reflection hyperplanes of the underlying lattice. Since they are lifts of reflective modular forms, the bounds in Theorem \ref{thebeatles} imply that at most 474 lattices can support such an automorphic form. Out of those
$132$ can satisfy the Eisenstein condition for singular weight.
Pairing with cusp forms we see that there are at most $11$ lattices carrying a reflective automorphic product of singular weight. They are naturally related to $11$ conjugacy classes in Conway's group $\text{Co}_0$. We use this correspondence to construct on each of the lattices a reflective automorphic product of singular weight. Applying again obstruction theory and some combinatorial arguments we show that they are unique.

\medskip

In this section we will denote the Fourier coefficients of a vector-valued modular form $F = \sum_{\gamma \in D} F_{\gamma} e^{\gamma}$ for the Weil representation or its dual representation by $[F_{\gamma}](m)$.

\subsection*{Reflective forms}

For automorphic forms on orthogonal groups there are several closely related notions of reflectivity. We explain the definition that we use in this paper.

\medskip

Recall that a root of a rational lattice $L$ is a primitive vector $\al \in L$ of positive norm such that the reflection $\sigma_{\al}$ is in $\text{O}(L)$.

Let $L$ be an even lattice of signature $(n,2)$, $n>2$, $n$ even. An automorphic form $\psi$ for a subgroup of $\text{O}(L)^+$ is called {\em geometrically reflective} if it is holomorphic and the support of its divisor is contained in $\bigcup \al^{\perp}$ where $\al$ ranges over the roots of $L$. If $L$ splits $I\!I_{1,1}$ and $\psi_F$ is a geometrically reflective automorphic product on $L$ then $F$ is weakly reflective in the sense of Section \ref{modforms}. For our purposes the following definition is adequate. We say that an automorphic product $\psi_F$ on $L$ is {\em reflective} if it is the theta lift of a reflective modular form $F = \sum_{\gamma \in D} F_{\gamma}$. Such an automorphic form is then also geometrically reflective as the following proposition shows.

\begin{prp} \label{mbig}
	Let $L$ be an even lattice of signature $(n,2)$, $n>2$, $n$ even and $\psi_F$ a reflective automorphic product on $L$. Let $\lambda \in L$ be a primitive vector of norm $\lambda^2 = 2d >0$. Then the divisor $\lambda^{\perp}$ has order
	\[ [F_{\lambda/d}](-1/d) + [F_{\lambda/2d}](-1/4d)  \]
	if $\lambda/2d \in L'$,
	\[ [F_{\lambda/d}](-1/d) \]
	if $\lambda/d \in L'$ but $\lambda/2d \notin L'$, and order $0$ otherwise.
	In particular the order of $\lambda^\perp$ is $0$, $1$ or $2$. If $\lambda^\perp$ has positive order, then $\lambda$ is a root of $L$. 
	Furthermore for $\lambda \in L$ of norm $\lambda^2=2$ the divisor $\lambda^{\perp}$ has positive order. 
\end{prp}

{\em Proof:} 
Write $\lambda = m \mu$ with $\mu \in L'$ primitive and $m \in \Z_{>0}$. Then $2d/m = (\lambda,\mu) \in \Z$, i.e.\ $m |2d$.

Suppose $[F_{k\mu}](-k^2\mu^2/ 2) \neq 0$ for some $k \in \Z_{>0}$. Then $ck \mu \in L$ and $k^2\mu^2/2 = 1/c$ for some $c \in \Z_{>0}$ by the reflectivity of $F$. It follows $ck \mu^2 \in \Z$ and $k|2$.

We determine the relation between $m$ and $k$. Since $ck\mu = (ck/m) \lambda \in L$ and $\lambda$ is primitive we have $m|ck$. The equation $2d = \lambda^2 = m^2 \mu^2 = 2m^2/ck^2$ implies $m = dk(kc/m)$ so that $dk|m$.

Altogether we have $dk|m|2d$. It follows $m=d$ or $2d$. If $m=d$ then $k=1$ and $\lambda^{\perp}$ has order $[F_{\mu}](-\mu^2/2)$. If $m=2d$ then $k=1$ or $2$ and
$\lambda^{\perp}$ has order $[F_{\mu}](-\mu^2/2) + [F_{2\mu}](-2\mu^2)$. In all other cases the order of $\lambda^{\perp}$ vanishes.

If the order of $\lambda^{\perp}$ is positive, then $\lambda/d \in L'$ so that $d$ divides the level of $L$. But then $\lambda$ is a root of $L$ (see Proposition 2.2 in \cite{S06}). The last statement follows from our definition of reflective (see Section \ref{modforms}). \eop

\subsection*{Necessary conditions}

Using the bounds for reflective modular forms given in Theorem \ref{thebeatles} and obstruction theory we show that there are at most $11$ lattices carrying a reflective automorphic product of singular weight.

\medskip

Let $D$ be a discriminant form of even signature and level $N$ carrying a reflective modular form $F = \sum_{\gamma \in D} F_{\gamma} e^{\gamma}$. For $d|N$ we define the sets of singular components 
\[  M_d = \{ \gamma \in D_{d,1/d} \, | \, \text{$F_{\gamma}$ singular} \}  \, . \]
We derive necessary conditions for the existence of $F$ by pairing $F$ with lifts $G$ of modular forms $g$ for $\Gamma_0(N)$ on $0$. Since the components $G_{\gamma}$ of $G$ depend only on the invariants $\{ (c,\q_c(\gamma)) \, | \, c|N, \gamma \in D^{c*} \}$ (cf.\ Section \ref{modforms}), we decompose $D_{d,1/d}$ with respect to these invariants.

\begin{prp}
Let $D$ be a discriminant form of even signature and level $N$ with $2^j||N$. Let $d|N$. If $d$ is odd then all $\gamma \in D_{d,1/d}$ have the same invariants. Suppose $d$ is even. Then $D_{d,1/d}$ decomposes with respect to the above invariants as 
\[  D_{d,1/d} = {\cal O}_{d,1,1/d} \, \cup \bigcup_{\substack{2|c|2^j \\ x \in \Q/\Z }} {\cal O}_{d,c,x}   \]
where 
${\cal O}_{d,c,x} = \{  \gamma \in D_{d,1/d} \, | \, \text{$\gamma \in D^{c*}_x$ and $\gamma \notin D^{m*}$ for $2c|m$}  \}$. Furthermore for $2|c|2^j$ the sets ${\cal O}_{d,c,x}$ are empty if $D$ contains no odd $2$-adic Jordan components.
\end{prp}

To simplify notations we define ${\cal O}_d = D_{d,1/d}$ if $d$ is odd and ${\cal O}_d = {\cal O}_{d,1,1/d}$ if $d$ is even.

\medskip

Next we describe upper bounds for the weights of reflective automorphic products.

\begin{prp} \label{sepultura}
Let $L$ be a regular even lattice of signature $(n,2)$, $n>2$ and even, splitting $I\!I_{1,1} \oplus I\!I_{1,1}$. Suppose $L$ carries a reflective automorphic product $\psi_F$ of weight $k$. Then $k$ is bounded above by
\[
\renewcommand{\arraystretch}{1.2}
\begin{array}{c|c|c|c|c|c|c|c|c|c|c|c|c}
  n  & 4 & 6 & 8 & 10 & 12 & 14 & 16 & 18 & 20 & 22 & 24 & 26  \\  \hline
  \text{bound}  & 54 & 96 & 120  & 252 &  53  & 64 & - & 132 & - & - & - & 12
\end{array}            
\]
\end{prp}
{\em Proof:} Since the proof is similar to the proof of Theorem \ref{classref}  we only sketch it.
Let $L$ be as specified and $\psi_F$ a reflective automorphic product of weight $k$ on $L$. Then $k \geq -1 + n/2$ and $[F_0](0) = 2k$. Pairing $F$ with the Eisenstein series $E = \sum_{\gamma \in D} E_{\gamma} e^{\gamma}$ of weight $1 + n/2$ for the dual Weil representation we obtain
\[  2k + \sum_{d|N} \, \sum_{\gamma \in D_{d,1/d}} [F_{\gamma}](-1/d)[E_{\gamma}](1/d) = 0  \, . \]
Since $[F_{\gamma}](-1/d) = 0$ or $1$ and $[E_{\gamma}](1/d) \in \Q_{\leq 0}$ we can bound $2k \in \Z$ by
\[
  n-2 \leq 2k \leq
  - \sum_{d|N} \, \sum_{\gamma \in {\cal O}_d} [E_{\gamma}](1/d)
  - \sum_{2|d|N} \, \sum_{\substack{2|c|2^j \\ x \in \Q/\Z }} \, \sum_{\gamma \in {\cal O}_{d,c,x}} [E_{\gamma}](1/d)  \, .  
\]
Note that the coefficients $[E_{\gamma}](1/d)$ in each sum are constant.
There are $474$ lattices satisfying the assumptions of the proposition and the condition of Theorem \ref{thebeatles}. We calculate the right-hand side for these lattices to get upper bounds for $k$. Then we check for which $k$ we can solve the first equation with $[F_{\gamma}](-1/d) = 0$ or $1$. For some $k$ this knapsack problem has no solution. The remaining weights satisfy the given bounds.  \eop

\medskip
\noindent
We give bounds for the individual levels in the Appendix. In signature $(8j+2,2)$, $j=1,2$ or $3$ the maximal weight is attained by the theta lift of $E_4^{3-j}/\Delta$ on $I\!I_{8j+2,2}$.

\begin{thm} \label{classref}
Let $L$ be a regular even lattice of signature $(n,2)$, $n>2$ and even, splitting $I\!I_{1,1} \oplus I\!I_{1,1}$. Suppose $L$ carries a reflective automorphic product $\psi_F$ of singular weight. Then $L$ is one of the following $11$ lattices:
\[ \renewcommand{\arraystretch}{1.2}
\begin{array}{c|c|c|c}
\text{lattice} & \text{decomposition}  & \text{cardinality} & \text{intersection} \\
               & \text{of $D_{d,1/d}$}  &                     & \text{with $M_d$}    \\ \hline
               &                       &                     &                      \\ [-4.5mm]
  I\!I_{26,2}   &   {\co}_1             &  1 &   *                                  \\ \hline
                &                       &                     &                     \\ [-4.5mm]
   I\!I_{18,2}(2_{I\!I}^{+10})   &   {\co}_1, {\co}_2    &  1, 496    &   *,*  \\ [0.1mm] \hline
               &                       &                     &                      \\ [-4.5mm]
    I\!I_{14,2}(3^{-8})             &  {\co}_1, {\co}_3       &  1, 2214      &  *,*  \\ [0.1mm] \hline
               &                       &                     &                      \\ [-4.5mm]
   I\!I_{14,2}(2_{I\!I}^{-10}4_{I\!I}^{-2})     &  {\co}_1, {\co}_2, {\co}_4  & 1,2112,6144  &   *,264,*   \\ [0.1mm] \hline
               &                       &                     &                      \\ [-4.5mm]
    I\!I_{12,2}(2_2^{+2} 4_{I\!I}^{+6})      &   {\co}_1, {\co}_{2,2,0},     &  1,36,          &  *,0,  \\  
               &   {\co}_{2,2,1/2}, {\co}_4   &  28,4160       &  *,4032  \\  \hline
                &                       &                     &                      \\ [-4.5mm]
  I\!I_{10,2}(5^{+6})                       &  {\co}_1, {\co}_5    &   1,   3100      &  *,* \\ [0.1mm] \hline
               &                       &                     &                      \\ [-4.5mm]
    I\!I_{10,2}(2_{I\!I}^{+6} 3^{-6})              & {\co}_1, {\co}_2, {\co}_3, {\co}_6 & 1, 28, 234, 6552 & *,*,*,* \\ [0.2mm] \hline
               &                       &                     &                      \\ [-4.5mm]
    I\!I_{8,2}(7^{-5})                             & {\co}_1, {\co}_7 & 1, 2352 & *,*  \\[0.3mm] \hline
               &                       &                     &                      \\ [-4.5mm]
I\!I_{8,2}(2_7^{+1} 4_7^{+1} 8_{I\!I}^{+4}) & {\co}_1, {\co}_{2,4,0}, {\co}_{2,4,1/2},         & 1, 10,6, & *,0,*, \\
               & {\co}_{4,2,1/4}, {\co}_{4,2,3/4}, {\co}_{8}         &  120,136,4032   & 0,120,3840    \\  \hline
               &                       &                     &                      \\ [-4.5mm]  
 I\!I_{8,2}(2_{I\!I}^{+4} 4_{I\!I}^{-2} 3^{+5})   & {\co}_1, {\co}_2, {\co}_3, & 1,24,72,      & *,12,*,  \\ 
                                              & {\co}_4, {\co}_6, {\co}_{12} & 96,2160,6912 & *,1080,* \\ \hline
               &                       &                     &                      \\ [-4.4mm] 
  I\!I_{6,2}(2_{I\!I}^{-2} 4_{I\!I}^{-2} 5^{+4}) & {\co}_1, {\co}_2, {\co}_4,            & 1,12,24,      & *,4,*,  \\
               & {\co}_5, {\co}_{10}, {\co}_{20}        & 120,1440,2880 & *,480,*
\end{array} \]

Here the second column describes the decomposition of $D_{d,1/d}$, $d|N$ where $N$ is the level of $L$, the third column gives the number of elements in each component and the fourth column the cardinality of the intersections with $M_d$. A $*$ indicates that the numbers are the same. 
\end{thm}

Before we proceed to the proof, we describe an example. If there is a reflective automorphic product $\psi_F$ on the lattice $L$ of genus $I\!I_{12,2}(2_2^{+2} 4_{I\!I}^{+6})$ then the singular sets $M_d$ are given by $M_1 = \{ 0 \}$, $M_2 = {\cal O}_{2,2,1/2} = D^{2*}_{1/2}$ and $M_4\subset D_{4,1/4}$. These sets have cardinalities $1$, $28$ and $4032$.

\medskip

{\em Proof:}
Theorem \ref{thebeatles} restricts the signature and the level of a regular lattice carrying a reflective automorphic product. For each possible signature and level we determine all regular lattices splitting $I\!I_{1,1} \oplus I\!I_{1,1}$ by writing down their genus symbol (cf.\ \cite{CS} and \cite{AGM}). We find $474$ lattices. 

Let $L$ be one of these lattices, $(n,2)$ its signature, $N$ its level, $2^j||N$ and suppose $L$ carries a reflective automorphic product $\psi_F$ of singular weight. Then $F$ is a reflective modular form for the Weil representation $\rho_D$ of weight $1-n/2$ and any modular form $G$ for the dual Weil representation $\ov{\rho}_D$ of weight $1 + n/2$ imposes restrictions on $F$ as explained in Section \ref{modforms}. We construct such forms as lifts of scalar-valued modular forms $g$ for $\Gamma_0(N)$ on $0$. We assume that $g$ vanishes at each cusp except possibly at $\infty$. If $G$ is such a lift then the condition coming from $G$ is
\[[F_0](0)[G_0](0) + \sum_{d|N} \, \sum_{\gamma \in D_{d,1/d}} [F_{\gamma}](-1/d)[G_{\gamma}](1/d) = 0 \, . \]
In order to evaluate this formula we need to compute the Fourier coefficients of $G$. This can be done as follows (cf.\ \cite{S09}). Let $s \in \Gamma_0(N) \backslash P = \Gamma_0(N) \backslash (\Q \cup \{\infty\})$ be a cusp of $\Gamma_0(N)$. Choose a representative $a/c$ of $s$ with $c|N$ and $(a,N) = 1$ and a matrix $M_s = \left(\begin{smallmatrix} a  & b \\ c & d \end{smallmatrix}\right) \in \text{SL}_2(\Z)$ such that $d = 0 \! \mod N/N_c$ where $N_c$ is the smallest Hall divisor of $N$ divisible by $c$. Then
\[ g|_{1+n/2, M_s}(\tau) = \sum_{n=0}^{\infty} b_s(n) q_{m_s t_s}^n \]
is an expansion of $g$ at $s$.
As above $t_s = N/(c^2,N)$ is the width of $s$ and
$m_s= 1$ or $2$ the order of $\ov{\chi}_D(T_s)$ (cf.\ Proposition 5.1 in \cite{S09}). We abbreviate $a_s(n) = \overline{\xi}(M_s^{-1}) b_s(n)$ where $\xi$ is a root of unity coming from the Weil representation (see Theorem 3.7 in \cite{S15}). Then the coefficient of $G_{\gamma}$ at $q^n$ for $n \in \Z + \q(\gamma)$ is given by
\[  [G_{\gamma}](n) = \sum_{\substack{s \in \Gamma_0(N) \backslash P \\ \gamma \in D^{c*}}} t_s
  \frac{\sqrt{|D_{c}|}}{\sqrt{|D|}} \, a_s(m_s t_s n) e(-d \q_{c}(\gamma))  \, . \]
Our choice of $d$ in the matrix $M_s$ implies $e(- d\q_{c}(\gamma)) = 1$ for $\gamma \in \mathcal{O}_{d}$. 

First we take the Eisenstein series $E_{1+n/2, \ov{\chi}_D}$ from Section \ref{modforms} for $g$. Then $[G_0](0) = 1$. We determine the sets of invariants $\mathcal{O}_d$ and $\mathcal{O}_{d,c,x}$ and for each such set the Fourier coefficient $E_\gamma(1/d)$ of the Eisenstein series $E$. These computations are based on a program written by Opitz (cf.\ \cite{O}) for SageMath (Version 8.1). We obtain the Eisenstein condition
\begin{multline*}
  (n-2) + [E_0](1)
  + \sum_{d|N,\, d>1} \; \sum_{\gamma \in {\cal O}_d} [F_{\gamma}](-1/d)[E_{\gamma}](1/d)  \\
  + \sum_{2|d|N} \, \sum_{\substack{2|c|2^j \\ x \in \Q/\Z }} \, \sum_{\gamma \in {\cal O}_{d,c,x}}  [F_{\gamma}](-1/d)[E_{\gamma}](1/d) = 0
\end{multline*}
where $[E_{\gamma}](1/d) \in \Q_{\leq 0}$ (see Proposition 5.3 in \cite{BrK}).
We have to determine whether this equation is solvable with $[F_{\gamma}](-1/d)=0$ or $1$. This is a bounded knapsack problem with maximum capacity $(n-2) + [E_0](1)$ where the variables are the cardinalities $| M_d \cap \mathcal{O}_{d}|$ and $|M_d \cap \mathcal{O}_{d,c,x}|$ with bounds $|\mathcal{O}_d|$ and $|\mathcal{O}_{d,c,x}|$, respectively, and the weights are the Eisenstein coefficients $-[E_\gamma](1/d)$. We reduce the problem to a $0-1$ knapsack problem by introducing binary variables (cf.\ Section 7.1.1 in \cite{KPP}) and find that exactly $132$ lattices can solve the Eisenstein condition. In general there are different solutions for the cardinalities $| M_d \cap \mathcal{O}_{d}|$ and $|M_d \cap \mathcal{O}_{d,c,x}|$ for a given lattice. There are no solutions in signature $(4,2)$.  

Next we lift cusp forms $g$ for $\Gamma_0(N)$ on $0$. We obtain conditions
\begin{multline*} 
\sum_{d|N} \, \sum_{\gamma \in {\cal O}_d} [F_{\gamma}](-1/d)[G_{\gamma}](1/d) \\
+ \sum_{2|d|N} \, \sum_{\substack{2|c|2^j \\ x \in \Q/\Z }} \, \sum_{\gamma \in {\cal O}_{d,c,x}}  [F_{\gamma}](-1/d)[G_{\gamma}](1/d) = 0  \, . 
\end{multline*}
Only the $11$ lattices and the corresponding cardinalities for $|{\cal O}_{d} \cap M_d|$ and $|{\cal O}_{d,c,x} \cap M_d|$ given in the theorem satisfy the additional restrictions.

We describe the restrictions coming from cusp forms in more detail for the two most complicated cases.

In signature $(8,2)$ and level $12$ there are $29$ lattices satisfying the Eisenstein condition. They fall into two classes depending on their character.

If the $3$-rank of $D$ is odd then
\[   \chi_D(M) = \Big( \frac{a}{3} \Big)  \]
for $M = \left(\begin{smallmatrix} a  & b \\ c & d \end{smallmatrix}\right) \in \Gamma_0(12)$. We lift the $4$ cusps forms
\[   g_1 = T_2^2 \eta_{1^3 3^3 6^2 12^2}, \quad g_2 = T_2( \eta_{4^1 8^4 12^3}\theta_{A_1}^2 ),
  \quad g_3 = T_3 g_1, \quad g_4 = T_3 g_2   \]
in $S_5(\Gamma_0(12), \overline{\chi}_D)$ and calculate the corresponding conditions using PARI/GP \cite{P}. Together with the Eisenstein condition they exclude the following discriminant forms:
\[ \renewcommand{\arraystretch}{1.2}
\begin{array}{l}
g_1           \\ \hline
              \\ [-4.5mm]
  2_{I\!I}^{-2} 4_{I\!I}^{+4} 3^{-1} \quad 2_{I\!I}^{+2} 4_{I\!I}^{-4} 3^{+1} \quad 4_{I\!I}^{-2} 3^{+5} \quad 4_{I\!I}^{+4} 3^{-3} \quad 2_{0}^{+2}4_{I\!I}^{+2}3^{+5} \quad 4_{I\!I}^{-4} 3^{-3} \\
  \\
  g_1, g_2  \\ \hline
            \\ [-4.5mm]  
2_{I\!I}^{-4} 4_{I\!I}^{+2} 3^{+3} \quad 2_{I\!I}^{+2} 4_{I\!I}^{+2} 3^{+5} \quad 2_{I\!I}^{+2} 4_{I\!I}^{-2} 3^{+5} \quad 4_{I\!I}^{-6} 3^{-3} \\
  \\
  g_1, g_2,g_3  \\ \hline
               \\ [-4.5mm]  
  2_{I\!I}^{+4} 4_{I\!I}^{-2} 3^{-3} \quad 2_{4}^{-4} 4_{I\!I}^{+2} 3^{-3} \quad 2_{0}^{-4} 4_{I\!I}^{+2} 3^{+3} \quad 2_{I\!I}^{-2} 4_{I\!I}^{-2} 3^{-5} \quad 2_{I\!I}^{-2} 4_{I\!I}^{+2} 3^{-5}  \\
  2_{I\!I}^{-2} 4_{I\!I}^{+4} 3^{+3} \quad 2_{I\!I}^{+2} 4_{I\!I}^{-4} 3^{-3} \quad 2_{0}^{-4} 4_{I\!I}^{+2} 3^{-5} \quad 2_{I\!I}^{-4} 4_{I\!I}^{+2} 3^{-5} \quad 2_{I\!I}^{-2} 4_{I\!I}^{+4} 3^{-5}
\end{array}
\]
For the discriminant form $2_{I\!I}^{+4} 4_{I\!I}^{-2} 3^{+5}$ pairing with $g_1, \ldots, g_4$ determines the cardinalities $|{\cal O}_d \cap M_d|$.

If the $3$-rank of $D$ is even, then $\chi_D(M) = (-1)^{(a-1)/2}$ for $M = \left(\begin{smallmatrix} a  & b \\ c & d \end{smallmatrix}\right) \in \Gamma_0(12)$. Here we lift the cusp forms
\[
\begin{array}{lll}
h_1 = \eta_{1^2 2^2 3^2 6^2} \theta_{A_1}^2,  & h_2 = T_3 h_1,  & h_3 = \eta_{2^2 4^2 6^2 12^2} \theta_{A_1}^2, \\[1mm]
h_4 = T_3 T_2 \eta_{2^6 3^4}, & h_5 = T_2\eta_{1^6 3^2 6^2}, & h_6 = T_3 h_3 \, . 
\end{array}
\]
The corresponding conditions exclude the following discriminant forms:
\[ \renewcommand{\arraystretch}{1.2}
\begin{array}{l|l}
\multicolumn{2}{l}{h_1}            \\ \hline
\multicolumn{2}{l}{}               \\ [-4.5mm]
\multicolumn{2}{l}{ 2_{2}^{+2} 4_{I\!I}^{+2} 3^{-4} \quad 2_{6}^{+2} 4_{I\!I}^{+2} 3^{+4}  \quad 2_{2}^{+2} 4_{I\!I}^{+4} 3^{+2} \quad 2_{6}^{+2} 4_{I\!I}^{+4} 3^{-2}} \\
\multicolumn{2}{l}{}              \\
h_1, h_2, h_3, h_4 & h_1, h_2, h_3, h_4, h_5 \\ \hline
           &    \\ [-4.5mm]
  2_{6}^{+4}4_{I\!I}^{+2}3^{+4} \quad 2_{2}^{+2} 4_{I\!I}^{+2} 3^{+6} &   2_{2}^{+4} 4_{I\!I}^{+2} 3^{-4}  \\
\multicolumn{2}{l}{}   \\
\multicolumn{2}{l}{h_1, h_2, h_3, h_4,h_5,h_6}   \\ \hline 
\multicolumn{2}{l}{}      \\ [-4.5mm]
\multicolumn{2}{l}{2_{2}^{+4} 4_{I\!I}^{+2} 3^{+6}}
\end{array}
\]

In signature $(6,2)$ and level $36$ there are $18$ lattices which satisfy the Eisenstein condition. If the exponent $n_3$ at the prime $3$ in the genus symbol is odd, then the character $\chi_D$ is given by
$\chi_D(M) = (-1)^{(a-1)/2} \left( \frac{a}{3} \right)$ for 
a matrix
$M = \left(\begin{smallmatrix} a  & b \\ c & d \end{smallmatrix}\right) \in \Gamma_0(36)$
and we can use the cusp forms $T_n T_2 \eta_{1^3 6^3 9^1 18^1}$, $n=1,5$ and $T_n T_2^2 \eta_{2^4 4^1 12^3}$, $n=1,2$ to eliminate these cases.
For $n_3$ even, $\chi_D$ is trivial and the cusp forms $T_n \eta_{2^3 6^2 18^3}$, $n=1,5$ and $\eta_{1^1 2^1 3^1 12^1 18^1 36^1} \theta_{A_1}^2$ exclude the corresponding genera.

The cusp forms for the remaining cases are described in the Appendix.
\eop

\subsection*{Existence}

The $11$ lattices in Theorem \ref{classref} are naturally related to certain conjugacy classes in $\text{Co}_0$. We use this correspondence to construct on each of the lattices a reflective automorphic product of singular weight.

\medskip

Let $\Lambda$ be the Leech lattice. The orthogonal group of $\Lambda$ is Conway's group $\text{Co}_0$. The quotient $\text{Co}_1 = \text{Co}_0/\la -1 \ra$ is a sporadic simple group. For $g \in \text{Co}_0$ of order $m$ and cycle shape $\prod_{d|m} d^{\, b_d}$ we define the eta product
\[  \eta_g(\tau) = \prod_{d|m} \eta(d\tau)^{b_d}  \, . \]
The level $N$ of $g$ is defined as the level of $\eta_g$. Then $m|N$ and we denote $h = N/m$. If the fixed-point sublattice $\Lambda^g$ is non-zero its level divides the level of $g$ and $\Lambda^g$ is the unique lattice in its genus with maximal minimal norm (see \cite{S04}, Theorem 5.2). The group $\text{Co}_0$ has $72$ conjugacy classes with non-trivial fixed-point lattice. They fall into $70$ algebraic conjugacy classes.  

\begin{thm} \label{cabedelo}
Let $L$ be one of the lattices in Theorem \ref{classref}. Then there is a unique class in $\text{Co}_0$ such that $L$ has index $
h$ in $\Lambda^g \oplus I\!I_{1,1} \oplus I\!I_{1,1}(m/h)$.
\end{thm}
{\em Proof:}
By going through the classes we easily see that for each $L$ in Theorem \ref{classref} there is a unique class in $\text{Co}_0$ such that $\rk(L)=\rk(M)$ and $|L'/L| = h^2 |M'/M|$ for $M = \Lambda^g \oplus I\!I_{1,1} \oplus I\!I_{1,1}(m/h)$. We verify that $L$ is a sublattice of $M$. \eop

\medskip

We list the classes in the following table. The names are taken from Table 4 in \cite{H17}.

\[ \renewcommand{\arraystretch}{1.2}
\begin{array}{c|c|c|c|c|c|c}
  \text{name} & \text{genus of $L$} & m &\text{cycle shape} & h & \text{genus of $\Lambda^g$} & \text{class} \\ \hline
              &                     &   &                   &   &                             &      \\[-4mm]
     A      &  I\!I_{26,2}                                  & 1 &     1^{24}           &  1  & I\!I_{24,0} & 1A   \\ 
     B      &  I\!I_{18,2}(2_{I\!I}^{+10})                    & 2 &    1^8 2^8          &  1  & I\!I_{16,0}(2_{I\!I}^{+8}) & 2A   \\ 
     C      &  I\!I_{14,2}(3^{-8})                           & 3 &     1^6 3^6          &  1  & I\!I_{12,0}(3^{+6}) & 3B   \\ 
     D      &  I\!I_{14,2}(2_{I\!I}^{-10} 4_{I\!I}^{-2})       & 2 &     2^{12}            &  2  & I\!I_{12,0}(2_{4}^{+12}) & 2C   \\ 
     E      &  I\!I_{12,2}(2_2^{+2} 4_{I\!I}^{+6})            & 4 &     1^4 2^2 4^4       &  1  & I\!I_{10,0}(2_2^{+2}4_{I\!I}^{+4}) & 4C  \\
     F      &  I\!I_{10,2}(5^{+6})                           & 5 &     1^4 5^4           &  1  & I\!I_{8,0}(5^{+4}) & 5B   \\ 
     G      &  I\!I_{10,2}(2_{I\!I}^{+6} 3^{-6})              & 6 &     1^2 2^2 3^2 6^2   &  1  & I\!I_{8,0}(2_{I\!I}^{+4}3^{+4}) & 6E   \\
     H      &  I\!I_{8,2}(7^{-5})                            & 7 &    1^3 7^3           &  1  & I\!I_{6,0}(7^{+3}) & 7B   \\ 
     I      &  I\!I_{8,2}(2_7^{+1} 4_7^{+1} 8_{I\!I}^{+4})     & 8 &     1^22^14^18^2      &  1  & I\!I_{6,0}(2_7^{+1} 4_7^{+1} 8_{I\!I}^{+2}) & 8E   \\ 
     J      &  I\!I_{8,2}(2_{I\!I}^{+4} 4_{I\!I}^{-2} 3^{+5})   & 6 &     2^3 6^3          &  2  & I\!I_{6,0}(2_4^{-6} 3^{-3}) & 6G   \\ 
     K      &  I\!I_{6,2}(2_{I\!I}^{-2} 4_{I\!I}^{-2} 5^{+4})   & 10 &    2^2 10^2         &  2  & I\!I_{4,0}(2_4^{+4} 5^{+2}) & 10F
\end{array} \]

\vspace*{2mm}
\noindent
Note that the classes all have balanced cycle shapes, i.e.\ $b_d = b_{N/d}$ if we set $b_d= 0$ for $d \, | \hspace{-0.4em}/ m$.

In the cases with $h = N/m = 2$ we can write
\[  L =  \Lambda^g_N \oplus I\!I_{1,1} \oplus I\!I_{1,1}(m/h)   \]
with $\Lambda^g_N = \la \, \al \in \Lambda^g \,| \, \al^2 = 4  \text{ or }  N \, \ra \subset \Lambda^g$.

We remark that for each class in $\text{Co}_0$ with non-trivial fixed-point sublattice there is a reflective automorphic product of singular weight on a sublattice
of index $h$ in $\Lambda^g \oplus I\!I_{1,1} \oplus I\!I_{1,1}(m/h)$.
The details will be presented elsewhere.

We now describe the constructions of the reflective automorphic products of singular weight.

\medskip

If $L$ has squarefree level $N$, we take an automorphism $g$ of $\Lambda$ of cycle shape $\prod_{d|N} d^{\, 24/\sigma_1(N)}$ and lift $f=1/\eta_g$ on $0$. Then $F = F_{f, 0}$ is a reflective modular form on the indicated discriminant form and $F_0$ has constant coefficient $24\sigma_0(N)/\sigma_1(N)$ (cf.\ \cite{S04}, \cite{S06}). The theta lift $\psi_F$ is a reflective automorphic product of singular weight on $L$.

\medskip

The cases $E$ and $I$ of level $N=4$ and $8$, respectively, are similar but more complicated. Here
\[  F = F_{f, 0} + \frac{N}{16} \, F_{T_2 f, D^{N/2}}  \]
with $f=1/\eta_g$ as above is the desired reflective modular form. For more details we refer to Section 7 in \cite{S09}. There, a different Jordan decomposition for the case $I$ was chosen. Sign walking gives the decomposition used here.

\medskip

We consider the case $D$ (cf.\ \cite{S15}). Let $g$ be an automorphism of the Leech lattice $\Lambda$ of cycle shape $2^{12}$. Then $g$ has level $4$. The fixed-point sublattice $\Lambda^g$ is isomorphic to $D_{12}^+(2)$ and is the unique lattice in the genus $I\!I_{12,0}(2_{4}^{+12})$ of minimum $4$. The group $\text{O}(\Lambda^g)$ has $7$ orbits on the discriminant form of $\Lambda^g$ which are described in the following tables:

\[ \renewcommand{\arraystretch}{1.2}
\begin{array}{c|c|c|ccc|c|c|c}
  \text{norm} & \text{length} & \text{order} & \text{name} & &
                \text{norm} & \text{length} & \text{order} & \text{name}      \\ \cline{1-4} \cline{6-9}  
  0  &  1     &    1     &  0_0   &   \quad &   1/4  &   1024   &    2   &  1    \\ \cline{1-4} \cline{6-9}  
  0  &  1     &    2     &  0_A   &   &   1/2  &    132   &    2   &  2_A   \\  
     &  990   &    2     &  0_B   &   &        &    924   &    2   &  2_B   \\  \cline{6-9}
\multicolumn{4}{}{}                &   &   3/4  &   1024   &    2   &  3        
\end{array}
\]

\vspace*{2mm}
\noindent
The orbit $2_A$ is generated by the $264$ elements $\al \in {\Lambda^g}' \cap (\Lambda^g/2) = {\Lambda^g}'$ of norm $\al^2 = 1$. 
The theta functions of the orbits $0_0$ and $0_B$ are given by
\begin{align*} 
\theta_{0_0}(\tau)  \, = \, \, & \theta_{\Lambda^g}(\tau) = \, 
 1 + 264q^2 + 2048q^3 + 7944q^4 + 24576q^5 + 64416q^6 + \ldots  \\
\theta_{0_B}(\tau)  \, = \, \,
 &8q + 256q^2 + 1952q^3 + 8192q^4 + 25008q^5 + 62464q^6 + \ldots 
\end{align*}
We define
\[  h(\tau) 
= \frac{1}{\eta_g(\tau/2)}
= q^{-1/2} + 12q^{1/2} + 90q^{3/2} + 520q^{5/2} + 2535q^{7/2} + \ldots \]

The lattice $\Lambda^g_4$ is isomorphic to $D_{12}(2)$ and has genus $I\!I_{12,0}(2_{I\!I}^{-10} 4_{I\!I}^{-2})$. Define $L = D_{12}(2) \oplus I\!I_{1,1}\oplus I\!I_{1,1}$. Then $L$ has genus $I\!I_{14,2}(2_{I\!I}^{-10} 4_{I\!I}^{-2})$.

We choose a subset $M$ of the discriminant form $D$ of $L$ with the following properties:
\begin{enumerate}[i)]
\item $|M|=264$
\item $M \subset D_{2,1/2}$ 
\item $M = D^2 + M$
\item Let $\gamma \in D_{2,1/2}$. Then
\[ | M \cap \gamma^{\perp} | =
\begin{cases} 
\, 184 & \text{if $\gamma \in M$,} \\
\, 120 & \text{otherwise.} 
\end{cases}  \]
\end{enumerate}
We can construct such a set as follows. Consider the embeddings
\[   \Lambda^g_4 \subset \Lambda^g \subset {\Lambda^g}' \subset {\Lambda^g_4}'    \]
Then $H = \Lambda^g/\Lambda^g_4$ is an isotropic subgroup of ${\Lambda^g_4}'/\Lambda^g_4$ with orthogonal complement $H^{\perp} = {\Lambda^g}'/\Lambda^g_4 \subset {\Lambda^g_4}'/\Lambda^g_4$. The quotient $H^{\perp}/H$ is naturally isomorphic to ${\Lambda^g}'/\Lambda^g$. The pullback of the orbit $2_A = \{ \al + \Lambda^g \, | \, \al \in {\Lambda^g}', \, \al^2 = 1 \} \subset {\Lambda^g}'/\Lambda^g$ under the projection $H^{\perp} \to H^{\perp}/H$ embeds naturally into the discriminant form $D$ of $L$. This set then has the desired properties.
We will see later that a subset of $D$ with the above properties is unique modulo $\text{O}(D)$.

Now define 
\[  F = F_{\theta_{0_0}/\Delta, 0} - F_{\theta_{0_B}/\Delta, 0}
  + \frac{1}{12}\sum_{\gamma \in M} F_{h,\gamma}             \]
Then
\begin{align*} 
F_0
& = q^{-1} + 12 + 300q + 5792q^2 + 84186q^3 + 949920q^4 + 8813768q^5 + \ldots \\
& = \frac{1}{\eta_g}  
+ \frac{1}{2} \left( \frac{\theta_{\Lambda^g}}{\Delta} - \frac{1}{\eta_g} \right)
\intertext{and}
F_{\gamma} & = \frac{\theta_{2_A}}{2\Delta} = 
q^{-1/2} + 44q^{1/2} + 1242q^{3/2} + 22216q^{5/2} + 287463q^{7/2} + \dots 
\intertext{if $\gamma \in M$,}
F_{\gamma} & = \frac{\theta_{3}}{2\Delta} = 
q^{-1/4} + 90q^{3/4} + 2535q^{7/4} + 42614q^{11/4} + 521235q^{15/4} + \dots 
\end{align*}
if $\gamma^2/2 = 1/4 \! \mod 1$. All the other components of $F$ are holomorphic at $\infty$. The singular sets $M_d$ are thus given by
\[
\renewcommand{\arraystretch}{1.2}
\begin{array}{c|c|c}
	M_1 & M_2 & M_4  \\  \hline
	D_{1,1/1} & M & D_{4,1/4}
\end{array}            
\]
It follows that $F$ is reflective.
For a different construction of $F$ see \cite{S15}.

\medskip

The construction generalises easily to the case $J$.
Here we choose an automorphism $g$ of $\Lambda$ with cycle shape $2^36^3$. Then $g$ has level $12$ and the fixed point sublattice $\Lambda^g$ is the unique lattice in the genus $I\!I_{6,0}(2_4^{-6} 3^{-3})$ with minimal norm $4$. The group $\text{O}(\Lambda^g)$ has $7$ orbits on the subgroup $D_2$ of $D$. They are described in the following tables:

\[ \renewcommand{\arraystretch}{1.2}
\begin{array}{c|c|c|ccc|c|c|c}
  \text{norm} & \text{length} & \text{order} & \text{name} & &
                \text{norm} & \text{length} & \text{order} & \text{name}      \\ \cline{1-4} \cline{6-9}  
  0  &  1     &    1     &  0_0   &   \quad &   1/4  &   16   &    2   &  1    \\ \cline{1-4} \cline{6-9}  
  0  &  1     &    2     &  0_A   &   &   1/2  &    6   &    2   &  2_A   \\  
     &  18    &    2     &  0_B   &   &        &    6   &    2   &  2_B   \\  \cline{6-9}
\multicolumn{4}{}{}               &   &   3/4  &   16   &    2   &  3        
\end{array}
\]

\vspace*{2mm}
\noindent
Again the orbit $2_A$ is generated by the $12$ elements $\al \in {\Lambda^g}' \cap (\Lambda^g/2)$ of norm $\al^2 = 1$. 
The theta functions of the orbits $0_0$ and $0_B$ are given by
\begin{align*} 
\theta_{0_0}(\tau)  \, = \, \,
  & 1 + 12q^2 + 32q^3 + 42q^4 + 96q^5 + 84q^6 + 96q^7 + 300q^8 + \ldots \\
  \theta_{0_B}(\tau)  \, = \, \,
  & 2q + 16q^2 + 26q^3 + 32q^4 + 96q^5 + 112q^6 + 100q^7 + 256q^8 + \ldots
\end{align*}
As above we define $h(\tau) = 1/\eta_g(\tau/2)$.

The lattice $\Lambda^g_{12}$ has genus $I\!I_{6,0}(2_{I\!I}^{+4} 4_{I\!I}^{-2} 3^{-3})$. Let $L = D_4(6) \oplus A_2(2) \oplus I\!I_{1,1} \oplus I\!I_{1,1}$. Then $L$ has genus $I\!I_{8,2}(2_{I\!I}^{+4} 4_{I\!I}^{-2} 3^{+5})$ and is isomorphic to $\Lambda^g_{12} \oplus I\!I_{1,1} \oplus I\!I_{1,1}(3)$.

We choose a subset $M$ of the discriminant form $D$ of $L$ with the following properties:
\begin{enumerate}[i)]
\item $|M|=12$
\item $M \subset D_{2,1/2}$
\item $M = D^6 + M$
\item Let $\gamma \in D_{2,1/2}$. Then
\[ | M \cap \gamma^{\perp} | =
\begin{cases} 
      \, 4  & \text{if $\gamma \in M$,} \\
      \, 12 & \text{otherwise.} 
\end{cases} \]
\end{enumerate}
Such a set can be constructed exactly in the same way as in the case $2^{12}$ as the pullback of the orbit $2_A$. Its uniqueness modulo $\text{O}(D)$ is easy to see.

Define
\[  F = F_{\theta_{0_0}/\eta_{g^2}, 0} - F_{\theta_{0_B}/\eta_{g^2}, 0}
  + \frac{1}{24} \sum_{\gamma \in M} F_{h,\gamma}  \, .            \]
Then
\begin{align*} 
F_0
& = q^{-1} + 6 + 480q + 20192q^2 + 472068q^3 + 7504260q^4 + 91178456q^5 + \ldots \\
& = \sum_{k|6} \sum_{d|k} \frac{\mu(k/d)}{k} \, \frac{\theta_{\Lambda^{g,d}}}{\eta_{g^d}}
\intertext{and}
F_{\gamma} & = q^{-1/2} + 44q^{1/2} + 3258q^{3/2} + 102280q^{5/2} + 1949277q^{7/2} + \dots 
\intertext{if $\gamma \in M$,}
F_{\gamma}
  & = q^{-1/3} + 104q^{2/3} + 6233q^{5/3} + 173448q^{8/3} + 3087720q^{11/3} + \dots
\intertext{if $\gamma \in D_3$ and $\q(\gamma) = 1/3 \! \mod 1$,}
F_{\gamma}
  & = q^{-1/4} + 144q^{3/4} + 8259q^{7/4} + 222398q^{11/4} + 3857625q^{15/4} + \dots
\intertext{if $\gamma \in D_4$ and $\q(\gamma) = 1/4 \! \mod 1$,}
F_{\gamma} & = q^{-1/6} + 220q^{5/6} + 11276q^{11/6} + 287584q^{17/6} + 4831653q^{23/6} + \dots
\intertext{if $\gamma \in (D_{2,1/2} \backslash M) + D_{3,2/3}$ and finally}
F_{\gamma}
  & = q^{-1/12} + 296q^{11/12} + 14829q^{23/12} + 366730q^{35/12} + 6013842q^{47/12} + \dots 
\end{align*}
if $\q(\gamma) = 1/12 \! \mod 1$. All the other components of $F$ are holomorphic at $\infty$. Hence the singular sets of $ F $ are

\[ \renewcommand{\arraystretch}{1.2}
\begin{array}{c|c|c|c|c|c}
	M_1 & M_2 & M_3 & M_4 & M_6 & M_{12} \\ \hline
	D_{1,1/1} & M & D_{3,1/3} & D_{4,1/4} & (D_{2,1/2} \backslash M) + D_{3,2/3} &D_{12,1/12}  
\end{array}
\]

\vspace*{2mm}
\noindent
and $F$ is again reflective.

\medskip

The case $K$ is slightly different from the previous two constructions. We choose an automorphism $g$ of $\Lambda$ with cycle shape $2^210^2$. Then $g$ has level $20$ and the fixed point sublattice $\Lambda^g$ is the unique lattice in the genus $I\!I_{4,0}(2_4^{+4} 5^{+2})$ with minimum $4$. The group $\text{O}(\Lambda^g)$ now has $6$ orbits on the subgroup $D_2$ of $D$. They are given by

\[ \renewcommand{\arraystretch}{1.2}
\begin{array}{c|c|c|ccc|c|c|c}
  \text{norm} & \text{length} & \text{order} & \text{name} & &
                \text{norm} & \text{length} & \text{order} & \text{name}      \\ \cline{1-4} \cline{6-9}  
  0  &  1     &    1     &  0_0   &   \quad &   1/4  &    4   &    2   &  1    \\ \cline{1-4} \cline{6-9}  
  0  &  1     &    2     &  0     &   &   1/2        &    2   &    2   &  2_A   \\  
 \multicolumn{4}{}{}              &   &              &    4   &    2   &  2_B   \\  \cline{6-9}
\multicolumn{4}{}{}               &   &   3/4        &    4   &    2   &  3        
\end{array}
\]

\vspace*{2mm}
\noindent
Analogously to the previous cases the orbit $2_A$ is generated by the $4$ elements $\al \in {\Lambda^g}' \cap (\Lambda^g/2)$ of norm $\al^2 = 1$.
The theta functions of the orbits $0_0$ and $2_B$ are given by
\begin{align*} 
\theta_{0_0}(\tau)  \, = \, \,
  & 1 + 4q^2 + 8q^3 +  4q^4+ 16q^5 + 16q^6 + 8q^7+ 4q^8 + 16q^9 + \ldots \\
  \theta_{2_B}(\tau)  \, = \, \,
  &  4q^{3/2} + 8q^{5/2} + 4q^{7/2} + 8q^{9/2} + 8q^{11/2} + 8q^{13/2} + 28q^{15/2} + \ldots
\end{align*}

\vspace*{2mm}
\noindent
The lattice $\Lambda^g_{20}$ has genus $I\!I_{4,0}(2_{I\!I}^{-2} 4_{I\!I}^{-2} 5^{+2})$. Let $L = D_4(10) \oplus I\!I_{1,1}\oplus I\!I_{1,1}$. Then $L$ has genus $I\!I_{6,2}(2_{I\!I}^{-2} 4_{I\!I}^{-2} 5^{+4})$ and is isomorphic to $\Lambda^g_{20} \oplus I\!I_{1,1} \oplus I\!I_{1,1}(5)$.

We choose a subset $M$ of the discriminant form $D$ of $L$ with the following properties:
\begin{enumerate}[i)]
\item $|M|=4$ 
\item $M \subset D_{2,1/2}$ 
\item $M = D^{10} + M$
\item Let $\gamma \in D_{2,1/2}$. Then
\[ | M \cap \gamma^{\perp} | =
\begin{cases} 
      \, 4  & \text{if $\gamma \in M$,} \\
      \, 0 & \text{otherwise.} 
\end{cases} \]
\end{enumerate}
The existence of such a set can be proved as before. The uniqueness up to $\text{O}(D)$ is easy to see.

Define
\[  F = F_{\theta_{0_0}/\eta_{g^2}, 0} - \frac{1}{32} \sum_{\gamma \in M} F_{\theta_{2_B}/\eta_{g^2}, \gamma}  \, .            \]
Then
\begin{align*} 
F_0
& = q^{-1} + 4 + 748q + 43040q^2 + 1197138q^3 + 21539168q^4 + \ldots \\
& = \sum_{k|10} \sum_{d|k} \frac{\mu(k/d)}{k} \, \frac{\theta_{\Lambda^{g,d}}}{\eta_{g^d}}
\intertext{and}
F_{\gamma} &= q^{-1/2} + 60q^{1/2} + 6386q^{3/2} + 242792q^{5/2} + 5303951q^{7/2} + \dots
\intertext{if $\gamma \in M$,}
F_{\gamma}
 & = q^{-1/4} + 210q^{3/4} + 16815q^{7/4} + 545294q^{11/4} + 10781475q^{15/4} + \dots 
\intertext{if $\gamma \in D_4$ and $\q(\gamma) = 1/4 \! \mod 1$,}
F_{\gamma}
  & = q^{-1/5} + 280q^{4/5} + 20558q^{9/5} + 641296q^{14/5} + 12413390q^{19/5} + \dots
\intertext{if $\gamma \in D_5$ and $\q(\gamma) = 1/5 \! \mod 1$,}
F_{\gamma} & = q^{-1/10} + 456q^{9/10} + 29830q^{19/10} + 878048q^{29/10} + 16375851q^{39/10} + \dots 
\intertext{if $\gamma \in M + D_{5,3/5}$ and finally}
F_{\gamma}
  & = q^{-1/20} + 558q^{19/20} + 35539q^{39/20} + 1022903q^{59/20} + 18768281q^{79/20} + \dots 
\end{align*}
if $\q(\gamma) = 1/20 \! \mod 1$. All the other components of $F$ are holomorphic at $\infty$, its singular sets are given by

\[ \renewcommand{\arraystretch}{1.2}
\begin{array}{c|c|c|c|c|c}
	M_1 & M_2 & M_4 & M_5 & M_{10} & M_{20} \\ \hline
	D_{1,1/1} & M & D_{4,1/4} & D_{5,1/5} & M + D_{5,3/5} &D_{20,1/20}  
\end{array}
\]

\vspace*{2mm}
\noindent
and $F$ is reflective.

\subsection*{Uniqueness}

We show that the $11$ automorphic products of singular weight constructed above are unique up to isomorphism of the underlying lattice. The proof is based on obstruction theory and some combinatorial arguments.

\medskip

We consider four different cases the most complicated being the one corresponding to the lattice $I\!I_{14,2}(2_{I\!I}^{-10} 4_{I\!I}^{-2})$. We start with the simplest case.

\begin{prp}
  Let $L$ be one of the lattices $I\!I_{26,2}$, $I\!I_{18,2}(2_{I\!I}^{+10})$, $I\!I_{14,2}(3^{-8})$, $I\!I_{10,2}(5^{+6})$, $I\!I_{8,2}(7^{-5})$ or $I\!I_{10,2}(2_{I\!I}^{+6} 3^{-6})$. Then $L$ carries a unique reflective automorphic product $\psi_F$ of singular weight.
\end{prp}
{\em Proof:} In these cases the sets $M_d = \{ \gamma \in D_{d,1/d} \, | \, \text{$F_{\gamma}$ singular} \}$ are completely fixed by the condition coming from the Eisenstein pairing (cf.\ Theorem \ref{classref}) and are given by $M_d = D_{d,1/d}$. This implies that the principal part of $F$ is uniquely determined and therefore also the modular form $F$. \eop

\medskip
In the remaining cases we need additional restrictions coming from pairing $F$ with cusp forms to prove the uniqueness.

\begin{prp}\label{uniquenessEI}
  Let $L$ be the lattice $I\!I_{12,2}(2_2^{+2} 4_{I\!I}^{+6})$ or $I\!I_{8,2}(2_7^{+1} 4_7^{+1} 8_{I\!I}^{+4})$. Then $L$ carries a unique reflective automorphic product $\psi_F$ of singular weight.
\end{prp}
{\em Proof:} We start with the case $I\!I_{12,2}(2_2^{+2} 4_{I\!I}^{+6})$.
We know from Theorem \ref{classref} that $M_2 = {\co}_{2,2,1/2} = D^{2*}_{1/2}$. In order to determine $M_4$ we pair $F$ with the lift of the cusp form $\eta_{1^4 2^2 4^4}\theta_{A_1^4}$ on $\mu \in D^2 \backslash \{ 0 \}$ (cf.\ Section \ref{modforms}). We find
\[ 64 + \sum_{\gamma \in D_{4,1/4,2}} [F_{\gamma}](-1/4)
    + \sum_{\gamma \in D_{4,1/4,4}} [F_{\gamma}](-1/4) e((\gamma,\mu)) = 0  \, . \]
Since $|M_4| = 4032$ we also have
\[ \sum_{\gamma \in D_{4,1/4,2}} [F_{\gamma}](-1/4)
  +  \, \sum_{\gamma \in D_{4,1/4,4}} [F_{\gamma}](-1/4) = 4032   \]
so that
\[  \sum_{\gamma \in D_{4,1/4,4}} [F_{\gamma}](-1/4) (1-e((\gamma,\mu))) = 4096 \, . \]
The inner product of $\gamma \in D_{4,1/4,4}$ and $\mu$ is $(\gamma,\mu) = 0$ or $1/2 \! \mod 1$ so that $1-e((\gamma,\mu)) = 0$ or $2$. It is easy to verify that
\[  \sum_{\gamma \in D_{4,1/4,4}} (1-e((\gamma,\mu))) = 4096 \, . \]
This implies $[F_{\gamma}](-1/4) = 1$ for $\gamma \in D_{4,1/4,4}$ with $(\gamma,\mu) = 1/2 \! \mod 1$. Since for each $\gamma \in D_{4,1/4,4}$ there is an element $\mu \in D^2$ such that $(\gamma,\mu) = 1/2 \! \mod 1$ we get $[F_{\gamma}](-1/4) = 1$ for all $\gamma \in D_{4,1/4,4}$. Now $|D_{4,1/4,4}| = |M_4| = 4032$ implies $D_{4,1/4,4} = M_4$. Hence $F$ is unique and therefore also $\psi_F$.

The argument for $I\!I_{8,2}(2_7^{+1} 4_7^{+1} 8_{I\!I}^{+4})$ is similar. Here we lift the cusp forms $(T_2-2)g$ and $(T_3+2)g$ with
$g = T_2^2 \eta_{1^2 2^1 4^3 8^4}$ on $\mu \in D^4 \backslash \{ 0 \}$ to get equations similar to the ones above. They can be used to show $M_4= D^{2*}_{3/4} \cap D_{4,1/4,4}$ and $M_8 = D_{8,1/8,8}$. \eop

\medskip

The cases with $N/m=2$ are more difficult. Here the singular sets are unique only up to $\text{O}(D)$.

\begin{prp}\label{uniquenessJK}
Let $L$ be the lattice $I\!I_{8,2}(2_{I\!I}^{+4} 4_{I\!I}^{-2} 3^{+5})$ or $I\!I_{6,2}(2_{I\!I}^{-2} 4_{I\!I}^{-2} 5^{+4})$. Then $L$ admits a unique reflective automorphic product of singular weight up to $\text{O}(L)^+$.
\end{prp}
{\em Proof:}
We start with the lattice $I\!I_{8,2}(2_{I\!I}^{+4} 4_{I\!I}^{-2} 3^{+5})$. Then
$M_j = D_{j,1/j}$ for $j=1$, $3$, $4$ and $12$. We have to determine the sets $M_2$ and $M_6$. Choose $\mu \in D_{2,1/2}$ and let $g \in S_5(\Gamma_1(12),\overline{\chi}_{\mu})$. We lift $g$ to a vector-valued modular form $G$ for the dual Weil representation $\overline{\rho}_D$. Then
\[  \sum_{j|12} \, \sum_{\gamma \in D_{j,1/j}} [F_{\gamma}](-1/j) [G_{\gamma}](1/j) = \sum_{j|12} R_j = 0  \, . \]
We will derive explicit formulas for the summands $R_j$. We represent the cusps of $\Gamma_1(12)$ by the rational numbers $1/1$, $1/5$, $1/2$, $1/3$, $2/3$, $1/4$, $3/4$, $1/6$, $1/12$, $5/12$. For each cusp $s=a/c$ we choose a matrix
$M_s=\left(\begin{smallmatrix}  a & b \\  c & d \\ \end{smallmatrix}\right) \in \text{SL}_2(\Z)$. To simplify the calculations we assume $b = 1 \! \mod 2$ and $d = 0 \! \mod 12/N_c$ where $N_c$ the smallest Hall divisor of $12$ such that $(c,12)|N_c$. This will simplify the calculations. Then
\[  g_s(\tau) = g|_{5,M_s}(\tau) = \sum_{n=1}^\infty b_s(n)q^n_{m_s t_s}  \]
is an expansion of $g$ at the cusp $s$. As before we denote by $t_s$ the width of $s$ and by $m_s$ the order of $\ov{\chi}_{\gamma}(T_s)$ (see Proposition 3.4 in \cite{S15}) and abbreviate $a_s(n) = \overline{\xi}(M^{-1}_s) b_s(n)$.
Then we get the following expressions for the $R_d$
\begin{align*}
R_1 =
  & - \bigg( \frac{1}{6\sqrt{3}} \big( a_{1/1}(12) + a_{1/5}(12) \big) + \frac{1}{2} \big( a_{1/3}(4)-a_{2/3}(4) \big) \bigg) \\
R_2 =
  & - \bigg( \frac{1}{6\sqrt{3}} \big( a_{1/1}(6) + a_{1/5}(6) \big) + \frac{1}{2} \big( a_{1/3}(2)-a_{2/3}(2) \big) \bigg) A_2\\
  & - \bigg( \frac{2}{3\sqrt{3}} \, a_{1/2}(3) + 2 a_{1/6}(1) \bigg)  B_2 \\
  & - \bigg( \frac{2}{3\sqrt{3}} \big( a_{1/4}(3)+a_{3/4}(3) \big) + 2 \big(a_{1/12}(1)+a_{5/12}(1) \big) \bigg) C_2 \\
R_3 =
  & - \frac{1}{6\sqrt{3}} \big( a_{1/1}(4) + a_{1/5}(4) \big) \sum_{\gamma \in D_{3,1/3}} e((\gamma,\mu))\\
  = & - \frac{12}{\sqrt{3}} \big( a_{1/1}(4) + a_{1/5}(4) \big) 
\end{align*}
\begin{align*}  
R_4=
  & - \bigg( \frac{1}{6\sqrt{3}} \big( a_{1/1}(3) + a_{1/5}(3) \big) + \frac{1}{2} \big( a_{1/3}(1) - a_{2/3}(1) \big) \bigg)
    \sum_{\gamma \in D_{4,1/4}}e((\gamma,\mu)) \\
  & = 0 \\
R_6 =
  & - \frac{1}{6\sqrt{3}} \big( a_{1/1}(2) + a_{1/5}(2) \big) \, A_6  - \frac{2}{3\sqrt{3}} a_{1/2}(1) \, B_6 \\
  & - \frac{2}{3\sqrt{3}} \big( a_{1/4}(1) + a_{3/4}(1) \big) \, C_6 \\
R_{12} =
  & -\frac{1}{6\sqrt{3}} \big( a_{1/1}(1) + a_{1/5}(1) \big) \sum_{\gamma \in D_{12,1/12}}e((\gamma,\mu)) = 0
\end{align*}
with
\begin{align*}
A_j & = \sum_{\gamma \in D_{j,1/j}} [F_{\gamma}](-1/j) \, e((\gamma,\mu)) \\
B_j & = \sum_{\gamma \in D_{j,1/j} \cap (\mu + D^2)} [F_{\gamma}](-1/j) \, e \big( 3 q_2(\gamma-\mu) \big) \\
C_j & = \sum_{\gamma \in D_{j,1/j} \cap (\mu + D^4)} [F_{\gamma}](-1/j)
\end{align*}
for $j=2$ and $6$. Note that $C_2 = [F_{\mu}](-1/2) = 1$ or $0$ depending on whether $\mu \in M_2$ or not. Choosing for $g$ the cusp forms
\begin{gather*}
\eta_{1^{-1}3^{-1}4^{10}6^8 12^{-6}}, \, \eta_{1^{4}2^{-7}4^{10}6^9 12^{-6}}, \, \eta_{1^{-1}2^{-2}3^{7}4^{2}6^{10} 12^{-6}}, \, \eta_{1^{-1}2^{-3}3^{7}4^{7}6^{5} 12^{-5}}, \\
\eta_{1^{1}2^{-1}3^{9}4^{1}6^{3}12^{-3}} \in S_5(\Gamma_1(12),\overline{\chi}_{\mu}) 
\end{gather*} 
we get a system of linear equations for $A_j$, $B_j$, $C_j$ of rank $5$ with solutions 
\begin{align*}
A_2 &= -4+a     & B_2 &= -2 + \frac{a}{8} & C_2 &= 1-\frac{a}{16}  \\[1mm]
A_6 &= 1080-90a & B_6 &=  -\frac{45}{2}\,a  & C_6 &= \frac{45}{4}\, a  
\end{align*}
for some $a \in \C$. The values of $A_j$, $B_j$ and $C_j$ are fixed by $C_2 = [F_{\mu}](-1/2)$.

We determine the structure of $M_2$. Let $\mu \in M_2$. Then
\begin{align*}
  -2 = B_2 &= \sum_{\gamma \in \mu + D^6} [F_{\gamma}](-1/2) \, e \big( q_2(\gamma-\mu) \big) \\
           &= [F_{\mu}](-1/2) - \sum_{\gamma \in D^6\backslash \{ 0 \}} [F_{\mu +\gamma}](-1/2)
             = 1 - \sum_{\gamma \in D^6\backslash \{ 0 \}} [F_{\mu + \gamma}](-1/2) \, .
\end{align*}
Since $|D^6|=4$, this implies $[F_{\mu + \gamma}](-1/2)=1$ for all $\gamma \in D^6$, i.e.\ $M_2$ is invariant under addition of $D^6$.
Now $|M_2|=12$, so that
\[   M_2 = \bigcup_{\gamma \in U} ( \gamma + D^6)  \]
for a $3$-element subset $U$ of $M_2$.
The equation for $A_2$ reads
\[  -4 = A_2 = \sum_{\gamma \in M_2} e((\gamma,\mu)) = \sum_{\gamma \in U, \, \bt \in D^6} e((\gamma+\bt,\mu))
  = 4 \sum_{\gamma \in U} e((\gamma,\mu))  \, . \]
This implies that two different elements in $U$ have inner product $1/2 \! \mod 1$. The quotient $D_2/D^6$ is a discriminant form of type
$2_{I\!I}^{+4} = 2_{I\!I}^{-2} \oplus 2_{I\!I}^{-2}$. Since $2_{I\!I}^{-2}$ contains no non-trivial isotropic elements, the image of $M_2$ under the projection $D_2 \to D_2/D^6$ generates one of the copies of $2_{I\!I}^{-2}$. It follows that two different sets $M_2$ are conjugate under $\text{O}(D)$.

Next we consider $M_6$.
Note that $D_{6,1/6} = D_{2,1/2} + D_{3,2/3}$. For $\mu \in M_2$ we have
\[  0 = C_6 = \sum_{\gamma \in D_{6,1/6} \cap (\mu + D^4)} [F_{\gamma}](-1/6)
  = \sum_{\gamma \in \mu + D_{3,2/3}} [F_{\gamma}](-1/6)   \]
so that $[F_{\gamma}](-1/6) = 0$ for all $\gamma \in \mu + D_{3,2/3}$. This implies $[F_{\gamma}](-1/6) = 0$ for all $\gamma \in M_2 + D_{3,2/3}$. This set has $12 \cdot 90 = 1080$ elements. Since $|D_{6,1/6}| = 2 |M_6| = 2160$ we must have $[F_{\gamma}](-1/6) = 1$ for the remaining elements in $D_{6,1/6}$, i.e.
\[  M_6 = ( D_{2,1/2} \backslash M_2 ) + D_{3,2/3} \, . \]
We remark that $D_{2,1/2}$ contains $6 \cdot 4 = 24$ elements.

In summary we have seen that the singular sets $M_d$ are unique up to $\text{O}(D)$. Hence $F$ is unique modulo $\text{O}(D)$. Since the map $\text{O}(L)^+ \to \text{O}(D)$ is surjective for $L$ (cf.\ \cite{N}, Theorem 1.14.2), it follows that $\psi_F$ is unique up to $\text{O}(L)^+$.

The argument for $I\!I_{6,2}(2_{I\!I}^{-2} 4_{I\!I}^{-2} 5^{+4})$ is similar. By lifting the cusp forms
\begin{gather*}
\eta_{2^{-1}4^{1}5^4 10^9 20^{-5}}, \, \eta_{2^{-2}4^{4}5^4 10^6 20^{-4}},
\eta_{1^{3}2^{-6}4^{9}5^{-3}10^{10} 20^{-5}}, \, \eta_{1^{-1}2^{1}4^{4}5^{1}10^{7} 20^{-4}},\\
\eta_{1^{-1}2^{2}4^{1}5^{1}10^{10}20^{-5}} \in S_4(\Gamma_1(20),\overline{\chi}_{\mu})  
\end{gather*}
on $\mu \in D_{2,1/2}$ we obtain $5$ relations for the numbers $A_j$, $B_j$, $C_j$, $j = 2, 10$ with $A_j$, $C_j$ as above and
\[  B_j = \sum_{\gamma \in D_{j,1/j} \cap (\mu + D^2)} [F_{\gamma}](-1/j) \, e \big( 5 q_2(\gamma-\mu) \big)  \, . \]
They can be written as
\begin{align*}
A_2 &= -4+a         & B_2   &= -\frac{a}{4}       & C_2    &= \frac{a}{8}  \\[1mm]
A_{10} &= 480+120a  & B_{10} &=  -30a              & C_{10} &= 15a  
\end{align*}
for some $a \in \C$. The values of $A_j$, $B_j$ and $C_j$ are fixed by $C_2=[F_{\mu}](-1/2)$. 
We consider the set $M_2$. For $\mu \in M_2$ we obtain 
\[-2 = B_2 = 1 - \sum_{\gamma \in D^{10} \backslash \{ 0 \}} [F_{\mu + \gamma}](-1/2)   \]
which implies $M_2= \mu + D^{10}$ because $|M_2|=|D^{10}|=4$.
  Finally we determine $M_{10}$. For $\mu \in M_2$ we have
\[  120 = C_{10} = \sum_{\gamma \in D_{10,1/10} \cap (\mu + D^4)} [F_{\gamma}](-1/10)
  = \sum_{\gamma \in \mu + D_{5,3/5}} [F_{\gamma}](-1/10)   \]
so that $[F_{\gamma}](-1/10) = 1$ for all $\gamma \in \mu + D_{5,3/5}$ because $|D_{5,3/5}|=120$. Hence $[F_{\gamma}](-1/10) = 1$ for all $\gamma \in M_2 + D_{5,3/5}$ . This set has $4 \cdot 120 = 480$ elements. Since $|M_{10}| = 480$, it follows
\[  M_{10} = M_2 + D_{5,3/5} \, . \]
This proves that $F$ is unique up to $\text{O}(D)$ in this case. \eop

\medskip

Finally we consider the case $I\!I_{14,2}(2_{I\!I}^{-10} 4_{I\!I}^{-2})$. In order to prove the uniqueness for this case, we need some preparation.

We consider a discriminant form $D$ of type $2_{I\!I}^{-10}$ and take a subset $U$ of $D$ such that
\[  | U \cap \gamma^{\perp} | =
\begin{cases}
  \, 66  & \text{if $\gamma =0$,} \\
  \, 34  & \text{if $\q(\gamma) = 0 \! \! \mod 1$ and $\gamma \neq 0$,} \\
  \, 46  & \text{if $\q(\gamma) = 1/2 \! \! \mod 1$ and $\gamma \in U$,} \\
  \, 30  & \text{if $\q(\gamma) = 1/2 \! \! \mod 1$ and $\gamma \notin U$}
\end{cases}  
\]
for all $\gamma \in D$. A set with these properties can be constructed starting from the lattice $D_{12}(2)$ as described in the previous subsection. We will show that $U$ is unique up to $\text{O}(D)$. For $\gamma \in D$ and $b = 0$ or $1/2 \! \mod 1$ we define the sets
\[  U_{\gamma,b}   = \{ \mu \in U \, | \, (\mu,\gamma) = b \! \!  \mod 1 \}    \]
and write
\[
  U_{\gamma_1,\ldots,\gamma_n,b_1,\ldots,b_n }  = \bigcap_{i=1}^n U_{\gamma_i,b_i}  \, .
\]
The cardinalities
\[  c(\gamma_1,\ldots,\gamma_n,b_1,\ldots,b_n) = |  U_{\gamma_1,\ldots,\gamma_n,b_1,\ldots,b_n } |
\]
satisfy
\begin{multline*}
  c(\gamma_1,\ldots,\gamma_{n-1},\gamma_n + \gamma_{n+1},b_1,\ldots,b_{n-1},b) =  \\
  c(\gamma_1,\ldots,\gamma_{n-1},\gamma_n,\gamma_{n+1},b_1,\ldots,b_{n-1},0,b) \\
  + c(\gamma_1,\ldots,\gamma_{n-1},\gamma_n,\gamma_{n+1},b_1,\ldots,b_{n-1},1/2,b+1/2)
\end{multline*}
and
\begin{multline*}
  c(\gamma_1,\ldots,\gamma_n,b_1,\ldots,b_n) = 
  c(\gamma_1,\ldots,\gamma_n,\gamma_{n+1},b_1,\ldots,b_n,0) \\
  + c(\gamma_1,\ldots,\gamma_n,\gamma_{n+1},b_1,\ldots,b_n,1/2) \, .
\end{multline*}
These equations imply the following reduction formula
\[  c(\gamma_1,\ldots,\gamma_n,b_1,\ldots,b_n)
  =  \frac{1}{2^{n-1}} \sum_{\substack{S \subset \{ 1,\ldots, n\} \\ S \neq \{ \} }}
  (-1)^{|S| + 1} \, c(\gamma_S,b_S + (|S| + 1)/2)
\]
with $\gamma_S = \sum_{j \in S} \gamma_j$ and $b_S = \sum_{j \in S} b_j \! \mod 1$. 
We apply this formula to show that $U$ is closed under addition.

\begin{prp} \label{sumoftwo}
Let $\gamma_1,\gamma_2 \in U$ with $(\gamma_1,\gamma_2) = 1/2 \! \mod 1$. Then $\gamma_1 + \gamma_2 \in U$. 
\end{prp}
{\em Proof:}
Suppose $\gamma_1 + \gamma_2 \notin U$. We show that the map $(U \cap \gamma_1^{\perp})\backslash (\{ \gamma_1 \} \cup \gamma_2^{\perp}) \to U$, $\mu \mapsto \mu + \gamma_2$ is well-defined. Let $\mu \in (U \cap \gamma_1^{\perp})\backslash (\{ \gamma_1 \} \cup \gamma_2^{\perp})$. Then $\q(\gamma_2+\mu) = 1/2 \! \mod 1$ and $\gamma_2 + \mu \in U$ because otherwise 
\begin{align*}
\MoveEqLeft   
0 \leq c(\gamma_1,\gamma_2,\gamma_2 + \mu,1/2,1/2,0) \\
  =& \frac{1}{4} \big( c(\gamma_1,1/2) + c(\gamma_2,1/2) + c(\gamma_2+ \mu,0) \\
   & - c(\gamma_1+\gamma_2,1/2) - c(\gamma_1 + \gamma_2 + \mu,0) - c(\mu,0) + c(\gamma_1 + \mu,0) \big) \\
  =& \frac{1}{4} \big( 20 + 20 + 30 - 36 - c(\gamma_1 + \gamma_2 + \mu,0) - 46 +34 \big) 
\end{align*}
would imply $c(\gamma_1 + \gamma_2 + \mu,0) \leq 22$ which is impossible. The image of the map lies in $U_{\gamma_1,1/2} \cap U_{\gamma_2,1/2} =  U_{\gamma_1,\gamma_2,1/2,1/2}$. It follows
\[   |U_{\gamma_1,0}| - (1+|U_{\gamma_1,\gamma_2,0,0}|)  \leq |U_{\gamma_1,\gamma_2,1/2,1/2}|   \, . \]
But this contradicts 
\[  c(\gamma_1,\gamma_2,0,0) = \frac{1}{2} \big( c(\gamma_1,0) + c(\gamma_2,0) - c(\gamma_1+\gamma_2,1/2)  \big) = 28  \]
and $c(\gamma_1,\gamma_2,1/2,1/2) =2$. \eop

\medskip
For three elements the situation is as follows.

\begin{prp} \label{sumofthree}
Let $\gamma_1, \gamma_2, \gamma_3 \in U$ be pairwise orthogonal and different. Then $\gamma_1 + \gamma_2 + \gamma_3 \notin U$. 
\end{prp}
{\em Proof:}
Suppose $\gamma_1 + \gamma_2 + \gamma_3 \in U$. Then the reduction formula implies
\begin{align*}
  c(\gamma_1, \gamma_2, \gamma_3,1/2,1/2,0) &= 8  \\
  c(\gamma_1,\gamma_2,1/2,1/2) &= 4
\end{align*}
which contradicts $c(\gamma_1,\gamma_2,\gamma_3,1/2,1/2,0) \leq c(\gamma_1, \gamma_2,1/2,1/2)$. \eop

\medskip

Next we construct from $U$ a basis of $D$ corresponding to a decomposition $2_{I\!I}^{-10} = (2_{I\!I}^{-2})^4 \oplus (2_{I\!I}^{-2})^1$.

\begin{prp} \label{basis}
There exists a basis $(\gamma_1, \dots, \gamma_{10})$ of $D$ such that 
\begin{enumerate}[i)]
\item the elements $\gamma_{2i-1},\gamma_{2i}$ , $i=1,\ldots,5$ generate pairwise orthogonal Jordan blocks $J_i$ of type $2_{I\!I}^{-2}$,
\item $J_i \backslash \{ 0 \} \subset U$ for $i=1,\ldots,4$,
\item $J_5 \cap U = \{ \}$.
\end{enumerate}
\end{prp}
{\em Proof:}
Let $\gamma_1 \in U$. Since $c(\gamma_1,1/2)= 20>0$ we can choose $\gamma_2 \in U_{\gamma_1,1/2}$. Then $\gamma_1+\gamma_2 \in U$ by Proposition \ref{sumoftwo}. Next we take $\gamma_3 \in U_{\gamma_1,\gamma_2,0,0}$ which is possible because
\[ c(\gamma_1,\gamma_2,0,0) = \frac{1}{2} \big( c(\gamma_1,0)+c(\gamma_2,0)-c(\gamma_1+\gamma_2,1/2) \big) = 36>0  \, . \]
We continue in this way. For $\gamma_8$ we have $2$ choices. We list the corresponding cardinalities.
\begin{align*}
  c(\gamma_1,1/2)
  &= 20 \\
  c(\gamma_1,\gamma_2,0,0)
  &= 36 \\
  c(\gamma_1,\gamma_2,\gamma_3,0,0,1/2)
  &= 14 \\
  c(\gamma_1,\gamma_2,\gamma_3,\gamma_4,0,0,0,0)
  &= 15 \\
c(\gamma_1,\gamma_2,\gamma_3,\gamma_4,\gamma_5,0,0,0,0,1/2) 
  &= 8 \\
c(\gamma_1,\gamma_2,\gamma_3,\gamma_4,\gamma_5,\gamma_6,0,0,0,0,0,0)
  &= 3 \\  
  c(\gamma_1,\gamma_2,\gamma_3,\gamma_4,\gamma_5,\gamma_6,\gamma_7,0,0,0,0,0,0,1/2)
  &= 2 \\
  c(\gamma_1,\gamma_2,\gamma_3,\gamma_4,\gamma_5,\gamma_6,\gamma_7,\gamma_8,0,0,0,0,0,0,0,0)
  &= 0
\end{align*}
We calculate these numbers with the reduction formula. In order to determine the numbers $c(\gamma_S,b_S + (|S| + 1)/2)$ which enter the reduction formula, we add the $\gamma_i$ blockwise and apply Propositions \ref{sumoftwo} and \ref{sumofthree}. We see that the recursion stops after $\gamma_8$. Choosing any two non-zero elements in the orthogonal complement of $\la \gamma_1, \ldots, \gamma_8 \ra$ we obtain the desired basis. \eop

\medskip

The previous propositions show that $U$ consists of the 12 elements in the sets $\Gamma_i = J_i \backslash \{ 0 \}$, $i = 1,\ldots,4$ and 54 elements of the form $\mu + \rho$ with $\mu \in \Gamma_5$ and $\rho \in \la \Gamma_1,\ldots,\Gamma_4 \ra$. We now determine the latter elements. Let ${\cal P}_2$ be the set of partitions of $\{ 1,\ldots, 4 \}$ into $2$-element subsets. Clearly $|{\cal P}_2|=3$. 

\begin{prp} \label{murtinheira}
There is a bijection $\Phi: \Gamma_5 \to {\cal P}_2$ such that $U$ consists of the elements
\begin{enumerate}[i)]
\item $\gamma \in \Gamma_i$ for $i = 1,\dots,4$,  
\item $\gamma = \mu + \rho_i + \rho_j$ with $\{i,j\} \in \Phi(\mu)$ and $\rho_i \in \Gamma_i$, $\rho_j \in \Gamma_j$.
\end{enumerate}
\end{prp}
{\em Proof:}
Let $\mu \in \Gamma_5$. Suppose $\mu + \rho \in U$ for some $\rho \in \Gamma_5^{\perp}$. Then $\rho$ is non-zero and isotropic. Write $\rho = \rho_1 + \rho_2 + \rho_3 + \rho_4$ with $\rho_i \in J_i$. Then either all $\rho_i$ are non-zero or exactly two are non-zero. In the first case all elements of the form $\mu + \rho_1 + \rho_2 + \rho_3 + \rho_4$ with $\rho_i \in \Gamma_i$ would be in $U$ by Proposition \ref{sumoftwo} which contradicts $|U| = 66$. Hence $\rho = \rho_i + \rho_j$ with $\rho_i \in \Gamma_i$, $\rho_j \in \Gamma_j$ for some $2$-element subset $\{ i,j \} \subset \{ 1, \ldots,4 \}$. Then for this subset $\{ i,j \}$ again all $9$ elements of the form $\mu + \rho_i + \rho_j$ with $\rho_i \in \Gamma_i$, $\rho_j \in \Gamma_j$ are in $U$. Let $\Gamma_5 = J_5\backslash \{ 0 \} = \{\mu, \mu', \mu''\}$. Then $\mu + \rho \in U_{\mu',\mu'',1/2,1/2}$. Since $c(\mu',\mu'',1/2,1/2) = 18 = (66-12)/3$ this implies that there are exactly $18$ elements $\rho \in \Gamma_5^{\perp}$ such that $\mu + \rho \in U$. Hence there are two subsets $\{ i,j \}, \{ i',j' \}  \subset \{ 1, \ldots,4 \}$ such that the elements $\mu + \rho_i + \rho_j$, $\mu + \rho_{i'} + \rho_{j'}$ are in $U$. We show that these subsets are disjoint. Suppose $i=i'$. Choose $\rho_i \in \Gamma_i$. Then
\begin{multline*}
  c(\mu,\rho_i,0,1/2) 
  \geq |\Gamma_i \backslash \{ \rho_i \}| + |(\mu + \Gamma_i + \Gamma_j)\backslash (\mu + \rho_i + \Gamma_j) | \\
   + |(\mu + \Gamma_i + \Gamma_{j'})\backslash (\mu + \rho_i + \Gamma_{j'}) | = 2 + 2(9-3) = 14
\end{multline*}
which contradicts $c(\mu,\rho_i,0,1/2) = 8$. Finally we show that $\Phi: \Gamma_5 \to {\cal P}_2$ is bijective. Fix a $2$-element subset $\{ i,j \} \subset \{ 1, \ldots,4 \}$. Choose $\rho_i \in \Gamma_i$, $\rho_j \in \Gamma_j$. Let $\mu \in \Gamma_5$ such that $\mu + \Gamma_i + \Gamma_j \subset U$. Then $\mu + \Gamma_i + \Gamma_j \subset U$ contributes $4$ elements to $U_{\rho_i,\rho_j,1/2,1/2}$. Hence
\[  | \{ \mu \in \Gamma_5 \, | \, \mu + \Gamma_i + \Gamma_j \subset U \} | \leq
  \frac{1}{4} c(\rho_i,\rho_j,1/2,1/2) = 1 \, .  \]
This proves the proposition. \eop

\medskip

Propositions \ref{basis} and \ref{murtinheira} imply

\begin{prp} \label{judaspriest}
Let $D$ be a discriminant form of type $2_{I\!I}^{-10}$ and $U$ a subset of $D$ such that
\[  | U \cap \gamma^{\perp} | =
\begin{cases}
  \, 66  & \text{if $\gamma =0$,} \\
  \, 34  & \text{if $\q(\gamma) = 0 \! \! \mod 1$ and $\gamma \neq 0$,} \\
  \, 46  & \text{if $\q(\gamma) = 1/2 \! \! \mod 1$ and $\gamma \in U$,} \\
  \, 30  & \text{if $\q(\gamma) = 1/2 \! \! \mod 1$ and $\gamma \notin U$}
\end{cases}  
\]
for all $\gamma \in D$. Then $U$ is unique modulo $\text{O}(D)$.
\end{prp}

\begin{prp}\label{uniquenessD}
The lattice $L = I\!I_{14,2}(2_{I\!I}^{-10} 4_{I\!I}^{-2})$ admits a unique reflective automorphic product of singular weight up to $\text{O}(L)^+$.
\end{prp}
{\em Proof:}
The sets $M_j$ are given by $D_{j,1/j}$ for $j=1,4$. We have to describe $M_2$. Let $\mu \in D_2$. We construct obstructions by lifting $g \in S_8(\Gamma(4))$ on $\mu$ to a modular form $G$ for the dual Weil representation $\ov{\rho}_D$. Pairing $F$ with $G$ we get
\[   R_1 + R_2 + R_4 = 0  \]
where as above 
\[   R_j = \sum_{\gamma \in D_{j,1/j}} [F_{\gamma}](-1/j) [G_{\gamma}](1/j)  \, . \]
We represent the cusps of $\Gamma(4)$ by the rational numbers $1$, $2$, $3$, $4$, $1/2$ and $1/4$. They all have width $4$. For each of these cusps $s=a/c$ we choose a matrix $M_s=\left(\begin{smallmatrix}  a & b \\  c & d \\ \end{smallmatrix} \right) \in \text{SL}_2(\Z)$. To simplify the calculations we assume $b =  1 \! \mod 2$ and additionally $d = 0 \! \mod 4$ for all cusps $s=a/c$ with $c=1$.

First we choose $\mu \in D_{2,1/2}$. Writing
\[  g_s(\tau) = g|_{8,M_s}(\tau) = \sum_{n=1}^\infty b_s(n) q_4^n  \]
we find
\begin{align*}
  R_1 =
  & - \frac{1}{16} \Big( a_{1}(4) - a_{2}(4) + a_{3}(4) - a_{4}(4) \Big) \\
  R_2 =
  & - \frac{1}{16} \Big( a_{1}(2) - a_{2}(2) + a_{3}(2) - a_{4}(2) \Big) A_2 \\
  & - 4 a_{1/2}(2) B_2 - 8 a_{1/4}(2) [F_{\mu}](-1/2) \\
  R_4 =
  & -\frac{1}{16} \Big( a_{1}(1) - a_{2}(1) + a_{3}(1) -a_{4}(1) \Big) \sum_{\gamma \in D_{4,1/4}} e( (\gamma,\mu)) = 0 
\end{align*} 
with $a_s(n) = \overline{\xi}(M^{-1}_s) b_s(n)$ and
\begin{align*}
A_2 & = \sum_{\gamma \in D_{2,1/2}} [F_{\gamma}](-1/2) \, e((\gamma,\mu)) \\
B_2 & = \sum_{\gamma \in D_{2,1/2} \cap (\mu + D^2)} [F_{\gamma}](-1/2) \, e \big(  q_2(\gamma-\mu) \big) \, .
\end{align*}
Note that $A_2 = 2 |M_2 \cap \mu^{\perp}| -|M_2|$. The eta quotients $\eta_{1^{-12}2^{44}4^{-16}}$, $\eta_{1^{20}2^{-4}}$ yield the relations
\begin{align*}
  24 + A_2 - 128 [F_{\mu}](-1/2) & = 0 \\
  B_2 + 2 [F_{\mu}](-1/2)           & = 0 \, .
\end{align*}
For $\mu \in M_2$ the second equation implies $[F_{\mu+\gamma}](-1/2) = 1$ for all $\gamma \in D^2$, i.e.\ $M_2$ is invariant under addition of $D^2$. Hence there is a subset $U$ of $M_2$ of cardinality $|U| = 66$ such that
\[   M_2 = \bigcup_{\gamma \in U} ( \gamma + D^2)  \, . \]
From the first equation we get $A_2=104$ for $\mu \in M_2$. This implies $|M_2 \cap \mu^{\perp}| = 184$. For $\mu \notin M_2$ we have $A_2 = - 24$ so that $|M_2 \cap \mu^{\perp}| = 120$ in this case. 

Lifting $\eta_{1^8 2^8}$ on $\mu \in D_{2,0}$ gives
\[  8 - \sum_{\gamma \in D_{2,1/2}} [F_{\gamma}](-1/2) \, e((\gamma,\mu)) = 0   \, . \]
This implies $|M_2 \cap \mu^{\perp}| = 136$.

The quotient $D_2/D^2$ is a discriminant form of type $2_{I\!I}^{-10}$. We denote the image of $M_2$ under the projection $D_2 \to D_2/D^2$ also by $U$. Then we have just proved the following properties of $U$:
\[  | U \cap \gamma^{\perp} | =
\begin{cases}
  \, \, 66  & \text{if $\gamma = 0$,} \\
  \, \, 34  & \text{if $\q(\gamma) = 0 \! \! \mod 1$ and $\gamma \neq 0$,} \\
  \, \, 46  & \text{if $\q(\gamma) = 1/2 \! \! \mod 1$ and $\gamma \in U$,} \\
  \, \, 30  & \text{if $\q(\gamma) = 1/2 \! \! \mod 1$ and $\gamma \notin U$.}
\end{cases}  
\]
Hence $U$ is unique up to automorphisms of $D_2/D^2$ by Proposition \ref{judaspriest}. This implies that $M_2$ is unique modulo $\text{O}(D)$. \eop

\subsection*{Classification}

We summarise the results from the previous subsections.

\begin{thm} \label{blazefoley}
There are exactly 11 regular even lattices $L$ of signature $(n,2)$, $n>2$ and even, splitting $I\!I_{1,1} \oplus I\!I_{1,1}$ which carry a reflective automorphic product of singular weight. They are given in the following table:
\[
\renewcommand{\arraystretch}{1.2}
\begin{array}{r|l}
n  & L \\[0.5mm] \hline 
   &   \\[-4mm]
  26 & I\!I_{26,2} \\[0.7mm]
  18 & I\!I_{18,2}(2_{I\!I}^{+10}) \\[0.7mm]
  14 & I\!I_{14,2}(2_{I\!I}^{-10} 4_{I\!I}^{-2}), \, I\!I_{14,2}(3^{-8}) \\[0.7mm]
  12 & I\!I_{12,2}(2_2^{+2} 4_{I\!I}^{+6}) \\[0.7mm]
  10 & I\!I_{10,2}(2_{I\!I}^{+6} 3^{-6}), \, I\!I_{10,2}(5^{+6}) \\[0.7mm]
  8  & I\!I_{8,2}(2_{I\!I}^{+4} 4_{I\!I}^{-2} 3^{+5}), \, I\!I_{8,2}(2_7^{+1} 4_7^{+1} 8_{I\!I}^{+4}), \, I\!I_{8,2}(7^{-5}) \\[0.7mm]
 6  & I\!I_{6,2}(2_{I\!I}^{-2} 4_{I\!I}^{-2} 5^{+4}) 
\end{array}
\]
For each of the lattices the corresponding automorphic product is unique up to $\text{O}(L)^+$. Its zeros are simple zeros orthogonal to reflection hyperplanes.
\end{thm}
{\em Proof:} The only thing left to prove is the statement about the order of the zeros. It is a simple consequence of Proposition \ref{mbig} and the explicit construction of the corresponding modular form for the Weil representation. \eop

\section{Cusps of orthogonal modular varieties}\label{chemicalbrothers}

We show that the
cusps of orthogonal modular varieties can be parametrised by certain double quotients. Using this description we classify the $1$-dimensional cusps of type $0$ of the modular varieties $\text{O}(L,F)^+\backslash {\cal H}$ corresponding to the $11$ re\-flec\-tive automorphic products $\psi_F$ of singular weight. 
They are in bijection with the root systems in Schellekens' list. If the root system is non-trivial, $\psi_F$ vanishes at the corresponding cusp. Since $\psi_F$ has singular weight, $\text{O}(L,F)^+\backslash {\cal H}$ possesses a
$1$-dimensional cusp on which $\psi_F$ is a non-trivial modular form. For each $\psi_F$ we determine such a cusp.

\subsection*{Classification of cusps}

We associate to a cusp of an orthogonal modular variety two invariants, the type and the associated lattice and introduce the notion of a splitting cusp. Then we show that the splitting cusps of a given type and associated lattice are parametrised by a certain double quotient.

\medskip

Let $L$ be an even lattice of signature $(n,2)$, $n>2$ splitting a hyperbolic plane, $V = L \otimes_{\Z} \Q$ and $\Gamma$ a subgroup of $\text{O}(L)$ containing the kernel $\Delta(L)$ of the projection $\text{O}(L) \to \text{O}(L'/L)$. Set $\Gamma^+ = \Gamma \cap \text{O}(L)^+$ as in Section \ref{brel}.

\medskip

Let $U$ be a $2$-dimensional isotropic subspace of $V$. We define the {\em type} of $U$ as the isotropic subgroup
\[  H = (L'\cap U)/(L \cap U) \subset D  \]
of the discriminant form $D = L'/L$. 
The group $H$ defines an even overlattice
\[  L^H = \bigcup_{\gamma \in H} (\gamma + L)  \]
of $L$. Note that $L^H = (L' \cap U) + L$.
The quotient
\[  K = (L^H \cap U^{\perp})/(L^H \cap U)  \, \]
carries a quadratic form and is called the {\em lattice associated with U}. We show that $K$ is a positive-definite even lattice. 

\begin{prp} \label{dt}
There is a basis $(e_1, e_2)$ of $L^H \cap U$ and primitive isotropic elements $e_1', e_2'$ in $L^H$ such that $e_1, e_1'$ and $e_2, e_2'$ generate orthogonal unimodular hyperbolic planes and
\[  L^H = M \oplus \langle e_1, e_1' \rangle \oplus \langle e_2, e_2' \rangle   \]
for a positive-definite lattice $M$ isomorphic to $K$.
\end{prp}
{\em Proof:}
Let $e \in S^H = L^H \cap U $ be primitive. First we show that $e$ has level $1$, i.e.\ $(e,L^H) = \Z$. Since $L^H = (L' \cap U) + L$, we find $ S^H = L' \cap U $, so that $ e $ is primitive in $ L' $. Hence there is $f \in L \subset L^H$ with $ (e, f) = 1 $. It follows $(e,L^H) \supset \Z $. Since $ L^H $ is even, we conclude $ (e,L^H) = \Z $. 

Next we prove the existence of the desired decomposition of $L^H$. Let $e_1$ be a primitive vector in $S^H$. Since $e_1$ has level $1$ we can choose $\tilde{e}_1 \in L^H $ such that $(e_1,\tilde{e}_1) = 1$. Define $e_1' = - a e_1 + \tilde{e}_1$ with $a = \tilde{e}_1^2/2$. Then $e_1'$ is isotropic and $(e_1, e_1') = 1$ so that the lattice $P_1 = \la e_1, e_1' \ra$ is isomorphic to $I\!I_{1,1}$. It follows $L^H = P_1 \oplus P_1^{\perp}$. Next we choose a primitive vector $e_2 \in S^H \cap P_1^{\perp}$. Then $e_2$ has level $1$ and we can construct as before $e_2' \in P_1^{\perp}$ such that $(e_2, e_2') = 1$ and $(e_2')^2= 0$. Then $P_2 = \la e_2, e_2' \ra \simeq I\!I_{1,1}$. We obtain the decomposition $L^H = M \oplus P_1 \oplus P_2$ with $M = (P_1 \oplus P_2)^{\perp}$.

Finally we note that $ L^H \cap U^\perp $ equals the direct sum of $ M $ and $ \langle e_1, e_2 \rangle = L^H \cap U $. This implies $ M \simeq (L^H \cap U^{\perp})/(L^H \cap U) = K $. 
\eop 

\medskip 

The proposition implies that the discriminant form of $K$ is isomorphic to ${L^H}'/L^H = H^{\perp}/H$ so that $K$ has genus $I\!I_{n-2,0}(H^{\perp}/H)$, i.e.\ $H$ determines the genus of $K$ and $K$ determines $H^{\perp}/H$.

We say that $U$ {\em splits} if there is an isotropic subspace $U'$ of $V$ dual to $U$ such that
\[  L = M \oplus ( L \cap U + L \cap U' )  \]
for some positive-definite even lattice $M$. Since in this case $L^H \cap U^{\perp} = M \oplus (L' \cap U)$ and $L^H \cap U = L' \cap U$, the lattice $M$ is isomorphic to the lattice $K$ associated with $U$. 

\begin{prp} \label{mdt}
Suppose $U$ splits. Then there exist bases $(e_1, e_2)$ of $L \cap U$ and $(e_1', e_2')$ of $L \cap U'$ such that $(e_i,e_j') = m_i \delta_{ij}$ with positive integers $m_1|m_2$. In particular $L \cap U + L \cap U'$ is isomorphic to $I\!I_{1,1}(m_1) \oplus I\!I_{1,1}(m_2)$.
\end{prp}
{\em Proof:}
Let $ (f_1, f_2) $ and $ (f_1', f_2') $ be bases of $ L \cap U $ and $ L \cap U' $, respectively, and define  
\[ A = \begin{pmatrix} (f_1, f_1') & (f_1, f_2') \\ (f_2, f_1') & (f_2, f_2') \end{pmatrix}. \]
We now put $ A $ into Smith normal form. More precisely, we choose matrices $ B, C \in \text{GL}_2(\Z) $ such that $ BAC $ has diagonal form $ \begin{psmallmatrix} m_1 & 0 \\ 0 & m_2 \end{psmallmatrix} $ with $m_1|m_2$. Note that $ B $ and $ C $ represent changes of basis of $ L \cap U $ and $ L \cap U' $. Choosing $ (e_1, e_2) $ and $ (e_1', e_2') $ to be the resulting bases completes the proof. 
\eop 

\medskip 

We can transfer the above notions to isotropic planes $S$ in $L$ by considering $S\otimes_{\Z} \Q \subset V$.

Fix a splitting primitive isotropic plane $S$ in $L$. For the proofs of the next two results we choose a decomposition $L = K \oplus I\!I_{1,1}(m_1) \oplus I\!I_{1,1}(m_2)$ with $S \subset K^{\perp}$ and an isotropic basis $(e_1,e_2,e'_1,e'_2)$ of $ K^{\perp}$ as in Proposition \ref{mdt}.

\begin{prp}
Let $T$ be a primitive isotropic plane in $L$. Then $T \in \text{O}(L)S$ if and only if
\begin{enumerate}[i)]
	\item $T$ splits,
	\item $S$ and $T$ have the same type modulo $\text{O}(D)$,
	\item the lattices associated with $S$ and $T$ are isomorphic.  
\end{enumerate}
\end{prp}
{\em Proof:}
If $ T \in \text{O}(L)S $, then the decomposition $ L = K \oplus I\!I_{1,1}(m_1) \oplus I\!I_{1,1}(m_2) $ of $L$ for $S$ gives a decomposition for $T$ and the statements follow.

Suppose $ T $ satisfies i)--iii) above. Proposition \ref{mdt} gives a decomposition $L = K_T \oplus \langle f_1,f_1' \rangle \oplus \langle f_2, f_2' \rangle$ with $f_1, f_1', f_2, f_2'$ isotropic, $(f_i,f_j') = n_i \delta_{ij} $, $ n_1|n_2 $ and $T = \langle f_1, f_2 \rangle$. By assumption $ K \simeq K_T $. It remains to show that $K^{\perp}$ and $K_T^{\perp}$ are isomorphic. The lattices $S$ and $T$ have types $\Z/m_1\Z \times \Z/m_2\Z $ and $\Z/n_1\Z \times \Z/n_2\Z$ with positive integers $m_1|m_2$ and $n_1|n_2$. Since the types are isomorphic, we conclude $m_i = n_i$. Hence the map defined by $ e_i \mapsto f_i $ and $ e_i' \mapsto f_i' $ gives the desired isomorphism. 
\eop 

\medskip 

We denote the stabiliser of $S$ in $\text{O}(L)$ by $\text{O}(L)_S$.

\begin{prp} \label{rb}
The stabiliser $\text{O}(L)_S$ and the discriminant kernel $\Delta(L)$ both contain an automorphism of $L$ not contained in $\text{O}(L)^+$. 
\end{prp}
{\em Proof:}
We first prove the statement for the stabiliser $ \text{O}(L)_S $. 
The automorphism $ \phi $ of $ L $ defined as $ -1 $ on $ \langle e_1, e_1' \rangle $ and as the identity on the orthogonal complement clearly preserves $ S $. Let $ z = (e_2 - e_2' / m_2) + i (e_1 - e_1' / m_1) \in \mathcal{K} $. Then  $ \phi$ maps $ [z] $ to $ [\bar{z}] $. Hence the automorphism $ \phi $ does not preserve $ \mathcal{H} $. 

Next we turn to the discriminant kernel and choose a decomposition $ L = M \oplus \la e_1, e_1' \ra $ where $ e_1, e_1' $ are isotropic with $ (e_1, e_1') = 1 $. Then the reflection $ \sigma_v $ in $ v = e_1 - e_1' $ is an automorphism of $ L $ that acts trivially on $ M $. Hence it lies in the discriminant kernel. Now $ (v,v) = -2 $ implies that $ \sigma_v $ has negative spinor norm so that $ \sigma_v $ is not contained in $ \text{O}(L)^+ $. \eop

\medskip 

In particular we have $\text{O}(L)S = \text{O}(L)^+ S$ and $\ov{\Gamma} = \ov{\Gamma^+}$. In order to describe the orbits of $\Gamma^+$ on $\text{O}(L)S$ we define the map
\begin{align*}
	j : \Gamma^+ \backslash \text{O}(L)S &\to \ov{\Gamma} \backslash \text{O}(D) / \ov{\text{O}(L)_S^+}  \\
	\Gamma^+ \phi(S)      &\mapsto \ov{\Gamma} \, \ov{\phi} \, \ov{\text{O}(L)_S^+}
\end{align*}
where we choose $ \phi $ in $ \text{O}(L)^+ $ and denote by $\ov{\vphantom{l}\, . \,}$ the projection to $\text{O}(D)$. 

\begin{prp}
The map $j$ is a bijection.
\end{prp}
{\em Proof:}
By Theorem 1.14.2 in \cite{N} the natural map $\text{O}(L) \to \text{O}(D)$ is surjective. Proposition \ref{rb} implies that this map stays surjective when restricted to $\text{O}(L)^+$. It follows that $ j $ is surjective. 

Now let $ \ov{\Gamma} \, \ov{\phi} \, \ov{\text{O}(L)_S^+} = \ov{\Gamma} \, \ov{\psi} \, \ov{\text{O}(L)_S^+} $ with $ \phi, \psi \in \text{O}(L)^+ $. Then $ \ov{\Gamma^+} \, \ov{\phi} \, \ov{\text{O}(L)_S^+} = \ov{\Gamma^+} \, \ov{\psi} \, \ov{\text{O}(L)_S^+} $ so that $ \Delta(L)^+ \Gamma^+ \phi(S) = \Delta(L)^+ \Gamma^+ \psi(S) $. Since $ \Gamma $ contains $ \Delta(L) $, we can omit $ \Delta(L)^+ $ in this equation and obtain $ \Gamma^+ \phi(S) =  \Gamma^+ \psi(S) $.
\eop 

\medskip 

In general the group $\ov{\text{O}(L)_S^+}$ is difficult to describe. We will see that $\ov{\text{O}(L)_S^+} = \ov{\text{O}(K)}$ if $S$ has type $0$. 

\medskip

Let ${\cal C}$ be a $1$-dimensional cusp of $\Gamma^+ \backslash {\cal H}$. Choose a $2$-dimensional isotropic subspace $U \subset V$ representing ${\cal C}$. We define the {\em type} of ${\cal C}$ as the orbit $\text{O}(D)H$ where $H=(L'\cap U)/(L \cap U)$ and the {\em lattice associated with ${\cal C}$} as the isomorphism class of the lattice $K$ associated with $U$. Sometimes we identify the type $\text{O}(D)H$ with a representative of this orbit. 

Suppose ${\cal C}$ or equivalently $U$ splits. Define $S=L \cap U$. Our considerations above give the following parametrisation of the splitting cusps of $\Gamma^+ \backslash {\cal H}$ which have the same type and associated lattice as ${\cal C}$.

\begin{thm} \label{classcuspsgen}
The $1$-dimensional splitting cusps of $\Gamma^+ \backslash {\cal H}$ of type $\text{O}(D)H$ with associated lattice $K$ are represented by the primitive isotropic planes $\phi(S)$ where $\phi \in \text{O}(L)$ ranges over a set of representatives of the double cosets $\ov{\Gamma} \backslash \text{O}(D) / \ov{\text{O}(L)_S^+}. $
\end{thm}

Finally we consider cusps of type $0$. Let ${\cal C}$ be a $1$-dimensional cusp of $\Gamma^+ \backslash {\cal H}$ of type $0$ represented by $S \subset L$ and $K$ its associated lattice. Since $H$ is trivial, we can write $L = L^H = K \oplus I\!I_{1,1} \oplus I\!I_{1,1}$ with $S \subset K^{\perp}$ (cf.\ Proposition \ref{dt}). 

\begin{prp}
The cusp ${\cal C}$ splits.
\end{prp}

The decomposition of $L$ defines an embedding

\[  \text{O}(K) \hookrightarrow \text{O}(L) \to \text{O}(D)  \, . \]

\begin{prp} \label{stab0}
The image of $\text{O}(L)_S^+$ in $\text{O}(D)$ is given by $\ov{\text{O}(L)_S^+} = \ov{\text{O}(K)}$.  
\end{prp}
{\em Proof:}
The map $ \phi $ constructed in the proof of Proposition \ref{rb} shows that $ \ov{\text{O}(L)_S^+} = \ov{\text{O}(L)_S} $. Now $\text{O}(K) \subset \text{O}(L)_S$ implies
$\ov{\text{O}(K)} \subset \ov{\text{O}(L)_S^+}$.
Next let $ \sigma \in \text{O}(L)_S $. Then $\sigma(S^{\perp}) = S^{\perp} = K \oplus S$ so that $\sigma(K) \subset K \oplus S$. Let $\pi : V \to K \otimes_{\Z} \Q$ be the orthogonal projection. The restriction $\pi|_{K \oplus S}$ preserves the inner product $(\, , \, )$ and therefore is injective on $\sigma(K)$. The lattice $(\pi \circ \sigma)(K)$ is a sublattice of $K$ isomorphic to $K$. Hence $(\pi \circ \sigma)(K) = K$ so that $\pi \circ \sigma \in \text{O}(K)$ and $\ov{\pi \circ \sigma} \in \ov{\text{O}(K)}$. The decomposition $L = K \oplus I\!I_{1,1} \oplus I\!I_{1,1}$ implies $ \ov{\pi \circ \sigma} = \ov{\sigma} $ and $\ov{\sigma} \in \ov{\text{O}(K)}$. \eop 

\medskip 

We obtain the following description of the cusps of type $0$ in $\Gamma^+ \backslash \cal{H}$.
\begin{thm} \label{classcusps}
For each isomorphism class in the genus $I\!I_{n-2,0}(D)$ choose a primitive representative $K \subset L$ and a primitive isotropic plane $S_K \subset K^{\perp}$. Then the $1$-di\-men\-sional cusps of $\Gamma^+\backslash \cal{H}$ of type $0$ are represented by the primitive isotropic planes $\phi_K(S_K)$ where 
\begin{enumerate}[i)]
\item $K$ ranges over the isomorphism classes of lattices in the genus $I\!I_{n-2,0}(D)$ and
\item for each such $K$, the automorphism $\phi_K \in \text{O}(L)$ ranges over a set of representatives of the double cosets $\ov{\Gamma} \backslash \text{O}(D) / \ov{\text{O}(K)}$.
\end{enumerate} 
\end{thm}

We remark that we can use the same approach to study $0$-dimensional cusps. If $L$ splits two unimodular hyperbolic planes, we find that $\Gamma^+ \backslash \cal{H}$ has a unique $0$-dimensional cusp of type $0$.

\medskip

Our results generalise those of Attwell-Duval \cite{At14, At16} and Kiefer \cite{Ki} who considered splitting cusps for rescaled maximal lattices and for arbitrary lattices under the assumption $\Gamma = \Delta(L)$. Note that in this case the double quotient in Theorem \ref{classcuspsgen} reduces to a quotient.

\subsection*{Reflective modular varieties}

We now consider the modular varieties $\text{O}(L,F)^+\backslash {\cal H}$ corresponding to the $11$ re\-flec\-tive automorphic products $\psi_F$ constructed in Theorem \ref{blazefoley}. We associate to a $1$-dimensional cusp ${\cal C}$ of $\text{O}(L,F)^+\backslash {\cal H}$ a set $R_{\cal C}$ which is determined by ${\cal C}$ and the singular coefficients of $F$. If ${\cal C}$ has type $0$, the set $R_{\cal C}$ is either empty or a scaled root system.
Using the parametrisation by double quotients given in Theorem \ref{classcusps} we show that the root systems which occur are exactly those from Schellekens' classification of meromorphic conformal field theories of central charge $24$. If $R_{\cal C}$ is non-empty, $\psi_F$ vanishes identically along ${\cal C}$. Since $\psi_F$ has singular weight, $\text{O}(L,F)^+\backslash {\cal H}$ possesses a $1$-dimensional cusp on which $ \psi_F $ is a non-zero modular form. We construct a cusp on which the restriction of $\psi_F$ is the eta product of the class in $\text{Co}_0$ corresponding to $\psi_F$.

\medskip

Let $L$ be an even lattice of signature $(n,2)$, $n>2$ and even, splitting $I\!I_{1,1} \oplus I\!I_{1,1}$ and $V = L \otimes_{\Z} \Q$. Suppose $L$ carries a reflective automorphic product $\psi_F$ of singular weight. We define $\Gamma = \text{O}(L,F)$ and denote the Fourier coefficients of $F_{\gamma}$, $\gamma \in D$ again by $c_{\gamma}(m)$. Let $U \subset V$ be a $2$-dimensional isotropic subspace. Then the quotient $(L' \cap U^{\perp})/(L \cap U)$ possibly has torsion.
We define $R_U \subset L'$ as a set of representatives of

\[  \{ \alpha \in (L' \cap U^{\perp})/(L \cap U) \, | \, (\al, \al) > 0, \, c_{\al}(-\al^2/2) = 1  \} \, .  \]
Different choices of $R_U$ are isometric. Note that $\psi_F$ vanishes on $\al^{\perp}$ for $\al \in R_U$. If $\phi \in \Gamma$, then $\phi(R_{U})$ is a possible choice for $R_{\phi(U)}$, i.e.\ the geometry of the set $R_U$ is invariant under $\Gamma$. We define $R_{{\cal C}(U)} = R_U$ for the $1$-dimensional cusp $\mathcal{C}(U)$ of $ \Gamma^+ \backslash \mathcal{H} $. The relevance of $R_{{\cal C}(U)}$ stems from the fact that together with the type of ${\cal C}(U)$ it determines the first coefficient $\psi_0$ in the sum expansion of $\psi_F$ at the cusp ${\cal C}(U)$ (see Theorem \ref{kudlaprod}).

Let $H$ be the type of $U$. We show that $R_U$ is a finite subset of ${L^H}' \subset L'$. Choose a decomposition $L^H = K \oplus \langle e_1, e_1' \rangle \oplus \langle e_2, e_2' \rangle $ where $K$ is the lattice associated with $U$ as in Proposition \ref{dt}. 

\begin{prp} \label{intersectionWithLDual}
  We can write the elements of $R_U$ as $\al + u$ where $\al$ ranges over $K'\backslash\{0\}$ and $u$ over a set of representatives of $H = (L' \cap U)/(L \cap U)$
such that $c_{\al + u}(-\al^2/2) = 1$. In particular $R_U \subset {L^H}'$ and $R_U$ is finite.
\end{prp}
{\em Proof:}
We have
\begin{align*}
	L' \cap U^{\perp}
	&= \{x \in V \, | \, (x,U) = 0 , \, (x,L) \subset \Z \} \\
	&= \{x \in V \, | \, (x,U) = 0 , \, (x,L^H) \subset \Z\} \\
	&= {L^H}' \cap U^{\perp}   
\end{align*}
because $L^H = (L' \cap U) + L$. The vectors $e_1,e_2$ span $U$ so that
\[  L' \cap U^{\perp} = {L^H}' \cap U^{\perp} = K' + (L' \cap U)  \, . \]
This implies the statement.  \eop

\medskip

Now we specialise to $1$-dimensional cusps of type $0$. In this case the sets $R_U$ are often root systems. So suppose $U$ has type $H = (L' \cap U)/(L \cap U) = 0$. Then $R_U$ can be chosen as
\[   R_U = \{\al \in K' \backslash \{ 0 \} \, | \, c_{\alpha}(-\alpha^2/2) = 1 \} \subset K'  \]
and
$L = K \oplus \la e_1, e_1' \ra \oplus \la e_2, e_2' \ra$.

\begin{prp} \label{ReflInAutGroup}
Let $\al \in V $. If $\psi_F$ vanishes on $\al^{\perp}$, then the reflection $\sigma_{\al}$ is in $\Gamma$.
\end{prp}
{\em Proof:}
We can assume that $\al$ is a root of $L$ because $\psi_F$ is reflective. Then $ \sigma_{\al} \in \text{O}(L)$. It suffices to show that the map on $D$ induced by $\sigma_{\al}$ preserves the principal part of $F$. Let $\bt \in L'$ and $x \in \Z - \bt^2/2 $ with $x < 0$. We have to show that $c_{\bt}(x) = c_{\sigma_{\al}(\bt)}(x)$. Since we can replace $\bt$ by any element in $\bt +L$ and $L$ splits a hyperbolic plane $I\!I_{1,1}$, we can assume that $\bt$ is primitive in $L'$ and has positive norm $\bt^2/2=-x$. Then the order of the rational quadratic divisor $\bt^{\perp}$ is given by
\[ \sum_{\substack{k \in \Q_{>0} \\k \bt \in L'}} c_{k \bt}(-k^2\bt^2/2)
= \sum_{k \in \Z_{>0}} c_{k \bt}(k^2x) = c_{\bt}(x) + c_{2\bt}(4x)
\]
(cf.\ Proposition \ref{mbig}). Analogously we find for the order of $\sigma_{\al}(\bt)^{\perp}$
\[ \sum_{k \in \Z_{>0}} c_{k \sigma_{\al}(\bt)}(k^2x) = c_{\sigma_{\al}(\bt)}(x) + c_{2\sigma_{\al}(\bt)}(4x) \, .
\]
Since $\psi_F$ vanishes on $\al^{\perp}$, the automorphic form $\psi_F$ transforms with $-1$ under the reflection in $\al$ (see Theorem 1.2 in \cite{WW23}). Hence the orders of the divisors $\bt^{\perp}$ and $\sigma_{\al}(\bt)^{\perp}$ agree, i.e.\
\[   c_{\bt}(x) + c_{2\bt}(4x) = c_{\sigma_\al(\bt)}(x) + c_{2\sigma_\al(\bt)}(4x) \, . \]
Next we choose a primitive representative of $2\bt + L$ in $L'$ of norm $-4x$. Repeating the above argument we obtain
\[   c_{2\bt}(4x) + c_{4\bt}(16x) = c_{2\sigma_{\al}(\bt)}(4x) + c_{4\sigma_{\al}(\bt)}(16x) \, . \]
We have already seen that $c_{4\bt}(16x) = c_{4\sigma_{\al}(\bt)}(16x) = 0$. This implies $c_{2\bt}(4x) = c_{2\sigma_{\al}(\bt)}(4x)$ so that
$c_{\bt}(x) = c_{\sigma_{\al}(\bt)}(x)$. \eop

\medskip 

It follows that the reflections $\sigma_{\al}$, $\al \in R_U$ preserve $R_U$. We will see that under weak assumptions $R_U$ is actually a root system. 

\begin{prp}\label{rootSys}
If $c_{2\alpha}(-4\alpha^2/2) = 0$ for all $ \alpha \in R_U $, then the set $R_{U}$ is either empty or a root system in $K'$.  
\end{prp}
{\em Proof:}
Suppose $R_{U}$ is non-empty. We show that $R_{U}$ is a root system.

By Proposition \ref{ReflInAutGroup} the set $R_{U}$ is invariant under the reflections $\sigma_{\al}$, $\al \in R_U$.

Let $\al \in R_{U}$. We verify that the only rational multiples of $\al$ in $R_{U}$ are $\pm \al$. The transformation behaviour of $F$ under $-1 \in \text{SL}_2(\Z)$ implies $F_{\al} = F_{-\al}$ so that $-\al \in R_{U}$. Suppose $x \al \in R_{U}$ for some $x\in \Q_{>0}$. Write $x = m/d$ with coprime positive integers $m,d$. Choose $a,b \in \Z$ such that $ad+bm=1$. Then $a\al + bx\al = (ad+bm)\al/d = \al/d$ so that $\bt = \al/d \in L'$. Note that $\bt + e_1$ is primitive in $L'$. The divisor $(\bt + e_1)^{\perp}$ has order
\[
  \sum_{k \in \Z_{>0}} c_{k(\bt +e_1)}(-k^2(\bt+e_1)^2/2)
  = c_{\bt}(-\bt^2/2) + c_{2\bt}(-4\bt^2/2) 
\]
(cf.\ Proposition \ref{mbig}). Since $c_{\al}(-\al^2/2) = c_{d\bt}(-d^2\bt^2/2)$ and $c_{x\al}(-x^2\al^2/2) = c_{m\bt}(-m^2\bt^2/2)$ are both $1$, we have $x = 1/2$, $1$ or $2$. Now
$c_{2\alpha}(-4\alpha^2/2) = c_{2x\alpha}(-4x^2\alpha^2/2) = 0$ implies $x=1$.

Let $\al, \bt \in R_{U}$. We show that $2(\al,\bt)/\bt^2$ is an integer. The divisor $(\bt + e_1)^{\perp}$ has order
\[  \sum_{k \in \Z_{>0}} c_{k(\bt +e_1)}(-k^2(\bt+e_1)^2/2) = c_{\bt}(-\bt^2/2) + c_{2\bt}(-4\bt^2/2) = 1 \, . \]
By the reflectivity of $\psi_F$ some multiple of $\bt + e_1$ is a root of $L$ and hence $\sigma_{\bt + e_1} \in \text{O}(L) = \text{O}(L')$. It follows
\[  \sigma_{\bt + e_1}(\al + e_1) = \sigma_{\bt}(\al) + (1 - 2(\al,\bt)/\bt^2) e_1 \in L' \, . \]
Hence $2(\al,\bt)/\bt^2 \in \Z$ because $e_1$ is primitive in $L'$.

Finally we prove that $R_{U}$ spans $K \otimes_{\Z} \R$. By Theorem 10.5 in \cite{B98} the numbers $c_{\al}(-\al^2/2)$ form a vector system in $K \otimes_{\Z} \R$. This implies that the function $\lambda \mapsto \sum_{\al \in K'} c_{\al}(-\al^2/2)(\lambda,\al)^2$ is constant on $\{ \lambda \in K \otimes_{\Z} \R \, | \, \lambda^2 = 1 \}$. Hence $ R_U $ spans $K \otimes_{\Z} \R$. \eop

\medskip

We remark that a similar result was proved by Wang (cf.\ Theorem 2.2 in \cite{W23}).

\medskip

We denote the irreducible components of $R_{U} = R_{{\cal C}(U)}$ of type $X_m$ with long roots of norm $2/k$ by $X_{m,k}$. Now we determine the root systems $R_{{\cal C}(U)}$ for the $11$ automorphic products $\psi_F$ constructed in Section \ref{DM}. We remark that the condition $ c_{2\alpha}(-4\alpha^2/2) = 0 $ for $ \alpha \in R_U $ is satisfied here since these automorphic products only have simple zeros.

\begin{thm} \label{heaven17}
  Let $\psi_F$ be one of the $11$ reflective automorphic products of singular weight described in Theorem \ref{blazefoley}. Then the root systems $R_{{\cal C}(U)}$ of the $1$-dimensional cusps ${\cal C}(U)$ of type $ 0 $ are given in the following table:
\vspace*{2mm}
\[
\renewcommand{\arraystretch}{1.2}
\begin{array}{c|lll}
		L & \multicolumn{3}{c}{R_{{\cal C}(U)}} \\ \hline
		& & & \\ [-4.5mm]
		I\!I_{26,2} & D_{24,1} & D_{16,1}E_{8,1} & E_{8,1}^3 \\ 
		& A_{24,1} & D_{12,1}^2 & A_{17,1} E_{7,1} \\ 
		& D_{10,1}E_{7,1}^2 & A_{15,1}D_{9,1} & D_{8,1}^3 \\ 
		& A_{12,1}^2 & A_{11,1}D_{7,1}E_{6,1} & E_{6,1}^4 \\ 
		& A_{9,1}^2 D_{6,1} & D_{6,1}^4 & A_{8,1}^3 \\ 
		& A_{7,1}^2D_{5,1}^2 & A_{6,1}^4 & A_{5,1}^4D_{4,1} \\ 
		& D_{4,1}^6 & A_{4,1}^6 & A_{3,1}^8 \\ 
		& A_{2,1}^{12} & A_{1,1}^{24} & \{\} \\[0.3mm]  \hline
		& & & \\ [-4.5mm]
		I\!I_{18,2}(2_{I\!I}^{+10}) & B_{8,1}E_{8,2} & B_{6,1} C_{10,1} &  C_{8,1} F_{4,1}^2 \\
		& B_{5,1}E_{7,2}F_{4,1} & A_{7,1} D_{9,2} & B_{4,1}^2 D_{8,2} \\
		& B_{4,1} C_{6,1}^2 & A_{5,1}C_{5,1} E_{6,2} &  A_{4,1} A_{9,2} B_{3,1} \\
		& B_{3,1}^2 C_{4,1} D_{6,2} & C_{4,1}^4 & A_{3,1} A_{7,2} C_{3,1}^2 \\
		& A_{3,1}^2 D_{5,2}^2 & A_{2,1}^2 A_{5,2}^2 B_{2,1} & B_{2,1}^4 D_{4,2}^2 \\
		& A_{1,1}^4 A_{3,2}^4 & A_{1,2}^{16} \\[0.3mm]  \hline
		& & & \\ [-4.5mm]
  I\!I_{14,2}(3^{-8}) & A_{5,1} E_{7,3} & A_{3,1} D_{7,3} G_{2,1} & E_{6,3} G_{2,1}^3 \\
                  & A_{2,1}^2 A_{8,3} & A_{1,1}^3 A_{5,3} D_{4,3} & A_{2,3}^6 \\[0.3mm]  \hline
		& & & \\ [-4.5mm]
  I\!I_{14,2}(2_{I\!I}^{-10} 4_{I\!I}^{-2}) & B_{12,2} & B_{6,2}^2 & B_{4,2}^3 \\  & B_{3,2}^4 & B_{2,2}^6 & A_{1,4}^{12} \\
                  & A_{8,2} F_{4,2} & A_{4,2}^2 C_{4,2} & A_{2,2}^4 D_{4,4} \\[0.3mm]  \hline
		& & & \\ [-4.5mm]
  I\!I_{12,2}(2_2^{+2} 4_{I\!I}^{+6}) & A_{3,1} C_{7,2} & A_{2,1} B_{2,1} E_{6,4} & A_{1,1}^3 A_{7,4} \\
                  & A_{1,1}^2 C_{3,2} D_{5,4} & A_{1,2} A_{3,4}^3 \\[0.3mm]  \hline
		& & & \\ [-4.5mm]
		I\!I_{10,2}(5^{+6}) & A_{1,1}^2 D_{6,5} & A_{4,5}^2 & \\[0.3mm]  \hline
		& & & \\ [-4.5mm]
		I\!I_{10,2}(2_{I\!I}^{+6} 3^{-6}) & A_{1,1} C_{5,3} G_{2,2} & A_{1,2} A_{5,6} B_{2,3} & \\[0.3mm]  \hline
		& & & \\ [-4.5mm]
		I\!I_{8,2}(7^{-5}) & A_{6,7} & & \\[0.3mm]  \hline
		& & & \\ [-4.5mm]
		I\!I_{8,2}(2_7^{+1} 4_7^{+1} 8_{I\!I}^{+4}) & A_{1,2} D_{5,8} & & \\[0.3mm]  \hline
		& & & \\ [-4.5mm]
		I\!I_{8,2}(2_{I\!I}^{+4} 4_{I\!I}^{-2} 3^{+5}) & A_{2,2} F_{4,6} & A_{2,6} D_{4,12} & \\[0.3mm]  \hline
		& & & \\ [-4.5mm]
		I\!I_{6,2}(2_{I\!I}^{-2} 4_{I\!I}^{-2} 5^{+4}) & C_{4,10} & &  
	\end{array}
	\]
      \end{thm}
      \vspace*{2mm}
{\em Proof:}
Let $L$ be the lattice corresponding to $\psi_F$, $n = \dim(L)$ and $D = L'/L$. The principal part of $F$ and in particular the set $M = \{ \gamma \in D \, | \, \text{$F_{\gamma}$ is singular} \}$ is described in Section \ref{DM}. We apply Theorem \ref{classcusps} to determine the cusps of type $ 0 $ of $\Gamma^+\backslash \cal{H}$ where $\Gamma = \text{O}(L,F)$. Let $K$ be a lattice in the genus $I\!I_{n-2,0}(D)$. Decompose $L = K \oplus K^{\perp}$ and choose a primitive isotropic plane $S$ in $K^{\perp}$. Let $\phi \in \text{O}(L)$. Then $L = \phi(K) \oplus \phi(K)^{\perp}$ and for $U = \phi(S) \otimes_{\Z} \Q$ the root system $R_{{\cal C}(U)}$ is given by
\begin{align*}
  R_{{\cal C}(U)}
  &= \{ \al \in \phi(K)' \backslash \{ 0 \} \, | \, \al^2/2 \leq 1 \text{ and } \al + \phi(K) \in M \}   \\
  &= \phi \big( \{ \al \in K' \backslash \{ 0 \} \, | \, \al^2/2 \leq 1 \text{ and } \al + K \in \phi^{-1}(M) \} \big)  \\
  & \simeq \{ \al \in K' \backslash \{ 0 \} \, | \, \al^2/2 \leq 1 \text{ and } \al + K \in \phi^{-1}(M) \} \, .
\end{align*}
Now we let $\phi$ range over representatives of $\ov{\Gamma} \backslash \text{O}(D) / \ov{\text{O}(K)}$ and $K$ over the lattices in the genus $I\!I_{n-2,0}(D)$ to obtain the root systems corresponding to the cusps of type $ 0 $ for this case. Note that $\ov{\Gamma} = \text{O}(D)_M$, the stabiliser of the set $M$ in $\text{O}(D)$.

The necessary computations were performed with the computer algebra system Magma \cite{BCP}. If $\Gamma = \text{O}(L)$, the quotient $ \ov{\Gamma} \backslash \text{O}(D) $ is trivial for each $ K $ and we can choose $ \phi $ as the identity. For the remaining cases $ I\!I_{14,2}(2_{I\!I}^{-10} 4_{I\!I}^{-2}) $, $ I\!I_{8,2}(2_{I\!I}^{+4} 4_{I\!I}^{-2} 3^{+5}) $ and $ I\!I_{6,2}(2_{I\!I}^{-2} 4_{I\!I}^{-2} 5^{+4}) $ this quotient is non-trivial. We give more details in the table below:
\vspace*{2mm}
\[
\renewcommand{\arraystretch}{1.2}
\begin{array}{c|c|c|cc}
	L & \vert \ov{\Gamma} \backslash \text{O}(D) \vert & K & \multicolumn{2}{c}{R_{{\cal C}(U)}} \\ [0.3mm] \hline
	& & & & \\ [-4.5mm]
	I\!I_{14,2}(2_{I\!I}^{-10} 4_{I\!I}^{-2}) & 2^{11} \cdot 3 \cdot 17 & D_{12}(2) & B_{12,2} & B_{6,2}^2 \\ 
	& & & B_{4,2}^3 & B_{3,2}^4 \\ [0.3mm] 
	& & & B_{2,2}^6 & A_{1,4}^{12} \\ [0.3mm] \cline{3-5}
	& & & & \\ [-4.5mm]
	& & D_4(2) \oplus E_8(2) & A_{8,2} F_{4,2} & A_{4,2}^2 C_{4,2} \\ [0.3mm] 
	& & & A_{2,2}^4 D_{4,4} & \\ \hline
	& & & & \\ [-4.5mm]
	I\!I_{8,2}(2_{I\!I}^{+4} 4_{I\!I}^{-2} 3^{+5}) & 2 & A_2(2) \oplus D_4(6) & A_{2,2} F_{4,6} & A_{2,6} D_{4,12} \\ [0.3mm] \hline
	& & & & \\ [-4.5mm]
	I\!I_{6,2}(2_{I\!I}^{-2} 4_{I\!I}^{-2} 5^{+4}) & 3 & D_4(10) & C_{4,10} & 
\end{array}
\]

\vspace*{2mm}
\noindent
We describe the entries for $ L = I\!I_{14,2}(2_{I\!I}^{-10} 4_{I\!I}^{-2}) $. The group $\text{O}(D)$ can be generated by reflections and Eichler transformations. We find that the stabilizer $ \ov{\Gamma} $ of $ M $ has index $ 2^{11} \cdot 3 \cdot 17 = 104448 $ in $ \text{O}(D) $. There are two isomorphism classes in the genus $ I\!I_{12,0}(2_{I\!I}^{-10} 4_{I\!I}^{-2}) $ represented by $ D_{12}(2) $ and $ E_8(2) \oplus D_4(2) $. For $ K = D_{12}(2) $ and $ K = E_8(2) \oplus D_4(2) $ the double quotient $\ov{\Gamma} \backslash \text{O}(D) / \ov{\text{O}(K)}$ has $ 6 $ and $ 3 $ elements, respectively, and we list the corresponding root systems in the last column. \eop  

\medskip

We see that no root system occurs twice so that the $1$-dimensional cusps of type $0$ are parametrised by their root system. Furthermore these root systems are exactly the same as those found by Schellekens in his classification of holomorphic vertex operator algebras of central charge $24$ (see \cite{ANS}, \cite{EMS}). We will explain this in the next section.

\begin{thm} \label{onedimexp}
	Let $\psi_F$ be one of the $11$ reflective automorphic products of singular weight given in Theorem \ref{blazefoley} and ${\cal C}(U)$ a $1$-dimensional cusp of type $0$. Then the expansion of $\psi_F$ at ${\cal C}(U)$ is given by
	\begin{align*}
		\MoveEqLeft q_1^{I_0} \sum_{m=0}^{\infty} \psi_m(w, \tau_2) q_1^m \\
		& = \, q_1^{I_0} \psi_0(w, \tau_2) \prod_{a \in \Z_{>0}} \, \prod_{b \in \Z} \, \prod_{\al \in K'}
		\big(1 - q_1^a q_2^b e(-(\al, w))\big)^{c_\al(-\al^2/2 + ab)}
	\end{align*}
	where $I_0 = \frac{n-2 + |R_{{\cal C}(U)}|}{24} - 1$ and
	\[ \psi_0(w, \tau_2) = \kappa \eta(\tau_2)^{n-2}
	\prod_{\al \in R_{{\cal C}(U)}^+} \frac{\vartheta(-(\al, w), \tau_2)}{\eta(\tau_2)}   \]
	for a set of positive roots $R_{{\cal C}(U)}^+ \subset R_{{\cal C}(U)} $ and a constant $\kappa$ of absolute value $1$.
\end{thm}
{\em Proof:}
We decompose $L = K \oplus \la e_1, e_1' \ra \oplus \la e_2, e_2' \ra$ with unimodular hyperbolic planes $\la e_i, e_i' \ra$ and $U = \Q e_1 + \Q e_2$ (see Proposition \ref{mdt}).
The expansion of $\psi_F$ at ${\cal C}(U)$ is described in Theorem \ref{kudlaprod}.
Since $L \cap \Q e_1 = L' \cap \Q e_1= \Z e_1$, the first product extends over the elements $\al + ae_2' - be_2 $ with $\al \in K'$, $a \in \Z_{>0}$ and $b \in \Z$ and is given by
\[  \prod_{a \in \Z_{>0}} \, \prod_{b \in \Z} \, \prod_{\al \in K'}
\big( 1 - q_1^a q_2^b e(-(\al, w)) \big)^{c_{\al}(-\al^2/2 + ab)}  \, . \]
The set $C$ specialises to a Weyl chamber of the root system $R_{{\cal C}(U)}$ (see Section 2.3 in \cite{K}) so that the second product reduces to
\[    \prod_{\al \in R_{{\cal C}(U)}^+} \frac{\vartheta(-(\al, w), \tau_2)}{\eta(\tau_2)}   \]
where $R_{{\cal C}(U)}^+$ is the set of positive roots in $R_{{\cal C}(U)}$ corresponding to $C$.
Since $(L' \cap U) \!\! \mod (L \cap U) = \{ 0 \}$, the third product is $1$.
Next we determine $I_0$. Using $(L' \cap U^{\perp}) \!\! \mod (L \cap U) = K'$ we obtain
\begin{align*}
	I_0 
	& = -\sum_{m \in \Q} \, \sum_{\al \in K'} c_{\al}(-m) \sigma_1(m - \al^2/2) \\
	& = - \sum_{m \in \Q} c_{0}(-m)\sigma_1(m) - \sum_{\al \in R_{{\cal C}(U)}} c_{\al}(-\al^2/2) \sigma_1(0) \\
	& = \frac{ n-2 + |R_{{\cal C}(U)}|}{24} - 1  \, . 
\end{align*}
The theorem now follows from the fact that $q_1^{I_0} \psi_0(w, \tau_2)$ is the product of the second, third and fourth factor in Theorem \ref{kudlaprod}. \eop

\medskip
\noindent
The constant coefficient $\psi_0$ in the expansion of $\psi_F$ at ${\cal C}(U)$ is essentially the denominator function of the affine Kac-Moody algebra associated to $R_{{\cal C}(U)}$.

\medskip

The case where $\psi_F$ is the theta lift of $1/\Delta$ on the unimodular lattice $I\!I_{26,2}$ was already studied by Gritsenko in \cite{G12}. 

\medskip

The automorphic products $\psi_F$ that we consider all have singular weight and therefore cannot be cusp forms. However with one exception they all vanish at the $1$-dimensional cusps of type $0$. In the following we will construct for each of the $11$ automorphic products $\psi_F$ a special $1$-dimensional cusp along which $\psi_F$ is given by the associated eta product $\eta_g$. Of course for $\psi_F$ the theta lift of $1/\Delta$ this cusp is the type-$0$ cusp with associated lattice the Leech lattice.   

\medskip

Let $\psi_F$ be one of the $11$ reflective automorphic products of singular weight given in Theorem \ref{blazefoley} and $g$ an element in the corresponding class in $\text{Co}_0$. Recall that
\[  L = \Lambda^g \oplus I\!I_{1,1} \oplus I\!I_{1,1}(m)   \]
if $N/m=1$ and
\[  L = \Lambda^g_N \oplus I\!I_{1,1} \oplus I\!I_{1,1}(m/2)   \]
if $N/m=2$.

\begin{prp}\label{egilstadir}
There exists a $2$-dimensional isotropic subspace $U$ of $L$ of type $H$ with associated lattice $\Lambda^g$ where $H \subset D^{N/m}$ is a cyclic isotropic subgroup of order $m$. 
\end{prp}
{\em Proof:}
First we consider the case $N/m = 1$. Then $L = \Lambda^g \oplus I\!I_{1,1} \oplus I\!I_{1,1}(m)$. Choose isotropic bases $(e_1,e_1')$ and $(e_2,e_2')$ of $I\!I_{1,1}$ and $I\!I_{1,1}(m)$, respectively, such that $(e_1,e_1') = 1$ and $(e_2,e_2') = m$ and define $U = \Q e_1 + \Q e_2$. Then
\[  H = (L' \cap U)/(L \cap U) = \Z (e_2/m)/\Z e_2 \simeq \Z/m\Z   \]
and
\[  L^H = \bigcup_{j = 0}^{m-1} (j(e_2/m) + L) = \Lambda^g \oplus \la e_1, e_1' \ra \oplus \la e_2/m, e_2' \ra \]
where $ \la e_2/m, e_2' \ra \simeq I\!I_{1,1} $. 

The case $N/m = 2$ is more complicated. Here $L = \Lambda^g_N \oplus I\!I_{1,1} \oplus I\!I_{1,1}(m/2)$. As above choose isotropic bases $(e_1,e_1')$ and $(e_2,e_2')$ of $I\!I_{1,1}$ and $I\!I_{1,1}(m/2)$, respectively, with $(e_1,e_1') = 1$ and $(e_2,e_2') = m/2$. There is a primitive element $x \in \Lambda^g \cap m{\Lambda^g}'$ such that
\[  \Lambda^g = \Lambda^g_N \cup (x + \Lambda^g_N)  \, .  \]
Define $U = \Q f_1 + \Q f_2$
with $f_1 = x + e_1 + a e_1'$, $a = -x^2/2$ and $f_2 = e_2$. Then
\[H = (L' \cap U)/(L \cap U) = \big( \Z f_1 + \Z (f_2/(m/2)) \big) \big/ \big( \Z (2 f_1) + \Z f_2 \big) \simeq \Z / m\Z \]
because $m/2$ is odd. Since $H$ is generated by the elements $x + L$ and $e_2 / (m/2) + L$, we have $H \subset D^2$. The lattice $L^H$ is given by
\[  L^H = \Lambda^g \oplus \la e_1, e_1' \ra \oplus \la e_2/(m/2), e_2' \ra \]
where $ \la e_2/(m/2), e_2' \ra \simeq I\!I_{1,1} $. However this decomposition does not yet determine the lattice associated with $U$ because $U \cap L^H$ is not orthogonal to the $\Lambda^g$ in the sum. We slightly modify the decomposition. Let $(b_1,\ldots,b_n)$ be a basis of $\Lambda^g$ with $x = b_1$ and define $K = \la b_1-(b_1,x)e_1',\, \ldots, \, b_n-(b_n,x)e_1' \ra$. Then $K$ is a sublattice of $L^H$ isomorphic to $\Lambda^g$ and orthogonal to $\la f_1,e_1' \ra \oplus \la f_2 / (m/2), e_2' \ra \simeq I\!I_{1,1} \oplus I\!I_{1,1}$. Hence
\[  L^H = K \oplus \la f_1,e_1' \ra \oplus \la f_2 / (m/2), e_2' \ra \, . \]
Now $U \cap L^H = \la f_1, f_2 / (m/2) \ra \subset K^\perp$ implies that the lattice associated with $U$ is isomorphic to $K$ and hence to $\Lambda^g$. \eop

\medskip

\noindent
Now choose $U$ as constructed in the proof of the previous proposition. Note that $U$ splits for $N/m = 1$ but not for $N/m =2$. We also remark that the order $m$ of $H$ is the maximal order of a cyclic subgroup of $D^{N/m}$. 

\medskip

We will see that $R_U$ vanishes if $ N/m = 1 $. Analogously $ R_{\phi(U)} $ vanishes for some $\phi \in \text{O}(L)$ if $ N/m = 2 $. This is necessary for the expansion of $\psi_F$ at the corresponding cusp to have a non-zero constant term (cf.\ Corollary 4.6 in \cite{K}).

\medskip

We fix a decomposition
\[  L^H = \Lambda^g \oplus I\!I_{1,1} \oplus I\!I_{1,1}  \]
with $U \cap L^H$ orthogonal to $\Lambda^g$ and with $(N/m) \Lambda^g \subset L$ (see the proof of Proposition \ref{egilstadir}). Then the composition
\[  {L^H}' \hookrightarrow L' \to L'/L  \]
maps ${\Lambda^g}' \subset {L^H}'$ to $H^{\perp}$ and $\Lambda^g$ to $H \cap D_{N/m}$. The projection ${\Lambda^g}' \to H^{\perp} \to H^{\perp}/H$
has kernel $\Lambda^g$ and therefore
defines an isomorphism ${\Lambda^g}'/\Lambda^g \to H^{\perp}/H$.

Recall that 
\[  R_U = \{ \al + u \, | \, \al \in {\Lambda^g}' \backslash\{0\}, \,  u \in H, \, c_{\al + u}(-\al^2/2) = 1 \} \subset {L^H}'  \, . \]
Here $u$ ranges over a set of representatives of $H$ (cf.\ Proposition \ref{intersectionWithLDual}).

The reflectivity of $F$ implies the following result.

\begin{prp} \label{jj}
If $\al + u \in R_U$ for $\al \in {\Lambda^g}'$ and $u \in L' \cap U$, then
$\al \in {\Lambda^g}' \cap (\Lambda^g/d)$ and $\al^2/2 = 1/d$ for some $d|N$.
\end{prp}
{\em Proof:}
Let $\al \in {\Lambda^g}'$ and $u \in L' \cap U$. Suppose $ \alpha + u \in R_U $ so $ c_{\alpha + u}(-\alpha^2/2) = 1 $. By the reflectivity of $ F $ this can only hold if $ (\alpha + u) + L \in D_{d,1/d} $ and $ \alpha^2/2 = 1/d $ for some positive divisor $ d $ of $ N $. Now $ (\alpha + u) + L \in D_{d} $ and $ u + L \in H $ imply $ d(\alpha + L) = 0 \!\! \mod H $. The isomorphism ${\Lambda^g}'/\Lambda^g \to H^{\perp}/H$ maps $ \alpha + \Lambda^g $ to $ (\alpha + L) + H $. Using $ (\alpha + L) + H \in (H^\perp/H)_d $ we deduce $ \alpha + \Lambda^g \in ({\Lambda^g}'/\Lambda^g)_d $ and $ \alpha \in \Lambda^g/d $ follows. \eop 

\medskip

For $d|N$ we define
\[  R_d = \{  \al + L \, | \, \al \in {\Lambda^g}' \cap  (\Lambda^g/d) \text{ and } \al^2/2 = 1/d  \, \} \subset H^\perp \cap D_{dN/m} \, . \]
We consider the image $R_U + L$ of $R_U$ in $L'/L$. Proposition \ref{jj} implies 
\[  R_U + L \subset \bigcup_{d|N} \big( (R_d + H) \cap M_d \big)   \]
where $M_d = \{\gamma \in D_{d,{1/d}} \, | \, F_{\gamma} \text{ singular} \}$. Hence $R_U = \{ \}$ if $(R_d + H) \cap M_d = \{ \}$ for all $d|N$. Using the computer algebra system Magma we can determine the sets $R_d$. We see that $R_d$ is empty for $d||N$. This implies that $R_U$ is empty if $N$ is squarefree. We consider now the remaining cases. 

\begin{prp}
  If $L = I\!I_{12,2}(2_2^{+2} 4_{I\!I}^{+6})$ or $L = I\!I_{8,2}(2_7^{+1} 4_7^{+1} 8_{I\!I}^{+4})  $, then $R_U = \{ \}$. 
\end{prp}
{\em Proof:}
Since the proofs are similar in both cases, we only describe the argument for $L = I\!I_{12,2}(2_2^{+2} 4_{I\!I}^{+6})$. In that case only $R_2$ is non-empty and consists of $10$ elements in $D^{2*}_0 $. The set $M_2$ is given by $D^{2*}_{1/2}$ (see Proposition \ref{uniquenessEI}). Let $\gamma \in R_2$ and $h \in H$ such that $\gamma + h \in D^{2*}$. Then we can write $h = 2h_4$ for a generator $h_4$ of $H$ and compute
\[  \q_2(\gamma + h) = \q_2(\gamma) + (\gamma, h_4) + 2q(h_4) = \q_2(\gamma)   \]
because $R_2 \subset H^{\perp}$ and $H$ is isotropic. It follows $\gamma + h \in D^{2*}_0$. This implies $(R_2 + H) \cap M_2 = (R_2 + H) \cap D^{2*}_{1/2} =  \{ \}$.
\eop

\medskip

Finally we consider the $3$ cases with $N/m =2$. Since $ R_d $ vanishes for $ d||N $, only $ R_2 $ and $ R_m $ may be non-empty. 

\begin{prp}\label{reykjavik}
For $d=2, m$ we have 
\[ (R_d + H) \cap M_d \subset (R_d \cap M_d) + H \, . \]
\end{prp}
{\em Proof:}
The description of the singular sets $M_d$ given in Section \ref{DM} (see the proofs of Proposition \ref{uniquenessJK} and \ref{uniquenessD}) implies
\[    ( \{ \gamma \} + H) \cap M_d \subset (\{\gamma\} \cap M_d) + H  \subset D   \]
for $\gamma \in H^{\perp} \cap D_d$. Since the elements of order $4$ in
$4_{I\!I}^{-2}$ have norm $1/4$ or $3/4 \! \mod 1$, the set $D_{2d} \backslash D_d $ contains no elements of norm $1/d \! \mod 1$. Hence $R_d \subset H^{\perp} \cap D_{2d}$ implies $R_d \subset H^{\perp} \cap D_d$. The statement now follows from the above inclusion.
\eop

\medskip

Proposition \ref{reykjavik} implies 
\[  R_U + L \subset \bigcup_{d = 2,m} \big( (R_d \cap M_d) + H \big). \]
Hence it would be enough to show $R_2 \cap M_2 = R_m \cap M_m = \{ \}$ in order to prove that $R_U$ vanishes. Since $ \psi_F $ for $ N/m = 2 $ is only unique up to $ \text{O}(L)^+ $, this will not be true in general for our specific choice of $ U $. In the following we construct an element in $\text{O}(L)^+$ such that the above intersections vanish after shifting $U$ by this automorphism.  

\begin{prp} \label{dcs}
There exists an element $\phi \in \text{O}(L)^+$ such that $\phi(M_d) \cap M_d = \{ \}$ for $d=2$ and $m$.
\end{prp}
{\em Proof:}
Since the projection $\text{O}(L)^+ \rightarrow \text{O}(D)$ is surjective,
it suffices to construct an automorphism of $D$ with the stated properties. We do this in a case-by-case analysis. 

Suppose $L = I\!I_{8,2}(2_{I\!I}^{+4} 4_{I\!I}^{-2} 3^{+5})$. Then $D_2/D^6 \simeq 2_{I\!I}^{-2} \oplus 2_{I\!I}^{-2}$. The singular set $M_2 \subset D_{2,1/2}$ generates one copy of the discriminant form $2_{I\!I}^{-2}$ in this quotient and $M_6 = (D_{2,1/2}\backslash M_2 + D_{3,2/3})$ (see the proof of Proposition \ref{uniquenessJK}). Hence the automorphism exchanging the two copies of $2_{I\!I}^{-2}$ in $D_2/D^6 $ gives the desired result. 

If $L = I\!I_{6,2}(2_{I\!I}^{-2} 4_{I\!I}^{-2} 5^{+4})$, then $M_2 = \mu + D^{10}$ for an element $\mu \in D_2$ of norm $1/2 \! \mod 1$ and $M_{10} = M_2 + D_{5,3/5}$.
Choosing an automorphism of $D$ permuting the non-zero elements of $ D_2/D^{10} \simeq 2_{I\!I}^{-2}$ yields the claim. 

Finally we consider the case $ L = I\!I_{14,2}(2_{I\!I}^{-10} 4_{I\!I}^{-2}) $. Let $ U $ be the projection of $ M_2 $ in $ D_2/D^2 \simeq 2_{I\!I}^{-10} $. We recall our description of $ U $ in the proof of the uniqueness of $ \psi_F $. There exists a basis $ (\gamma_1, \dots, \gamma_{10}) $ of $ D_2/D^2 $ with $ J_i = \la \gamma_{2i-1}, \gamma_{2i} \ra \simeq 2_{I\!I}^{-2} $ pairwise orthogonal such that the set $ \Gamma_i = J_i \backslash \{0\} $ is contained in $ U $ for $ i = 1, \dots, 4 $ and disjoint from $ U $ for $ i = 5 $ (cf.\ Proposition \ref{basis}). Let $ \mathcal{P}_2 $ be the set of partitions of $ \{1,2,3,4\} $ into 2-element subsets. Then there is a bijection $ \Phi : \Gamma_5 \rightarrow \mathcal{P}_2 $ such that the elements in $ U $ not contained in $ \bigcup_{i = 1, \dots, 4} \Gamma_i $ are precisely those in $ \mu + \Gamma_i + \Gamma_j $ where $ \{i,j\} \in \Phi(\mu) $ and $ \mu $ ranges over $ \Gamma_5 $ (cf.\ Proposition \ref{murtinheira}). 

We begin with the construction of $ \phi $. There are exactly eight elements $ \mu = \gamma_{i_1} + \gamma_{i_2} + \gamma_{i_3} $ with $ i_1 < i_2 < i_3 \leq 8 $ and $ i_1 = i_2 = i_3 \mod 2 $. We arrange these elements in a tuple $ (\mu_1, \mu_2, \dots, \mu_8) $ in such a way that $ (\mu_{2i-1}, J_i) = (\mu_{2i}, J_i) = 0 $. Let now $ \varphi : D_2/D^2 \to D_2/D^2 $ be the group homomorphism defined by $ \gamma_i \mapsto \mu_i $ on $ J_1 + J_2 + J_3 + J_4 = J_5^\perp $ and by a fixed-point free permutation of $ \Gamma_5 $ on $ J_5 $. We easily check that $ \varphi $ preserves the scalar product on $ D_2/D^2 $ and deduce that $ \varphi $ is an automorphism of $ D_2/D^2 $. Now choose an automorphism $ \phi \in \text{O}(D)$ inducing $ \varphi $. 

We show that $ \phi(M_2) \cap M_2 = \{\} $ or equivalently $ \varphi(U) \cap U = \{\} $. Suppose $\gamma \in U \cap \varphi(U)$. In order to derive a contradiction, we distinguish two cases. 

First we consider the case $ \gamma \in \Gamma_i $ for $ i \in \{1, \dots, 4\} $. Then $ \gamma \in J_1 + J_2 + J_3 + J_4 $ and hence $ \varphi^{-1}(\gamma) \in J_1 + J_2 + J_3 + J_4 $. Since by assumption $ \varphi^{-1}(\gamma) \in U $, we deduce $ \varphi^{-1} (\gamma) \in \Gamma_j $ for some $ j $ again satisfying $ 1 \leq j \leq 4 $. We find $ \gamma = \varphi(\varphi^{-1}(\gamma)) = \epsilon_{k_1} + \epsilon_{k_2} + \epsilon_{k_3} $ for non-zero elements $ \epsilon_{k} \in \Gamma_k $ and distinct indices $ k_1, k_2, k_3 $. But this clearly contradicts $ \gamma \in \Gamma_i $. 

Now suppose $ \gamma = \mu + \rho_i + \rho_j$ with $ \mu \in \Gamma_5 $ and $ \rho_i \in \Gamma_i $, $ \rho_j \in \Gamma_j $ where $ \{i,j\} \in \Phi(\mu) $. Then $\gamma$ is the image under $\varphi$ of some element of the same form, i.e.\ $ \gamma = \varphi(\lambda + \sigma_k + \sigma_l) $ with $ \lambda \in \Gamma_5 $ and $ \sigma_k\in \Gamma_k $, $ \sigma_l \in \Gamma_l $ where $ \{k,l\} \in \Phi(\lambda) $. 
Since $ \varphi $ preserves $ J_1 + J_2 + J_3 + J_4 $ as well as $ J_5 $, we obtain $ \varphi(\lambda) = \mu $ and $ \varphi(\sigma_k + \sigma_l) = \rho_i + \rho_j $. The images $ \varphi(\sigma_k) $ and $ \varphi(\sigma_l) $
are contained in $ \sum_{\nu \in \{1,\dots,4\}\backslash\{k\}} \, \Gamma_\nu $ and $ \sum_{\nu \in \{1,\dots,4\}\backslash\{l\}} \, \Gamma_\nu $. Note that if the projections of $ \varphi(\sigma_k) $ and $ \varphi(\sigma_l) $ agree on some $ J_\nu $
then the projections must agree for all $ \nu \neq k,l $.
Hence either
\[  \varphi(\sigma_k) + \varphi(\sigma_l) \in \Gamma_1 + \Gamma_2 + \Gamma_3 + \Gamma_4 \quad \text{or} \quad
    \varphi(\sigma_k) + \varphi(\sigma_l) \in \Gamma_k + \Gamma_l \, . \]
  Now $ \varphi(\sigma_k + \sigma_l) = \rho_i + \rho_j \in \Gamma_i + \Gamma_j $ implies $ \varphi(\sigma_k + \sigma_l) \in \Gamma_k + \Gamma_l $ so that $ \{k,l\} = \{i,j\} \in \Phi(\mu) $. Since $ \Phi $ is a bijection and $ \{k,l\} \in \Phi(\lambda) $, we conclude $ \lambda = \mu = \varphi(\lambda) $. Yet $ \varphi|_{\Gamma_5} $ is fixed-point free.
  \eop 

\medskip 

In Section \ref{DM} we have used the set
\[ M = \{\al + L \, | \, \al \in {\Lambda^g}' \cap (\Lambda^g/2), \ \al^2 = 1 \} + D^m \subset D  \]
to construct the reflective modular form on $L$. The uniqueness of $\psi_F$ up to $\text{O}(L)^+$ implies $ \psi(M) = M_2 $ for some $ \psi \in \text{O}(L)^+ $. Clearly $R_2 \subset M$. In summary this implies
\[   (\phi \circ \psi)(R_2) \cap M_2 \subset \phi(M_2) \cap M_2 = \{ \}    \]
for $ \phi $ as given in Proposition \ref{dcs}. The isotropic subspace $(\phi \circ \psi )(U)$ of $V$ has the same type and associated lattice as $U$.

\begin{prp}
We have $R_{(\phi \circ \psi)(U)}  = \{ \}$.
\end{prp}  
{\em Proof:}
We have seen above that $ (\phi \circ \psi)(R_2) \cap M_2 = \{\} $. It remains to consider $ R_6 $ for $ L = I\!I_{8,2}(2_{I\!I}^{+4} 4_{I\!I}^{-2} 3^{+5}) $ and $ R_{10} $ for $ L =  I\!I_{6,2}(2_{I\!I}^{-2} 4_{I\!I}^{-2} 5^{+4}) $. In the first case we use Magma to show that $ R_6 \subset D_{2,1/2}\backslash M + D_{3,2/3}$.
Then
\[  (\phi \circ \psi)(R_6) \subset \phi( D_{2,1/2}\backslash \psi(M) + D_{3,2/3} )
  = \phi( D_{2,1/2} \backslash M_2 + D_{3,2/3} ) = \phi(M_6)   \]
(see the proof of Proposition \ref{uniquenessJK} for the last equality). Finally Proposition \ref{dcs} implies
\[  (\phi \circ \psi)(R_6) \cap M_6 \subset \phi(M_6) \cap M_6 = \{ \}  \, . \]
In the same way we can show that $ (\phi \circ \psi)(R_{10}) \cap M_{10} = \{ \}$  for the lattice $ L =  I\!I_{6,2}(2_{I\!I}^{-2} 4_{I\!I}^{-2} 5^{+4}) $. \eop

\medskip

Now we return to the general case and replace $U$ by $(\phi \circ \psi)(U)$ if $N/m = 2$. Then $R_U = \{ \}$. Next we show that $\psi_F$ indeed has order $0$ at the cusp ${\cal C}(U)$. Recall that the class corresponding to $\psi_F$ is of cycle shape $\prod_{d|m} d^{\,b_d}$.

\begin{prp}\label{coeffsH}
For $ \alpha \in H $ we have $ c_{\alpha}(0) = \sum_{d|m, \, \alpha \in D_d} b_d $.
\end{prp}
{\em Proof:} The statement follows from a case-by-case analysis of the $ 11 $ reflective modular forms constructed in Section \ref{DM}. 
\eop 

\medskip
Now this observation and $ R_U = \{\} $ immediately imply that $ \psi_F $ does not vanish on $ \mathcal{C}(U) $. 

\begin{prp}
The expansion of $ \psi_F $ at $ \mathcal{C}(U) $ has order $ I_0 = 0 $.  
\end{prp}
{\em Proof:}
We evaluate the formula for the order $ I_0 $ given in Theorem \ref{kudlaprod}. We have seen in the proof of Proposition \ref{intersectionWithLDual} that $ L' \cap U^\perp = {\Lambda^g}' + (L' \cap U)$. Hence a vector in $ L' \cap U^\perp $ has non-negative norm and is isotropic if and only if it lies in $ L' \cap U $. We deduce
\[  I_0 =  - \sum_{r \in \Q} \, \sum_{\al \in H} c_\alpha(-r) \sigma_1(r)
  - \sum_{\al \in R_U} c_{\alpha}(-\alpha^2/2)\sigma_1(0) \, . \]
The second sum vanishes because $ R_U = \{\}$ and the first sum can be evaluated using Proposition \ref{coeffsH}. We obtain
\[  I_0 =  \frac{1}{24} \Bigg( \, \sum_{\alpha \in H} c_\alpha(0) \Bigg) - 1
        =  \frac{1}{24} \Bigg( \, \sum_{d|m} d b_d \Bigg) - 1  =  0  \, . 
      \]
This finishes the proof.  \eop 

\medskip 

Since $ H \simeq \Z/m\Z $ is generated by a single element, we can choose a basis $ (e_1, e_2) $ of $ L \cap U $ such that $ L' \cap U = \langle e_1/m, e_2 \rangle $. Let $ U' $ be an isotropic subspace dual to $ U $ and $ (e_1', e_2') $ a basis of $ U' $ such that $ (e_i, e_j') = \delta_{ij} $. Note that the basis $ (e_1, e_2) $ of $ L \cap U $ need not be equal to the basis $ (2f_1, f_2) $ chosen in the construction of $ U $ in Proposition \ref{egilstadir}.
We now calculate the constant term of the expansion of $ \psi_F $ at $ \mathcal{C}(U) $ with respect to this basis. 

\begin{prp}
The constant term $ \psi_0 $ in the expansion of $ \psi_F $ at $ \mathcal{C}(U) $ relative to the basis $ (e_1,e_2, e_1', e_2') $ of $ U + U' $ is given by
\[  \psi_0(w, \tau_2) = \kappa N^{(n-2)/4} \eta_g(\tau_2)   \]
where $ \kappa $ is a constant of absolute value $ 1 $. 
\end{prp}
{\em Proof:}
The first coefficient $ \psi_0 $ of the expansion of $ \psi_F $ at $ \mathcal{C}(U) $ is described in Theorem \ref{kudlaprod}. Since $ R_U = \{\} $, the second product is $ 1 $ so that 
\[ \psi_0(w, \tau_2) = \, \kappa \eta(\tau_2)^{c_0(0)} \prod_{\substack{\al \in H \\ \al\ne 0}}
  \left(\frac{\vartheta(-(\alpha,z_L),\tau_2)}{\eta(\tau_2)}e((\alpha,z_L)/2)^{(\alpha,e'_2)}\right)^{c_\alpha(0)/2}. \]
By our choice of basis we can represent $ \al \in H $ by an element $ je_1/m $ with $ 1 \leq j \leq m-1 $ so that $ (\alpha, z_L) = (j e_1/m, z_L ) = j/m $ and $ (\alpha, e_2') = (j e_1/m, e_2') = 0 $. Hence
\[  \psi_0(w,\tau_2) = \kappa \eta(\tau_2)^{c_0(0)} \prod_{j=1}^{m-1}
  \left(\frac{\vartheta(-j/m, \tau_2)}{\eta(\tau_2)}\right)^{c_{j e_1 / m}(0)/2}  \, . \]
We have $ c_0(0) = n-2 $. Using Proposition \ref{coeffsH} we find $ c_{j e_1 / m}(0) = \sum_{d \in I_j}b_d $ where $ d \in I_j $ if and only if $ j/m = k/d $ for some integer $ k $. It follows
\[  \psi_0(w,\tau_2) = \kappa \eta(\tau_2)^{n-2} \,\prod_{d|m} \Bigg( \prod_{k = 1}^{d-1}\frac{\vartheta(-k/d,\tau_2)}{\eta(\tau_2)}\Bigg)^{b_d/2} \, . \]
By the continuity of $ \vartheta/\eta $ we find for each $d|m$
\begin{align*}
  \prod_{k = 1}^{d-1}\frac{\vartheta(-k/d,\tau_2)}{\eta(\tau_2)} 
&= \, \lim_{z \rightarrow 0} \Bigg( \left( \frac{\vartheta(-z,\tau_2)}{\eta(\tau_2)} \right)^{-1} \; \prod_{k = 0}^{d-1}\frac{\vartheta(-z-k/d,\tau_2)}{\eta(\tau_2)} \Bigg)  \\
  &= \, (-1)^{d-1} \, \lim_{z \rightarrow 0} \Bigg( \left( \frac{\vartheta(-z,\tau_2)}{\eta(\tau_2)} \right)^{-1} \, \frac{\vartheta(-dz,d\tau_2)}{\eta(d\tau_2)} \Bigg)  \\
  &= \, (-1)^{d-1} d \, \eta(\tau_2)^{-2} \eta(d\tau_2)^2 \, .
\end{align*}
Here we used the product expansion of $ \vartheta/\eta $ (cf.\ Section \ref{brel}). Replacing $ \kappa $ by another constant of absolute value $ 1 $ we obtain
\[   \psi_0(w, \tau_2) \; = \; \kappa \Bigg( \, \prod_{d|m} d^{\, b_d/2} \Bigg) \; \eta(\tau_2)^{n-2 - \sum_{d|m} b_d} \; \eta_g(\tau_2) \, . \]
Now $ b_d = b_{N/d} $ implies
\[  \prod_{d|m} d^{\, b_d/2} = \prod_{d|N} d^{\, b_d/2} = N^{\sum_{d|N} b_d /4}   \, . \]
Inserting $ \sum_{d|N} b_d = n-2 $ finally yields the claim. \eop 

\medskip 

We can easily determine the constant term for any basis $ (f_1, f_2) $ of $ L \cap U $. After possibly replacing some $ f_i $ with $ -f_i $, we have $ f_1 = ae_1 - ce_2 $ and $ f_2 = -be_1 + de_2 $ with $ M = \begin{psmallmatrix} a & b \\ c & d \end{psmallmatrix} \in \text{SL}_2(\Z) $. We compute 
\[\psi_0^{(f_1,f_2)}(\tau) = \psi_0^{(e_1,e_2)}|_{M}(\tau)\] 
where $ \psi_0^{(e_1,e_2)} $ and $ \psi_0^{(f_1,f_2)} $ denote the constant terms of the expansions with respect to $ (e_1, e_2) $ and $ (f_1, f_2) $.

We summarise the properties of the cusp $ \mathcal{C}(U) $ in the following theorem. 

\begin{thm}\label{falco}
	There exists a $1$-dimensional cusp $ \mathcal{C} $ of $ O(L,F)^+ \backslash \mathcal{H} $ of type $ H \subset D^{N/m} $ of order $ m $ with associated lattice $ \Lambda^g $ such that $ \psi_F $ does not vanish on $ \mathcal{C} $. The constant term $ \psi_0 $ of the expansion at $ \mathcal{C} $ is given by
	\[\psi_0(w,\tau_2) = \kappa N^{(n-2)/4} \, \eta_g|_M(\tau_2)\]
	where $ \kappa $ is a constant of absolute value $ 1 $ and $ M \in \text{SL}_2(\Z) $ depends on the parametrisation of the neighbourhood of $ \mathcal{C} $. For $ N/m = 1 $ the cusp $ \mathcal{C} $ is the unique split cusp of type $ H $ with associated lattice $ \Lambda^g $.
\end{thm}

The expansions of $\psi_F$ at the $1$-dimensional cusps of type $ 0 $ and the special cusp constructed above can be visualised as follows. We choose $ L = I\!I_{12,2}(2_2^{+2} 4_{I\!I}^{+6}) $ as an example. 

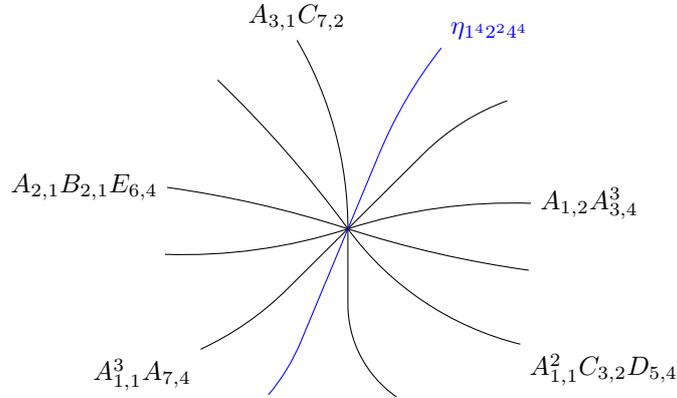
\begin{center}
	\begin{tikzpicture}
		\draw (0,0) -- (0,-1);
		\draw (0,-1) arc (180:235:1.5cm); 
		\draw (0,0) arc (0:30:5cm) node[anchor = south] {$A_{3,1}C_{7,2}$};
		
		\draw (0,0) arc (72:82:14cm) node[anchor = east] {$A_{2,1}B_{2,1}E_{6,4}$};
		\draw (0,0) arc (252:262:14cm);
		
		\draw (-0.8,-0.8) -- (1,1); 
		\draw (1,1) arc (135:110:3cm);
		\draw (-0.8,-0.8) arc (-45:-65:4cm) node[anchor = north east] {$A_{1,1}^3A_{7,4}$};
		
		\draw (0,0) arc (216:256:4cm) node[anchor = north west] {$A_{1,1}^2C_{3,2}D_{5,4}$};
		\draw (0,0) arc (36:46:15cm);
		
		\draw (0,0) arc (108:88:7cm) node[anchor = west] {$A_{1,2}A_{3,4}^3$};
		\draw (0,0) arc (288:268:7cm);
		
		\draw[color = blue] (-0.6213,-1.5) -- (0,0);
		\draw[color = blue] (-0.6213,-1.5) arc (-22.5:-40:2.7cm);
		\draw[color = blue] (0,0) -- (0.4142,1);
		\draw[color = blue] (0.4142,1) arc (157.5:142:6cm) node[anchor = south west] {$\eta_{1^4 2^2 4^4}$};
	\end{tikzpicture}
\end{center}

\medskip

The black lines illustrate the $1$-dimensional cusps of type $ 0 $ with their respective root systems indicated next to them. The blue line illustrates the special cusp with associated lattice $ \Lambda^g $. The restriction of $ \psi_F $ to this cusp is up to a constant equal to the eta product $ \eta_{1^4 2^2 4^4} $. The $1$-dimensional cusps of type $ 0 $ and the special cusp all share the unique $ 0 $-dimensional cusp of type $ 0 $ as a common boundary point depicted as the intersection of the black and blue lines.

\section{Holomorphic vertex operator algebras of central charge $24$} \label{AMK}

We show that the character of a holomorphic vertex operator algebra of central charge $24$ with non-trivial, semisimple weight-$1$ space defines a reflective modular form which lifts to a reflective automorphic product of singular weight. The corresponding modular variety has a canonical $1$-di\-men\-sio\-nal cusp of type $0$ whose root system determines the affine structure of $V$. It follows that under certain regularity assumptions the holomorphic vertex operator algebras of central charge $24$ with non-trivial, semisimple weight-$1$ space fall into at most $11$ classes with at most $69$ affine structures.

\subsection*{Affine Kac-Moody algebras and vertex operator algebras}

We recall some results on affine Kac-Moody algebras and the corresponding vertex operator algebras from \cite{CKS,DS,FMS,FZ,JF,Kac,KW,KP}.

\medskip

Let $\g$ be a finite-dimensional simple complex Lie algebra of rank $l$, $\h$ a Cartan subalgebra of $\g$ and $\Phi \subset \h'$ the corresponding set of roots. We normalize the non-degenerate, invariant, symmetric bilinear form $( \, , \,)$ on $\g$ such that the long roots have norm $2$. The bilinear form induces an isomorphism $\nu : \h \to \h'$. For $\al \in \Phi$ we denote by $\al^{\vee}$ the inverse image of $2\al/\al^2$. The root lattice $Q$ is the $\Z$-module in $\h'$ generated by $\Phi$. The coroot lattice $Q^{\vee}$ is the $\Z$-module in $\h$ generated by $\Phi^{\vee} = \{ \al^{\vee} \, | \, \al \in \Phi \}$. It is a positive-definite even lattice of rank $l$.
The weight lattice $P = \{ \lambda \in \h' \, | \, \lambda(Q^{\vee}) \subset \Z \} \subset \h'$ is the dual of $Q^{\vee}$ and analogously $P^{\vee} \subset \h$ the dual of $Q$. Let $ \Delta = \{ \al_1, \ldots, \al_l\} \subset \Phi$ be a set of simple roots. Then \[   Q = \sum_{i=1}^l \Z \al_i \, . \]

\medskip

The {\em untwisted affine Kac-Moody algebra} corresponding to $\g$ is the Lie algebra 
\[  \hat{\g} = \g \otimes \C[t,t^{-1}] \oplus \C K \oplus \C d \]
where $K$ is central and
\begin{align*}
	[ \, a \otimes t^m , b \otimes t^n \, ] &= 
	[a,b] \otimes t^{m+n} + m \delta_{m+n} (a,b) K \, ,    \\
	[ \, d , a \otimes t^n \, ] & = n a \otimes t^n   \, .
\end{align*}
The vector space $\hat{\h} = \h + \C K + \C d$ is a commutative subalgebra of $\hat{\g}$.

We extend a linear function $\lambda$ on $\h$ to $\hat{\h}$ by setting $\lambda(K) = \lambda(d) = 0$. Furthermore we define linear functions $\Lambda_0$ and $\delta$ on $\hat{\h}$ by $\Lambda_0|_{\h} = \delta|_{\h} = 0$ and 
\begin{align*}
	\qquad & & \Lambda_0(K) &= 1 \, , & \Lambda_0 (d) &= 0 \, , & & \qquad\\
	\qquad & & \delta(K)    &= 0 \, , & \delta(d)     &= 1 \, . & & \qquad
\end{align*}
Then 
\[    \hat{\h}' = \h' + \C \Lambda_0 + \C \delta     \]
and we have a natural projection 
$\hat{\h}' \to \h'$, $\lambda \mapsto \ov{\lambda}$ 
with $\ov{\Lambda_0} = \ov{\delta} = 0$. A linear function $\lambda$ in $\hat{\h}'$ can be written as $\lambda = \ov{\lambda} + \lambda(K) \Lambda_0 + \lambda(d) \delta$ and $\lambda(K)$ is called the {\em level} of $\lambda$. 

We also extend the bilinear form from $\g$ to $\hat{\g}$ by setting
\begin{gather*}
	( a \otimes t^m , b \otimes t^n ) = \delta_{m+n} (a,b) \, , \quad 
	( a \otimes t^m , K ) = ( a \otimes t^m , d ) = 0 \, , \\
	(K,d) = 1 \, , \quad  ( K, K) = ( d , d ) = 0 \, . 
\end{gather*}
The isomorphism $\hat{\h} \to \hat{\h}'$ induced by $( \, , \,)$ extends the map $\nu : \h \to \h '$.

Define $\al_0 = \delta - \theta$ where $\theta = \sum_{i=1}^{l} a_i \al_i$ with Coxeter labels $a_i$ (see \cite{Kac}, p.\ 54, Table Aff 1) is the highest root of $\g$ and $\al_0^{\vee} = K - \theta^{\vee}$. Then $\hat{\Delta} = \{ \al_0, \al_1,\ldots,\al_l \}$ is a set of simple roots of $\hat{\g}$ and $a_i^{\vee} \nu(\al^{\vee}_i) = a_i \al_i$. The fundamental weights $\Lambda_i \in \hat{\h}'$, $i = 0, \ldots, l$ are defined by
\[  \Lambda_i ( \al^{\vee}_j ) = \delta_{ij} \, ,   
\quad  \Lambda_i ( d ) = 0  \, . 
\]
Then $\Lambda_i = \ov{\Lambda_i} + a^{\vee}_i \Lambda_0$ and the projections $\ov{\Lambda_1}, \ldots, \ov{\Lambda_l}$ are the fundamental weights of $\g$.
We can write the weight lattice
\[  \hat{P} = \{ \lambda \in \hat{\h}' \, | \,  \lambda(\al^{\vee}_i) \in \Z \, \text{ for } i = 0,\ldots, l  \}  \]
as $\hat{P} = \sum_{i=0}^{l} \Z \Lambda_i + \C \delta$.
For $k \in \Z_{>0}$ let
\begin{align*}
	\hat{P}^k &= \{ \lambda \in \hat{P} \, | \, \lambda(K) = k   \},  \\
	\hat{P}_+ &= \{ \lambda \in \hat{P} \, | \, \lambda(\al^{\vee}_i) \geq 0 \, \text{ for } i = 0, \ldots, l  \}
\end{align*}
and $\hat{P}^k_+ = \hat{P}^k \cap \hat{P}_+$. Note that
$\hat{P}^k_+ \! \mod \C \delta = \{ \lambda \in \sum_{i=0}^l \Z_{\geq 0} \Lambda_i \, | \, \lambda(K) = k \}$ is finite.

The Weyl group $\hat{W}$ of $\hat{\g}$ is the subgroup of $\text{GL}(\hat{\h}')$ generated by the reflections $\sigma_{\al_i}$, $i = 0,\ldots,l$. Let $M = \nu(Q^{\vee})$. Then the translations $t_{\al}$, $\al \in M$ defined by
\[ t_{\al}(\lambda) = \lambda + \lambda(K)\al- ( (\lambda,\al) + \lambda(K) \al^2/2 )\delta   \]
form a normal subgroup of $\hat{W}$ isomorphic to $M$ and $\hat{W}$ is the semidirect product of $M$ and the subgroup $W$ generated by the $\sigma_{\al_i}$, $i = 1,\ldots,l$, i.e.
\[  \hat{W} \simeq M \rtimes W \, . \]
The set $\hat{\Phi}$ of roots of $\hat{\g}$ is invariant under $\hat{W}$. It is invariant even under the larger group $\hat{W}_0$ generated by the translations in $\nu(P^{\vee})$ (recall that $Q^{\vee} \subset P^{\vee}$) and the reflections $\sigma_{\al_i}$, $i = 1,\ldots,l$ (see Section 1.3 in \cite{KW}). We denote by $\af$ the affine action of $\hat{W}_0$ on $\h'_{\R} = Q \otimes_{\Z} \R$ (see \cite{Kac}, \S$\,$6.6). Let $\hat{W}_0^+$ be the subgroup of $\hat{W}_0$ preserving the positive roots $\hat{\Phi}^+$ of $\hat{\g}$ and
\[  J = \{ j \, | \, a_j = 1 \} \subset \{0, \ldots,l \} \, .   \]
For each $\sigma \in \hat{W}_0^+$ there is a unique $j \in J$ such that $\ov{\sigma(\Lambda_0)} = \af(\sigma)(0) = \ov{\Lambda_j}$. It follows that $\sigma = t_{\ov \Lambda_j} w$ for some $w \in W$. The induced map $\hat{W}_0^+ \to J$ is a bijection and endows $J$ with a group structure. Then the map $J \to P/Q$, $j \mapsto \ov{\Lambda_j} + Q$ is an isomorphism and we have the following sequence of group isomorphisms
\[  \hat{W}_0^+ \longrightarrow J \longrightarrow P/Q  \, . \]
We denote the element in $\hat{W}_0^+$ corresponding to $j \in J$ under the above map by $\sigma_j$. The simple roots $\hat{\Delta}$ are invariant under $\hat{W}_0^+$ so that $\hat{W}_0^+$ acts on the Dynkin diagram $\hat{\Gamma}$ of $\hat{\g}$. We have
\[  \Aut(\hat{\Gamma}) = \hat{W}_0^+ \ltimes \Aut(\Gamma) \]
(\cite{KW}, Proposition 1.3). We list the groups $P/Q$ in the following tables:
\[
\renewcommand{\arraystretch}{1.2}
\begin{array}{c|c|c|c|c}
	A_n & B_n & C_n & D_{2n} & D_{2n+1}  \\  \hline
	\Z/(n+1)\Z & \Z/2\Z & \Z/2\Z & \Z/2\Z \times \Z/2\Z & \Z/4\Z \\
	\multicolumn{5}{c}{} \\
	E_6 & E_7 & E_8 & F_4 & G_2 \\  \hline
	\Z/3\Z & \Z/2\Z & 1 & 1 & 1 
\end{array}            
\]

The Weyl vector $\rho = \sum_{i=0}^{l} \Lambda_i \in \hat{\h}'$ satisfies the strange formula of Freudenthal-de Vries (see \cite{Kac}, (12.1.8)) 
\[  \frac{\rho^2}{2h^{\vee}} = \frac{\dim(\g)}{24}   \, . \]
We remark that $\rho^2 = \ov{\rho}^2$.

Let $k \in \Z_{>0}$. Then $k+h^{\vee} \neq 0$. We consider the irreducible highest-weight module $L(\Lambda)$ associated with $\Lambda \in \hat{P}^k_+$.

A weight $\lambda$ of $L(\Lambda)$ is called maximal if $\lambda + \delta$ is not a weight. Then the set $\hat{P}(\Lambda)$ of weights can be decomposed into the disjoint union
\[  \hat{P}(\Lambda) = \bigcup_{\lambda \in \text{max}(\Lambda)} \{ \lambda - n \delta \, | \, n \in \Z_{\geq 0}  \}  \]
where $\text{max}(\Lambda)$ is the set of maximal weights.

We denote by 
\[ 
  m_{\Lambda}  = \frac{(\Lambda+\rho)^2}{2(k+h^{\vee})} - \frac{\rho^2}{2h^{\vee}} \, , \qquad
h_{\Lambda}  = \frac{(\Lambda + 2 \rho, \Lambda)}{2(k+h^{\vee})} \, , \qquad 
c_k = \frac{k \dim(\g)}{k+h^{\vee}}
\]
the {\em modular}, {\em vacuum} and {\em conformal anomaly} of $L(\Lambda)$ (see \cite{Kac}, \S$\,$12.7 and \S$\,$12.8). They are related by
\[  m_{\Lambda} = h_{\Lambda} - \frac{1}{24} c_k \, .  \]

The string functions
\[    c_{\lambda}^{\Lambda} = e^{-m_{\Lambda, \lambda} \delta}  \sum_{n \in \C} \mult_{L(\Lambda)}(\lambda-n\delta) e^{-n\delta}  \]
where $m_{\Lambda, \lambda} = m_{\Lambda} - \lambda^2/2k$ satisfy
\[   c_{\lambda}^{\Lambda} = c_{w(\lambda) + k\bt + b\delta}^{\Lambda + a \delta}  \]
for all $a,b \in \C$, $w \in W$ and $\bt \in M$. They are also invariant under $\hat{W}_0^+$, i.e.\
\[   c_{\sigma(\lambda)}^{\sigma(\Lambda)} = c_{\lambda}^{\Lambda}  \]
for all $\sigma \in \hat{W}_0^+$ (cf.\ Proposition 5.1 in \cite{CKS}).
The normalised character (\cite{Kac}, (12.7.12)) 
\[    \chi_{\Lambda} = e^{-m_{\Lambda} \delta} \ch_{L(\Lambda)}   \]
of $L(\Lambda)$ can be written as 
\[    \chi_{\Lambda} = \sum_{\lambda \in \hat{P}^k \bmod (kM + \C \delta)} c_{\lambda}^{\Lambda} \theta_{\lambda}  \]
with theta functions
\[    \theta_{\lambda} = e^{k\Lambda_0} \sum_{\bt \in M + \ov{\lambda}/k } e^{ - \delta k (\bt,\bt)/2 + k \bt }  \, . \] 

Recall that $P \subset \h'$ is the dual of $M$ with respect to $( \, , \, )$. We introduce a new bilinear form $( \, , \, )_k = k ( \, , \, )$ on $\h$ and on $\h'$ via $\nu$. We define the lattice $M_k \subset \h'$ which as a set is equal to $M$ but has $( \, , \, )_k$ as its bilinear form. Then $M_k$ is a positive-definite even lattice and the dual of $M_k$ in $\h'$ is given by $M_k' = (1/k)M'$ as a set. We obtain bijections
\[
\renewcommand{\arraystretch}{1.2}
\begin{array}{clclc}
	\hat{P}^k \, \bmod (kM + \C \delta)     & \longrightarrow & M'/kM             & \longrightarrow  & M_k'/M_k  \\
	\lambda \; \bmod (kM + \C \delta) & \longmapsto     & \ov{\lambda} + kM & \longmapsto      & \ov{\lambda}/k + M_k \, . 
\end{array}            
\]
We denote the composition of these maps by $\pi^k$.
The group $\hat{W}_0^+$ preserves the set $\hat{P}^k$ as well as $kM+\C\delta$. Hence we can define an action of $\hat{W}_0^+$ on $\hat{P}^k \, \bmod (kM + \C \delta)$ and on $M_k'/M_k$. It satisfies
\[ \pi^k \big( \sigma(\lambda) \, \bmod (kM + \C \delta) \big) =   \af(\sigma)(\overline{\lambda}/k) + M_k \, . \]

We can now rewrite the character as 
\[    \chi_{\Lambda} = \sum_{\lambda \in M_k'/M_k} c_{\lambda}^{\Lambda} \theta_{\lambda}  \]
with
\[ c_{\lambda}^{\Lambda}
= q^{m_{\Lambda, \mu}} \sum_{n \in \C} \mult_{L(\Lambda)}(\mu-n\delta) q^n \]
and
\[    \theta_{\lambda} = e^{k\Lambda_0} \sum_{\bt \in \lambda + M_k} q^{ (\bt, \bt)_k/2} e^{k \bt}   \]
where $\mu \in \hat{P}^k$ is such that $\pi^k \big( \mu  \, \bmod (kM + \C \delta) \big) = \lambda$ and $q = e^{-\delta}$.

Under the specialisation $ e^\al \mapsto e(\al(v)) $ for $ \alpha \in \hat{\h}' $, $ v \in \hat{\h} $, the character $\chi_{\Lambda}$ defines a holomorphic function on $Y = \h \oplus \C K \oplus H(-d)$ with the complex upper halfplane $H$. It transforms as a Jacobi form under a suitable Jacobi group (see Theorem 13.8 in \cite{Kac}).

The modular group $\text{SL}_2(\Z)$ acts on the linear span of the characters $\chi_{\Lambda}$, $\Lambda \in \hat{P}^k_+ \! \mod \C \delta$. The action of $S = \left( \begin{smallmatrix} 0 & -1 \\ 1 & 0 \end{smallmatrix} \right)$ is described by the ${\cal S}$-matrix (see Theorem 13.8 in \cite{Kac}).

\medskip

The irreducible $\hat{\g}$-module $L(k\Lambda_0)$ carries the structure of a vertex operator algebra of central charge $c_k$ \cite{FZ}. It is simple and strongly rational and called the {\em simple affine vertex operator algebra of level $k$}. The irreducible modules of $L(k\Lambda_0)$ are the $\hat {\g}$-modules $L(\Lambda)$, $\Lambda \in \hat{P}^k_+ \! \mod \C \delta$ of conformal weight $h_{\Lambda}$. The ${\cal S}$-matrix of $L(k\Lambda_0)$ as a vertex operator algebra \cite{Z} is exactly the ${\cal S}$-matrix described above. The map
$\hat{P}^k_+  \, \bmod \C \delta  \, \to \,  \hat{P}^k \, \bmod (kM + \C \delta)$ is injective (see the proof of Proposition 4.1.2 in \cite{DS}) so that the composition with $\pi^k$ defines an injective map
\[  \pi_+^k : \hat{P}^k_+  \, \bmod \C \delta  \; \longrightarrow  \; M_k'/M_k \, .   \]
In particular we can identify the irreducible modules of $L(k\Lambda_0)$ with a subset of $M_k'/M_k$. 

The {\em cominimal simple currents} of $L(k\Lambda_0)$ are the irreducible modules $L(k\Lambda_j)$, $j\in J$. The fusion product $\boxtimes$ of modules is closed on the set $S_J$ of cominimal simple currents and turns it into an abelian group (cf.\ \cite{DLM}, \cite{JF}). Note that the Weyl group $\hat{W}$ stabilises pointwise the image of $S_J$ under $\pi^k_+$. (This can be verified easily for the reflections in simple roots.) 

The group $\hat{W}_0^+$ acts on $\hat{P}^k_+ \! \mod \C \delta$ and preserves $S_J$. Furthermore the $\cal S$-matrix has the following symmetry
\[  {\cal S}_{\sigma_j(\lambda),\mu} = e^{-2 \pi i(\overline{\Lambda_j},\overline{\mu})} {\cal S}_{\lambda,\mu}  \]
for all $\lambda, \mu \in \hat{P}^k_+ \! \mod \C\delta$, $j\in J$ (\cite{FMS}, (14.255)). 
This implies
\[  L(k\Lambda_j) \boxtimes_{L(k\Lambda_0)} L(\lambda) = L(\sigma_j(\lambda)) \]
for all $\lambda \in \hat{P}^k_+ \! \mod \C\delta$, $j\in J$ (\cite{DS}, Proposition 4.1.4). It follows
\begin{align*}
	L(k\Lambda_i) \boxtimes_{L(k\Lambda_0)} L(k\Lambda_j) &= L(\sigma_i(k\Lambda_j)) = L((\sigma_i\sigma_j)(k\Lambda_0))
	= L(\sigma_{i+j}(k\Lambda_0))  \\
	&= L(k\Lambda_{i+j})  
\end{align*}
i.e.\ the natural map $\hat{W}_0^+ \to S_J, \, \sigma_j \mapsto L(k\Lambda_j)$ is a group isomorphism. Since
\[  \pi_+^k(k\Lambda_{i+j})
= \ov{\sigma_i(\Lambda_j)} + M_k = \af(\sigma_i)(\ov{\Lambda_j}) + M_k
= (\ov{\Lambda_i} + \ov{\Lambda_j}) + M_k  \, ,  \]
the restriction $\pi_+^k|_{S_J}$ is a group homomorphism. We denote $\pi_+^k(S_J) = H_J$. Then we have the following group isomorphisms
\[  \hat{W}_0^+ \longrightarrow S_J \longrightarrow H_J  \, . \]

Finally the conformal weights of the cominimal simple currents $L(k\Lambda_j)$ are given by
\[   h_{k\Lambda_j} = \Lambda_j^2/2k   \, . \]
They satisfy
\[  \q( \pi_+^k(k\Lambda_j)) = h_{k\Lambda_j}  \! \mod 1    \]
where $\q : M_k'/M_k \to \Q/\Z$ denotes the quadratic form on $M_k'/M_k$ and
\[  (\pi_+^k(k\Lambda_j),\pi_+^k(\lambda)) = h_{\sigma_j(\lambda)} - h_{k\Lambda_j} - h_{\lambda}  \! \mod 1  \]
for all $\lambda \in \hat{P}^k_+ \! \mod \C\delta$ (\cite{DS}, (615)).

\subsection*{Holomorphic vertex operator algebras of central charge $24$}

We associate to a holomorphic vertex operator algebra $V$ of central charge $24$ with non-trivial, semisimple weight-$1$ space a reflective automorphic product of singular weight. The corresponding modular variety has a canonical $1$-di\-men\-sio\-nal cusp of type $0$ whose root system determines the affine structure of $V$. 
By our classification results (Theorems \ref{blazefoley} and \ref{heaven17}) this implies that $V$ falls into one of H\"ohn's $11$ classes and the affine structure of $V$ is one of the $69$ Lie algebras described by Schellekens.

\medskip

Let $V$ be a strongly rational, holomorphic vertex operator algebra of central charge $24$. We normalize the non-degenerate, invariant, symmetric bilinear form $\la \, , \, \ra$ on $V$ such that $\la \vac, \vac \ra = -1$. The weight-$1$ subspace $V_1$ has the structure of a reductive Lie algebra. More precisely $V_1$ is either $0$ or abelian and of rank $24$ or non-zero and semisimple \cite{DM1}. In the second case $V$ is isomorphic to the vertex operator algebra associated with the Leech lattice.
Suppose $V_1 = \g \neq 0$ is semisimple with simple components $\mathfrak{g}_i$. Let $\h \subset \g$ be a Cartan subalgebra of $\g$. Then $\h_i = \g_i \cap \h$ is a Cartan subalgebra of $\g_i$ and the restriction of $\la \, , \, \ra$ to $\g_i$ satisfies $\la \, , \, \ra|_{\g_i} = k_i( \, , \, )$ where $k_i$ is a positive integer and $( \, , \, )$ the non-degenerate, invariant, symmetric bilinear form on $\g_i$ normalised such that the long roots have norm $2$ \cite{DM2}. To simplify the notation we identify $ \h_i $ with its dual $ \h_i' $ through $ (\,,\,) $. The decomposition
\[   \g = \g_{1,k_1} \oplus \ldots \oplus \g_{r,k_r}   \]
is called the {\em affine structure} of $V$. Schellekens \cite{ANS} showed that 
\[  \frac{h_i^{\vee}}{k_i} = \frac{\dim(\g)-24}{24}  \, .  \]
The subalgebra $V_{\g} = \la V_1 \ra$ generated by $V_1$ in $V$ is isomorphic to 
\[  V_{\g} \simeq L_{\g_1,k_1} \otimes \ldots \otimes L_{\g_r,k_r} \]
where $L_{\g_i,k_i}$ is the simple affine vertex operator algebra of level $k_i$ associated with $\g_i$. The irreducible modules of $V_{\g}$ are of the form
\[  V_{\g}(\lambda) \simeq L_{\g_1,k_1}(\lambda_1) \otimes \ldots \otimes L_{\g_r,k_r}(\lambda_r)   \]
where $L_{\g_i,k_i}(\lambda_i)$ is an irreducible module of $L_{\g_i,k_i}$ and $\lambda = (\lambda_1, \ldots, \lambda_r)$. Using the projection $\pi = (\pi_+^{k_1},\ldots,\pi_+^{k_r})$ we can embed the set $P_{\g}$ of weights $\lambda$ into the discriminant form of the lattice
\[  M = M_{k_1} \oplus \ldots \oplus M_{k_r}  \]
where
\[   M_{k_i} = Q^{\vee}_i(k_i)   \]
is the coroot lattice of $\g_i$ rescaled by $k_i$, i.e.\ the bilinear form of $ M_{k_i} $ is given by the restriction of $ \la \, , \, \ra $. A module $V_{\g}(\lambda)$ is called a {\em cominimal simple current} if it is the tensor product of cominimal simple currents. We denote the corresponding set of weights by $S_{\g}$. (We will sometimes identify a module with its highest weight.) The cominimal simple currents form an abelian group under fusion and act on the irreducible modules of $V_{\g}$ \cite{JF, DLM}. For $\sigma \in S_{\g}$ we denote the action of $\sigma$ on $P_{\g}$ again by $\sigma$, i.e.\
\[  V_{\g}(\sigma) \boxtimes V_{\g}(\lambda) \simeq V_{\g}(\sigma(\lambda))   \, . \]
Since $V_{\g}$ is rational, $V$ decomposes into finitely many irreducible modules under the action of $V_{\g}$. We write
\[  V \simeq \bigoplus_{\lambda \in P_{\g}} m_{\lambda} V_{\g}(\lambda)   \]
with multiplicities $m_{\lambda}\in \Z_{\geq 0}$ and set $P_V = \{  \lambda \in P_{\g} \, | \, m_{\lambda} \neq 0 \}$ and $S_V = S_{\g} \cap P_V$. Clearly the modules $V_{\g}(\lambda)$, $\lambda \in P_V$ have conformal weight $h_{\lambda} = 0 \! \mod 1$. Schellekens observed that the multiplicities are invariant under $S_V$ (\cite{ANS}, Section 3).

\begin{prp} \label{schellekens}
	Let $\sigma \in S_V$ and $\lambda \in P_V$. Then 
	\[  m_{\sigma(\lambda)} = m_{\lambda}   \, .   \]
	In particular $S_V$ is a group and $m_\sigma=1$ for all $\sigma\in S_V$.
\end{prp}

We denote by $ H_V $ the image of $ S_V $ under $ \pi $.

\begin{prp} \label{sotif}
The group $H_V$ is an isotropic subgroup of $M'/M$ and $\pi(P_V) \subset H_V^{\perp}$.
\end{prp}
{\em Proof:}
The first statement is clear. The second follows from the last formula in the previous section and the fact that $P_V$ is stable under $S_V$. \eop

\medskip

The {\em character} $\chi_V : \h \times H \to \C$ of $V$ is
defined as the trace
\[  \chi_V(v,\tau) = \tr_V \, e^{2\pi iv_0}q^{L_0-1}   \, . \]

\begin{prp} \label{VesperLynd}
	The character of $V$ can be written as
	\[  \chi_V(v,\tau) = \sum_{\gamma \in M'/M} F^M_{\gamma}(\tau) \theta^M_{\gamma}(v,\tau)   \]
	where 
	\[  F^M_{\gamma}(\tau) =  \sum_{\lambda \in P_V} m_{\lambda}
	\prod_{i=1}^r c_{\gamma_i}^{\lambda_i}(\tau)   \]
	and
	\[  \theta^M_{\gamma}(v,\tau) = \prod_{i=1}^r \theta_{\gamma_i}(v_i,\tau)   \]
	is the classical Jacobi theta function of $\gamma + M$.
\end{prp}
{\em Proof:}
We have 
\begin{align*}
	\chi_V (v,\tau)
	& = \sum_{\lambda \in P_V} m_{\lambda} \prod_{i=1}^r \chi_{L_{\g_i,k_i}(\lambda_i)}(v_i,\tau)   \\
	& =  \sum_{\lambda \in P_V} m_{\lambda} \prod_{i=1}^r \sum_{\gamma_i \in M_{k_i}'/M_{k_i}}
	c_{\gamma_i}^{\lambda_i}(\tau) \theta_{\gamma_i}(v_i,\tau)  \\
	& =  \sum_{\gamma \in M'/M} \sum_{\lambda \in P_V} m_{\lambda}
	\prod_{i=1}^r c_{\gamma_i}^{\lambda_i}(\tau) \theta_{\gamma_i}(v_i,\tau) 
\end{align*}
with
\[  \theta_{\gamma_i}(v_i,\tau) = \sum_{\al \in \gamma_i + M_{k_i}} q^{\la \al,\al \ra/2} e(\la \al,v_i \ra)  \]
(cf.\ \cite{Kac}, (13.2.5)). \eop

\medskip
Although it is not necessary for many of the following results, we assume from now on that $\g$ has even rank.

\begin{thm}\label{radiohead}
	The function $F^M(\tau) = \sum_{\gamma \in M'/M} F^M_{\gamma}(\tau) e^{\gamma}$ has the following properties:
	\begin{enumerate}[i)]
		\item $F^M$ is a modular form of weight $-\rk(\g)/2$ for the Weil representation of $M'/M$. 
		\item $F^M_0(\tau) = q^{-1} + \rk(\g) + \ldots$
		\item $F^M_{\gamma}(\tau) = 0$ if $\gamma \notin H_V^{\perp}$.
		\item $F^M_{\gamma}(\tau) = F^M_{\gamma+s}(\tau)$ for all $\gamma \in M'/M$, $s \in H_V$.
		\item The pole order of $F^M_{\gamma}(\tau)$, $\gamma \in M'/M$ is bounded by $1$ with equality if and only if $\gamma \in H_V$.
	\end{enumerate}
\end{thm}
{\em Proof:}
The character $\chi_V$ is a weakly holomorphic Jacobi form of lattice index $M$ and weight $0$ (see Theorem 1.1 in \cite{KM}). This implies that the decomposition in Proposition \ref{VesperLynd} is the theta decomposition \cite{G91} of $\chi_V$. This proves the first claim. 

Let $\al \in \h'$. Decompose $\al = \al_1 + \ldots + \al_r$ with $\al_i \in \h_i'$. Then for $\lambda_i \in \h_i'$  the weight space of $L_{\mathfrak{g}_i,k_i}(\lambda_i)$ of degree $\al_i$ is given by
\[   L_{\g_i,k_i}(\lambda_i)_{\al_i}
= \bigoplus_{n \in \C} L_{\mathfrak{g}_i,k_i}(\lambda_i)_{\al_i + k_i \Lambda_0 - n\delta}  \]
(to keep the notation simple we suppress the index $i$ at $\Lambda_0$ and $\delta$). Since $L_0$ acts as $h_{\lambda_i}-d$ on this space (see \cite{Kac}, Corollary 12.8), we obtain
\begin{align*}
	\tr_{L_{\g_i,k_i}(\lambda_i)_{\al_i}} q^{L_0-c_{k_i}/24}
	&= q^{h_{\lambda_i}-c_{k_i}/24} \sum_{n \in \C}
	\mult_{L_{\g_i,k_i}(\lambda_i)}(\al_i + k_i \Lambda_0 - n \delta ) q^n  \\
	&= q^{(\al_i,\al_i)/2k_i} c_{\al_i + k_i \Lambda_0}^{\lambda_i}(\tau)  \, .
\end{align*}
Hence for a weight $\lambda = \lambda_1 + \ldots + \lambda_r$ and $\al = \al_1 + \ldots + \al_r$ such that
$\al_i + k_i \Lambda_0 \in \hat{P}^{k_i} \bmod \C \delta$
(this means that $\al_i$ is in the weight lattice of $\g_{i}$)
we have
\[  \tr_{V_{\g}(\lambda)_{\al}} q^{L_0-1} = q^{\la \beta, \beta \ra / 2} \prod_{i=1}^r c^{\lambda_i}_{\beta_i}(\tau)\]
where $\beta_i = \al_i/k_i$ and
$\beta = \beta_1 + \ldots + \beta_r \in M'$. Note that $\beta_i + M_{k_i} = \pi^{k_i} (\al_i + k_i \Lambda_0 + k_i Q_i^{\vee} + \C \delta)$. We can compute 
\begin{align*}
	\tr_{V_{\al}} q^{L_0-1} 
	& = \sum_{\lambda \in P_V} m_{\lambda} \tr_{V_\mathfrak{g}(\lambda)_\al} q^{L_0-1} \\
	& = q^{\la \beta, \beta \ra / 2} \sum_{\lambda \in P_V} m_{\lambda} \prod_{i=1}^r c^{\lambda_i}_{\beta_i}(\tau)  \\
	& = q^{\la \beta, \beta \ra / 2} F^M_\beta(\tau). 
\end{align*}
The case $\beta = 0$ of this equation implies the second assertion.

Let $\gamma \in M'/M$ such that $F^M_{\gamma} \neq 0$. Then there is an element $\lambda \in P_V$ such that $c^{\lambda_i}_{\gamma_i}\neq 0$ for all $i= 0, \ldots,r$. For each $i$ we choose a weight $\mu_i$ of $L_{\mathfrak{g_i},k_i}(\lambda_i)$ in the inverse image of $\gamma_i$ under $\pi^{k_i}$. Then
\[   \lambda_i - \mu_i = \sum_l r_l \al_l  \]
where the $r_l$ are suitable integers and the $\al_l$ are simple roots of $\hat{\g}_i$. It follows
\begin{multline*}
  ( k_i \Lambda_{j}, \lambda_i - \mu_i )
	= k_i \sum_l r_l (\Lambda_{j},\al_l)
	= k_i \sum_l r_l \frac{a^{\vee}_l}{a_l} (\Lambda_{j},\al^{\vee}_l)   \\ 
	= k_i r_j \frac{a^{\vee}_j}{a_j} \
	= 0 \mod k_i
\end{multline*}
for all fundamental weights $\Lambda_j$ of $\hat{\g}_i$ with $a_j=1$. This implies
\[  (\pi(\sigma),\pi(\lambda)-\gamma) = 0 \mod 1  \]
for all $\sigma \in S_\g $. For a cominimal simple current $\sigma \in S_{\g}$ we have 
\[  (\pi(\sigma),\pi(\lambda)) = h_{\sigma(\lambda)} - h_{\sigma} - h_\lambda\mod 1  \]
(see the last formula in the previous section). In the special case $\sigma\in S_V$ all these conformal weights have to be integers because they correspond to irreducible modules contained in  the vertex operator algebra $V$. We obtain
\[  (\pi(\sigma),\pi(\lambda)) = 0 \mod 1  \, . \]
This implies $(s,\gamma) = 0 \mod 1$ for all $s \in H_V$ and the third item follows.

Let $\gamma \in M'/M$ and $s \in H_V$. Denote by $s_i$ the $i$-th component of $s$. In the following we again suppress some indices $i$ to simplify the notation. The element $\sigma_i$ corresponding to $s_i$ in the group $\hat{W}_0^+$ associated with $\hat{\g}_i$ can be written as $\sigma_i = t_{\ov{\Lambda_{j}}} w_i$ for some $j \in J$ and $w_i \in W_i$. Choose $\mu \in M'$ such that $\gamma = \mu + M$. Then
$\af(t_{\ov{\Lambda_{j}}} w_i)(\mu_i)
= \ov{\Lambda_j} + w_i(\mu_i)$
so that $\sigma_i(\gamma_i) = s_i + w_i(\gamma_i)$. Since the Weyl group stabilises the cominimal simple currents, we obtain
\[  \sigma_i(\gamma_i) = w_i(s_i+\gamma_i)   \,  . \]
We can now compute
\begin{align*}
	F_{\gamma+s}(\tau)
	& = \sum_{\lambda\in P_V}m_\lambda\prod_{i=1}^rc^{\lambda_i}_{s_i+\gamma_i}(\tau) 
	= \sum_{\lambda\in P_V}m_\lambda\prod_{i=1}^rc^{\lambda_i}_{w_i(s_i+\gamma_i)}(\tau)\\
	& = \sum_{\lambda\in P_V}m_{\sigma(\lambda)}\prod_{i=1}^rc^{\sigma_i(\lambda_i)}_{\sigma_i(\gamma_i)}(\tau) 
        = \sum_{\lambda\in P_V}m_{\lambda}\prod_{i=1}^rc^{\lambda_i}_{\gamma_i}(\tau) =F_\gamma(\tau) .
\end{align*}
The second equality uses the invariance of the string functions under the Weyl group, the third is a resummation and the fourth follows from the invariance of the string functions under cominimal simple currents and Proposition \ref{schellekens}. This is the fourth claim.

Let $\prod_{i=1}^rc^{\lambda_i}_{\gamma_i}$ be the summand with the highest pole order in the expression for $F_{\gamma}$. For each $i$ let $\mu_i \in (\pi^{k_i})^{-1}(\gamma_i)$ be a maximal weight. Then
\[ c^{\lambda_i}_{\gamma_i} (\tau) = c_{\mu_i}^{\lambda_i}(\tau)
= q^{m_{\lambda_i, \mu_i}} \sum_{n \in \Z_{\geq 0}} \mult_{L_{\g_i,k_i}(\lambda_i)}(\mu_i-n\delta) q^n \]
so that the pole order of $\prod_{i=1}^rc^{\lambda_i}_{\gamma_i}$ is bounded by
\begin{multline*}
  	- \sum_{i=1}^rm_{\lambda_i,\mu_i}
	\leq \sum_{i=1}^r \frac{\rho_i^2}{2h^{\vee}_i} \frac{k_i}{k_i+ h^{\vee}_i} 
	= \sum_{i=1}^r \frac{\dim(\g_i)}{24} \frac{k_i}{k_i+ h^{\vee}_i} \\ 
	= \sum_{i=1}^r \frac{\dim(\g_i)}{24} \frac{24}{\dim(\g )} = 1 \, .
\end{multline*}
Here we first used Proposition 13.11 in \cite{Kac}, then the strange formula of Freuden\-thal-de Vries and finally Schellekens' equation.

Applying again Proposition 13.11 in \cite{Kac} we see that equality holds if and only if for each $i$ we have $\lambda_i = k_i \Lambda_{j}\! \mod \C\delta$ for some $j$ with $a_j=1$ and $\mu_i=w_i(\lambda_i)$ for some $w_i \in \hat{W}_i$. The last equation implies $\gamma_i= w_i(s_i)$ for a simple current $s_i$ of $\g_i$. Since the Weyl group stabilizes the simple currents, we have $\gamma_i= s_i$. It follows that if $F_\gamma$ has pole order $1$, then $\gamma \in H_V$. Conversely by ii) and iv) the component $F_\gamma$ has pole order $1$ for $\gamma \in H_V$.
\eop

\medskip

We can rewrite the character as
\begin{align*}
	\chi_V
	& = \sum_{\gamma \in M'/M} F^M_{\gamma} \theta^M_{\gamma}  
	= \sum_{\gamma \in H_V^{\perp}} F^M_{\gamma} \theta^M_{\gamma}  \\
	& = \sum_{\mu \in H_V^{\perp}/H_V} \sum_{\gamma \in H_V} F^M_{\mu + \gamma} \theta^M_{\mu + \gamma} 
	= \sum_{\mu \in K'/K} F_{\mu} \theta_{\mu} 
\end{align*}
where 
\[  K  = \bigcup_{\gamma \in H_V} (\gamma + M) \subset M'  \]
is the {\em lattice associated} with $V$ (for a different approach to this lattice cf.\ \cite{M14} and \cite{H17}). It is a Jacobi form of lattice index $K$ and weight $0$ with Fourier expansion
\[  \chi_V(v,\tau) = \sum_{\substack{\al\in K'\\ \!\! n \in \Z}} [F_\al](n-\al^2/2)e(\la \al,v \ra)q^{n}.  \]

\medskip

In order to show that $F$ is reflective, we construct a Lie algebra $\g(V)$ corresponding to $V$.  
The $b,c$-ghost system of the bosonic string can be described by the vertex operator superalgebra $V_{\Z\sigma}$ with $\sigma^2 =1$. The tensor product
\[   V \otimes V_{I\!I_{1,1}} \otimes V_{\Z\sigma}    \]
is acted on by the BRST-operator $Q$ with $Q^2 =0$. The cohomology group of ghost number $1$ is a Lie algebra \cite{LZ} which we denote by $\g(V)$. It is graded by $K' \oplus I\!I_{1,1}$, i.e.
\begin{gather*}
	\g(V) = \bigoplus_{\al \in K' \oplus I\!I_{1,1}  }  \g(V)_{\al}\, ,   \\
	[ \g(V)_{\al} , \g(V)_{\bt} ] \subset \g(V)_{\al+\bt}\,. 
\end{gather*}
The Lie algebra $\g(V)$ carries a non-degenerate, invariant symmetric bilinear form $( \, , \, )$ satisfying
\[   ( \g(V)_{\al} , \g(V)_{\bt} )  = 0   \]
if $\al + \bt \neq 0$. We define $\h(V)$ as the subspace of degree $\al =0$. Then $\h(V)$ is isometric to $(K' \oplus I\!I_{1,1})\otimes_{\Z}\C$ under $( \, , \, )$ and
\[  \g(V)_{\al} = \{ x \in \g(V) \, | \, [h,x] = (h,\al) x \, \text{ for all } h \in \h(V)  \}   \, . \]
The No-ghost theorem implies that the multiplicity of $\al \in (K' \oplus I\!I_{1,1})\setminus \{ 0 \}$ is given by
\[  \mult (\al) = \dim (\g(V)_{\al}) = [ F_{\al} ] (-\al^2/2)\,.   \]
We call $\al \in (K' \oplus I\!I_{1,1})\setminus \{ 0 \}$ a \textit{root} if $\mult (\al) \neq 0$ and remark that the roots generate $K' \oplus I\!I_{1,1}$ (and not just a sublattice). The Lie algebra $\g(V)$ is almost a Kac-Moody algebra but not quite.
Suppose $V$ is unitary. Then Theorem $1$ in \cite{B95} (cf.\ also Lemma 3.4.2 in \cite{C12}) can be used to show

\begin{thm}
The Lie algebra $\g(V)$ is a generalised Kac-Moody algebra.
\end{thm}

Since the real roots (the roots of positive norm) of $\g(V)$ have multiplicity $1$, the coefficients of the principal part of $F$ are $0$ or $1$. We consider the theta lift $\psi_F$ of $F$ on $L = K \oplus I\!I_{1,1} \oplus I\!I_{1,1}$. 

\begin{thm}\label{everythingcounts}
	The function $\psi_F$ is a reflective automorphic product of singular weight.
\end{thm}
{\em Proof:}
It suffices to show that $\psi_F$ is geometrically reflective. Let $\lambda \in L$ be primitive and of norm $\lambda^2 > 0$. Suppose that the divisor $\lambda^{\perp}$ has positive order. Define a positive integer $m$ by $(\lambda,L) = m\Z$. Then $\gamma = \lambda/m$ is primitive in $L'$. We choose $\mu \in K'\oplus I\!I_{1,1}$ primitive such that $\gamma^2 = \mu^2$ and $\mu = \gamma \, \bmod L$. The order of $\lambda^{\perp}$ is given by
\[  \sum_{n=1}^{\infty} [ F_{n\gamma} ](-n^2\gamma^2/2) =
\sum_{n=1}^{\infty} [ F_{n \mu} ](-n^2\mu^2/2)   \, . \]
Let $k$ be a positive integer such that $[ F_{k\mu} ](-k^2\mu^2/2)$ is non-zero. Then $\al = k \mu$ is a real root of $\g(V)$. Hence
\[  [ F_{k\mu} ](-k^2\mu^2/2) = [ F_{\al} ](-\al^2/2) = \mult(\al) = 1  \]
and $k\mu$ is the only positive integral multiple of $\mu$ that is a root of $\g(V)$. It follows that $\lambda^{\perp}$ has order $1$. Furthermore the reflection corresponding to $\al$ preserves the root lattice $K'\oplus I\!I_{1,1}$ of $\g(V)$. This implies that $\mu$ is a root of $K'\oplus I\!I_{1,1}$ so that $2\mu/\mu^2 \in (K' \oplus I\!I_{1,1})' = K \oplus I\!I_{1,1}$ and $2/\mu^2 \in \Z$. Write $\gamma = \mu + x$ with $x \in L$. Then
\[  2\gamma/\gamma^2 = 2(\mu + x)/\mu^2 = 2\mu/\mu^2 + (2/\mu^2)x \in L  \, . \]  
Hence $\gamma$ is a root of $L'$. From this we derive that $\lambda = m \gamma$ is a root of $L$. \eop

\medskip

The modular variety $\text{O}(L,F)^+\backslash {\cal H}$ has a unique $0$-di\-men\-sional cusp ${\cal C}$ of type $0$ with associated lattice $K \oplus I\!I_{1,1}$. At this cusp the expansion of $\psi_F$ is given by
\[  e((\rho,z_L)) \prod_{\substack{\al \in K' \oplus I\!I_{1,1}\\[0.03mm] (\al,C) > 0}}\big(1-e(-(\al,z_L))\big)^{[F_{\al}](-\al^2/2)} \, . \]
This is the denominator function of the generalised Kac-Moody algebra $\g(V)$.

We can recover the affine structure of $V$ as follows.

\begin{thm}  \label{aw}
	The decomposition $L = K \oplus I\!I_{1,1} \oplus I\!I_{1,1}$ defines a $1$-dimensional cusp ${\cal C}$ of $\text{O}(L,F)^+\backslash {\cal H}$ of type $0$ with associated lattice $K$. The scaled root system 
	\[  R_{\cal C} = \{ \al \, | \, \al \in K'\backslash\{0\}, \, [F_{\al}](-\al^2/2) = 1 \}  \]
	associated with ${\cal C}$ is the root system of the affine structure of $V$ together with its scaling.
\end{thm}
{\em Proof:} The Cartan subalgebra $\h \subset V_1$ acts semisimply on the spaces $V_n$. The Fourier expansion of $\chi_V$ given above shows that the weight space of degree $\al \in K'$ in $V_n$ has dimension $[F_{\al}](n-1-\al^2/2)$. The action of $\h$ on $V_1$ is the restriction of the adjoint representation of $V_1$ and its weights are precisely the roots of $V_1$. Hence the claim of the theorem follows by taking $n=1$. \eop

\medskip

\noindent
In particular the first non-vanishing coefficient of the expansion of $\psi_F$ at ${\cal C}$ is the denominator function of $\hat{\g}_{1} \oplus \ldots \oplus \hat{\g}_{r}$.

\medskip

Combining the results of this section with our classification result we obtain

\begin{thm}
Let $V$ be a holomorphic vertex operator algebra of central charge $24$ with non-trivial, semisimple weight-$1$ space. Suppose $V$ is unitary and the lattice associated with $V$ is regular and of even rank. Then the affine structure of $V$ is one of the $69$ non-trivial structures given in Theorem \ref{heaven17}.
\end{thm}

We relate our results to the current state of research in vertex operator algebra theory. The weight-$1$ subspace $V_1$ of a holomorphic vertex operator algebra $V$ of central charge $24$ is a reductive Lie algebra. This Lie algebra is either trivial or abelian and of rank $24$ or non-trivial and semisimple. In the third case the isomorphism type of $V_1$ is called the affine structure of $V$. Schellekens \cite{ANS} showed, using the theory of Jacobi forms and extensive computer calculations, that there are at most $69$ possibilities for this structure (cf.\ also \cite{EMS}). He conjectured that each of these Lie algebras is realised by a unique vertex operator algebra (\cite{ANS}, Section $1$). He also  asked whether the list has a natural substructure (\cite{ANS}, Section $5$). In 2017 H\"ohn \cite{H17} observed that the $69$ possible affine structures can be related to the $11$ classes in Conway's group $\text{Co}_0$ described in Section \ref{DM}. He showed that the simple current extensions of the vertex operator algebras $V_{\Lambda^g} \otimes (V_{\Lambda^{g \perp}})^{\hat g}$ where $g$ ranges over the aforementioned classes realise the affine structures in Schellekens' list. The assumptions made in \cite{H17} where proved by Lam in \cite{L19}. Shortly after H\"ohn's paper Schellekens' conjecture was confirmed. The proof combined the efforts of many authors and was based on a case-by-case analysis. A discrete-geometric proof of the classification independent of Schellekens results was obtained in \cite{MS23, MS24}.
In the present paper we give a natural explanation of the $11$ classes described by H\"ohn in terms of automorphic forms and a complex-geometric derivation of Schellekens' list, both under the stated conditions. It is surprising to us that the only reflective automorphic products of singular weight are those coming from holomorphic vertex operator algebras of central charge $24$ and that Schellekens' list accounts for all type-$0$ cusps of the corresponding reflective modular varieties.

\section*{Appendix}

In the tables below we state more precisely the bounds given in Proposition \ref{sepultura} 
and list the cusp forms used in the proof of Theorem \ref{classref}.
\[ 
\begin{array}{c|c|c|c|c|c|c}
n & \text{level} & N_R & \text{bound} & N_E & \text{cusp forms} & \text{candidate} \\ [0.3mm] \hline
& & & & & & \\ [-3.5mm]
26 & 1 & 1 & 12 & 1 & - & I\!I_{26,2} \\ [0.3mm] \hline
& & & & & & \\ [-3.5mm]
18 & 1 & 1 & 132 & - & - & - \\ [0.2mm]  
  & 2 & 8 & 33 & 4 & {T_2 \eta_{1^8 2^4 4^8}} & I\!I_{18,2}(2_{I\!I}^{+10}) \\ [0.3mm] \hline
  & & & & & & \\ [-3.5mm]
14 & 2 & 5 & 64 & - & - & - \\ [0.2mm] 
& 3 & 5 & 27 & 2 & \eta_{1^6 3^6} \theta_{A_2}^2 & I\!I_{14,2}(3^{-8}) \\ [0.2mm] 
& 4 & 35 & 33 & 20 & {\eta_{1^8 2^8}} & I\!I_{14,2}(2_{I\!I}^{-10} 4_{I\!I}^{-2}) \\ [0.3mm] \hline
& & & & & & \\ [-3.5mm]
12 & 3 & 5 & 53 & 1 & {\eta_{1^6 3^6} \theta_{A_2}} & - \\ [0.2mm] 
& 4 & 10 & 22 & 6 & T_n\eta_{1^4 2^2 4^4} \theta_{A_1}^4, \, n = 1,2 & I\!I_{12,2}(2_{2}^{+2} 4_{I\!I}^{+6}) \\ [0.5mm] \hline
& & & & & & \\ [-3.5mm]
10 & 1 & 1 & 252 & - & - & - \\ [0.2mm] 
& 2 & 4 & 128 & - & - & - \\ [0.2mm] 
& 3 & 4 & 30 & 1 & {\eta_{1^6 3^6}} & - \\ [0.2mm] 
& 4 & 22 & 66 & 6 & {\eta_{2^{12}}} & - \\ [0.2mm] 
& 5 & 8 & 47/2 & 3 & \eta_{1^3 5^9}, \, \eta_{1^4 5^4} \theta_K & I\!I_{10,2}(5^{+6}) \\ [0.2mm] 
& 6 & 25 & 49 & 12 &  \eta_{1^6 3^6}, \, \eta_{2^6 6^6}, \, \eta_{1^1 2^1 3^5 6^5} & I\!I_{10,2}(2_{I\!I}^{+6} 3^{-6}) \\ [0.2mm] 
& 9 & 36 & 57/2 & 6 & {\eta_{1^6 3^6}} & - \\ [0.3mm] \hline
& & & & & & \\ [-3.5mm]
8 & 3 & 3 & 120 & - & - & - \\ [0.2mm] 
& 4 & 3 & 32 & - & - & - \\ [0.2mm] 
& 6 & 15 & 135/2 & 1 & \eta_{1^2 2^2 3^2 6^2} \theta_{A_2} & - \\ [0.2mm] 
& 7 & 3 & 75 & 1 & - & I\!I_{8,2}(7^{-5}) \\ [0.2mm]
& 8 & 18 & 16 & 8 & \eta_{1^4 2^2 4^4}, \, \eta_{2^4 4^4} \theta_{A_1}^2, & I\!I_{8,2}(2_{7}^{+1} 4_{7}^{+1} 8_{I\!I}^{+4}) \\ [0.2mm]
&   &    &    &   & T_n T_2^2 \eta_{1^2 2^4 8^2}, \, n = 1,2,3 & \\ [0.2mm] 
& 9 & 15 & 41 & 1 & T_3 \eta_{3^3 9^7} & - \\ [0.2mm] 
& 12 & 66 & 59/2 & 29 & \text{see text} & I\!I_{8,2}(2_{I\!I}^{+4} 4_{I\!I}^{-2} 3^{+5}) \\ [0.3mm] \hline
& & & & & & \\ [-3.5mm]
6 & 2 & 1 & 96 & - & - & - \\ [0.2mm] 
& 3 & 1 & 48 & - & - & - \\ [0.2mm] 
& 4 & 1 & 42 & - & - & - \\ [0.2mm] 
& 5 & 3 & 62 & - & - & - \\ [0.2mm] 
& 6 & 4 & 67/2 & - & - & - \\ [0.2mm] 
& 7 & 1 & 35/2 & - & - & - \\ [0.2mm] 
& 8 & 4 & 18 & - & - & - \\ [0.2mm] 
& 9 & 3 & 21 & - & - & - \\ [0.2mm] 
& 10 & 10 & 30 & - & - & - \\ [0.2mm] 
& 11 & 1 & 11 & - & - & - \\ [0.2mm] 
& 12 & 11 & 17 & 1 & {\eta_{1^2 2^2 3^2 6^2}} & - \\ [0.2mm] 
& 14 & 4 & 23/2 & 1 & {\eta_{1^2 2^2 7^2 14^2}} & - \\ [0.2mm] 
& 15 & 10 & 23 & 2 & \eta_{1^2 3^2 5^2 15^2}, \, \eta_{1^1 3^6 5^1} & - \\ [0.2mm] 
& 18 & 16 & 33 & 3 & \eta_3^8, \, \eta_{1^2 2^2 3^2 6^2} & - \\ [0.2mm] 
& 20 & 19 & 33/2 & 5 & \eta_{1^4 5^4}, \, \eta_{2^7 10^1}, \, \eta_{2^4 10^4} & I\!I_{6,2}(2_{I\!I}^{-2} 4_{I\!I}^{-2} 5^{+4}) \\ [0.2mm] 
& 25 & 9 & 6 & 1 & \eta_{1^4 5^4} & - \\ [0.2mm] 
& 36 & 37 & 25/2 & 18 & \text{see text} & - 
\end{array}
\]

\[ \renewcommand{\arraystretch}{1.1}
\begin{array}{c|c|c|cp{1cm}c|c|c|c} 
n & \text{level} & N_R & \text{bound} & & n & \text{level} & N_R & \text{bound} \\ \cline{1-4} \cline{6-9}
& & & & & & & & \\ [-3.5mm]
4 & 3 & 1 & 54 & & & 28 & 1 & 11/2 \\ 
  & 6 & 1 & 30 & & & 30 & 3 & 9 \\ 
& 7 & 1 & - & & & 33 & 2 & 7/2 \\ 
& 9 & 1 & 22 & & & 35 & 3 & 4 \\ 
& 11 & 1 & 41/2 & & & 36 & 1 & 11/2 \\ 
& 12 & 1 & 27/2 & & & 42 & 2 & 5 \\ 
& 14 & 1 & 13 & & & 44 & 1 & 4 \\ 
& 15 & 3 & 33/2 & & & 45 & 3 & 13/2 \\ 
& 18 & 1 & 11 & & & 49 & 1 & 7/2 \\ 
& 19 & 1 & 33/2 & & & 60 & 3 & 7/2 \\ 
& 21 & 2 & 17/2 & & & 63 & 4 & 13/2 \\ 
& 22 & 1 & 27/2 & & & 72 & 1 & 1 \\ 
& 23 & 1 & - & & & 75 & 5 & 11 \\ 
& 24 & 1 & 4 & & & 98 & 1 & 2 \\  
& 27 & 2 & 15 & & & 121 & 1 & 2  
\end{array}
\]

\vspace*{2mm}
\noindent
We describe the entries in the third line of the first table to explain our notation. In signature $(18,2)$ and level $2$ there are $8$ regular lattices splitting two hyperbolic planes. The bound for the weight of a reflective automorphic product on such a lattice is $33$. Four out of the $8$ lattices satisfy the Eisenstein condition. Pairing with the cusp form $T_2 \eta_{1^8 2^4 4^8}$ of weight $10$ leaves one lattice which might carry a reflective automorphic product of singular weight, the lattice $I\!I_{18,2}(2_{I\!I}^{+10})$, and shows that $M_2 \cap {\cal O}_2$ has cardinality $496$.

In signature $(10,2)$ and level $5$ the lattice $K$ is the unique lattice in the genus $I\!I_{4,0}(5^{-2})$.

In signature $(4,2)$ there are no lattices satisfying the Eisenstein condition so that we omitted the corresponding columns.

\end{document}